\documentclass[twoside,a4paper,reqno]{amsart}

\usepackage{style}

\usepackage{style}
\usepackage{graphicx}

\usepackage{pdfsync}

\usepackage{xcolor}
\definecolor{webgreen}{rgb}{0,.5,0}
\definecolor{webbrown}{rgb}{.6,0,0}
\definecolor{RoyalBlue}{cmyk}{1, 0.50, 0, 0}
\usepackage[colorlinks=true, breaklinks=true, urlcolor=webbrown, linkcolor=RoyalBlue, citecolor=webgreen,backref=page]{hyperref}

\usepackage{pgfplots}
\pgfplotsset{every tick label/.append style={font=\tiny}, compat=1.17}
\usepackage{mathrsfs}
\usetikzlibrary{arrows}
\usepackage{courier}

\begin{document}

\title[Szeg\H{o} condition, scattering, \ldots]{Szeg\H{o} condition, scattering, and vibration of Krein strings}

\author[R. Bessonov, S. Denisov]{R. Bessonov, S. Denisov}

\address{
\begin{flushleft}
bessonov@pdmi.ras.ru\\\vspace{0.1cm}
St.\,Petersburg State University\\  
Universitetskaya nab. 7-9, 199034 St.\,Petersburg, RUSSIA\\
\vspace{0.1cm}
St.\,Petersburg Department of Steklov Mathematical Institute\\ Russian Academy of Sciences\\
Fontanka 27, 191023 St.Petersburg,  RUSSIA\\
\vspace{0.1cm}
denissov@wisc.edu\\
\vspace{0.1cm}
Department of Mathematics, University of Wisconsin-Madison,\\ 480 Lincoln Dr., Madison, WI 53706, USA
\end{flushleft}
}

\thanks{
Research of RB in Sections 2.1-2.4, 4 is supported by the Russian Science Foundation grant 19-71-30002. Research of SD in Sections 2.5-2.7, 3 is  supported by  NSF DMS-1764245, NSF  DMS-2054465, and Van Vleck Professorship Research Award. RB is a Young Russian
Mathematics award winner and would like to thank its sponsors and jury.
}

\begin{abstract}
We give a dynamical characterization of  measures on the real line with finite logarithmic integral. The general case is considered in the setting of evolution groups generated by de Branges canonical systems. Obtained results are applied to the Dirac operators and Krein strings. 
\end{abstract} \vspace{1cm}

\subjclass[2010]{42C05}
\keywords{Szeg\H{o} class, canonical systems, Krein strings, Dirac operators, evolution, wave operators, scattering}

\maketitle

\setcounter{tocdepth}{2}

\tableofcontents

\newpage 

\section{Introduction}
\subsection{Nonstationary scattering}  \label{sect11}
In nonstationary scattering theory \cite{LPbook}, \cite{RSbook3}, \cite{Yafbook}, one  studies an evolution group of unitary operators $\{U_t\}_{t\in \R}$ that act on Hilbert space $H$. Elements $X \in H$ are called states and  the set $\{U_t X\}_{t\in \R}$ is a trajectory of $X$ under the evolution $U_t$. Motivated by physical heuristics, one sometimes expects that the long-time  behavior of $U_t X$, $t\to \pm \infty$, is asymptotically close to that of $U_t^0 Y_{\pm}$, $t\to \pm \infty$, for certain states $Y_{\pm} \in H$, which depend on $X$. Here,  $\{U^0_t\}_{t\in \R}$ is another evolution group of unitary operators that act on the same Hilbert space $H$. The scattering operator $S$ defined by $S\colon Y_{-} \mapsto Y_{+}$ then recovers the remote future of the process $U_t X$ from its past.

\medskip

For concrete evolution groups $U_t$ and $U^0_t$, proving the existence of the scattering operator might be a nontrivial problem. In this paper, we address this question in the context of $2\times 2$ canonical Hamiltonian systems on the positive half-axis $\R_+ = [0,+\infty)$. They provide a convenient framework for a unifying treatment of classical equations/operators of mathematical physics, such as one-dimensional Dirac systems, Krein strings, Jacobi matrices, and Schr\"odinger operators. In fact, any self-adjoint operator with a simple spectrum can be realized as a canonical Hamiltonian system on $\R_+$ although such a realization is often somewhat implicit. 
We refer to \cite{Remlingb} and \cite{Romanov} for an introduction to the spectral theory of canonical Hamiltonian systems.

\medskip 

A canonical Hamiltonian system is defined  by its Hamiltonian, which is a $2\times 2$ matrix-valued function on $\R_+$ that has the form
$$
\Hh = \begin{pmatrix}h_1 & h\\h&h_2\end{pmatrix}, \qquad \Hh(\tau) \ge 0, \quad \trace\Hh(\tau) >0, \quad \mbox{for a.e.\ }\tau \in\R_+.
$$
The real-valued functions $h_1$, $h_2$, $h$ satisfy $h_{1}, h_{2}, h \in L_{\loc}^1(\R_+)$. If $h = 0$ almost everywhere on $\R_+$, we say that $\Hh$ is diagonal. Let $J = \jm$. With each Hamiltonian $\Hh$ one can associate the canonical system,
\begin{equation}\label{cs-def}
JX'(\tau, z) = z\Hh(\tau) X(\tau, z).
\end{equation}
Here, $z \in \C$ is a spectral parameter and the derivative is taken in $\tau \in \R_+$. Canonical Hamiltonian system \eqref{cs-def} can be considered as a  generalized eigenvalue problem $\Di_{\Hh} X = z X$ for a self-adjoint differential operator 
$$
\Di_\Hh\colon X \mapsto Y, \qquad Y\colon \; JX'=\Hh Y,
$$
densely defined on a certain Hilbert space of functions on $\R_+$ (we give more details in Section \ref{cs-basics}). The self-adjoint operator $\Di_{\Hh}$ has a simple spectrum. It, therefore, admits a spectral representation as the multiplication operator by the independent variable in $L^2(\mu)$, with a canonical choice of the scalar non-negative spectral measure $\mu$ satisfying
\begin{equation}\label{pf}
\int_{\R}\frac{d\mu(x)}{1+x^2} < \infty.
\end{equation}
Remarkably, any non-negative Borel measure $\mu$ satisfying \eqref{pf} is a spectral measure of $\Di_{\Hh}$ for some Hamiltonian $\Hh$. In particular, the usual Lebesgue measure on $\R$ is the spectral measure of $\Di_{\Hh_0}$ for the constant Hamiltonian $\Hh_0 = \idm$ on $\R_+$.

\medskip

In our paper, we study dynamics of the unitary evolution group $U_t = e^{it\Di_{\Hh}}$ generated by a general  Hamiltonian $\Hh$ by comparing it to the ``unperturbed'' dynamics governed by $U^0_t = e^{it\Di_{\Hh_0}}$. The latter can be easily reduced to the shift operator on $L^2(\R)$, see Section \ref{app2}. Informally, one of our central results can be summarized as follows: ``scattering for the pair $U_t$, $U_t^0$ takes place if and only if the spectral measure $\mu$ of $\Di_{\Hh}$ belongs to the  Szeg\H{o} class $\sz$''. That class consists of Borel non-negative measures on the real line $\R$ whose density with respect to the Lebesgue measure on $\R$ has a finite logarithmic integral:
\begin{equation}\label{sz}
\sz = \left\{\mu = w\,dx + \mus :  \;\; \int_{\R}\frac{d\mu(x)}{x^2 + 1} < +\infty, \; \int_{\R}\frac{\log w(x)}{x^2+1}\,dx > -\infty\right\}. 
\end{equation}  
The Szeg\H{o} class is prominent in complex analysis  \cite{KoosisLI2}, \cite{KoosisLI}, theory of stationary processes  \cite{DymMcKean}, \cite{ibr}, orthogonal polynomials \cite{Simonbook}, \cite{Simonbook2}, \cite{SzegoBook}, and statistical physics \cite{DIK13}, \cite{SimonDes}. We will discuss its appearance in various aspects of scattering for canonical Hamiltonian systems: propagation of a single wavepacket (Theorem~\ref{t2} and Theorem~\ref{t1}), existence and completeness of wave operators (Theorem~\ref{t3}), dynamical classification of spectral types (Theorem~\ref{t4} and  Theorem~\ref{c3}). We also use our previous work to provide an explicit description of Hamiltonians corresponding to spectral measures in $\sz$ (Proposition~\ref{str1}).    

\medskip

Below we discuss how our general theory, summarized in Section \ref{css}, applies to two important classes of operators: Dirac systems and Krein strings. 
In these two cases, the connection between Szeg\H{o} condition on the spectral measure, propagation of the wavepacket, and the existence of the wave operators becomes particularly transparent. The Dirac systems give rise to locally absolutely continuous Hamiltonians $\Hh$ with $\det \Hh = 1$ on $\R_+$, and Krein strings are in one-to-one correspondence with diagonal Hamiltonians on $\R_+$. Some figures demonstrating a numerical simulation of the propagation of waves can be found in Sections \ref{s1p4} and \ref{s1p5}.

\medskip

We believe that suitable modifications of Szeg\H{o} class will appear naturally in other criteria for scattering in different settings (say, for CMV matrices, Jacobi matrices, Dirac operators with the positive mass, etc.) if one chooses  the unperturbed dynamics properly (see, e.g., \cite{Den_B1},\cite{VolYu2002}).

\medskip

\subsection{Dirac equation}\label{s1point2} Define $\Psi$ as the solution to the following Cauchy problem for one-dimensional Dirac equation on the positive half-line $\R_+ = [0, +\infty)$: 
\begin{equation}\label{cp-i}
J\Psi'(\tau, z) + Q(\tau)\Psi(\tau, z) = z \Psi(\tau,z), \qquad \Psi(0, z) = \oz, \qquad \tau \in \R_+, \quad z \in \C.
\end{equation} 
Here, again, $J = \jm$ and the derivative is taken with respect to $\tau$. The potential $Q$ is $2\times 2$ matrix-function with  real entries. It is symmetric, has zero trace, and satisfies $Q\in L^1_{\loc}(\R_+)$. We will write $Q$ in the form 
\begin{equation}\label{potential}
Q = \begin{pmatrix}q_1 & q_2\\ q_2  & -q_1\end{pmatrix},
\end{equation}
where real-valued functions $q_1$ and $q_2$ on $\R_+$ satisfy $q_{1}, q_2 \in L^1_{\loc}(\R_+)$. For each value of the spectral parameter $z \in \C$, the solution, $\Psi(\cdot,z)$, is a locally absolutely continuous function on $\R_+$ with values in $\C^2$. It can be considered as the generalized eigenvector of the Dirac operator 
\begin{equation}\label{do-i}
\Dd_Q \colon Z \mapsto JZ' + Q Z. 
\end{equation}
The operator $\Dd_{Q}$ is a densely defined self-adjoint operator on the Hilbert space 
$$L^2(\R_+,\C^2) = \Bigl\{Z\colon \R_+ \to \C^2:  \|Z\|_{L^2(\R_+,\C^2)}^{2} = \int_{\R_+}\|Z(\tau)\|_{\C^2}^2\,d\tau < \infty \Bigr\}.
$$
Its domain consists of locally absolutely continuous functions $Z \in L^2(\R_+,\C^2)$ that satisfy two conditions: $JZ' + QZ \in L^2(\R_+,\C^2)$ and $\langle Z(0), \zo \rangle_{\C^2} = 0$.  The ``free'' Dirac operator, corresponding to the potential $Q = 0$, will be denoted by $\Dd_0$. 

\medskip

The one-dimensional Dirac operator \eqref{do-i} has its origin in  relativistic quantum mechanics. It appears after the separation of variables as the radial part of a ``full'' Dirac operator that describes a relativistic particle in a radially symmetric external field in $\R^3$ (see Section 4.6 in \cite{Thaller}). In the one-dimensional model, $\Dd_{Q}$ describes the particle of unit mass moving in a field defined by a potential $q\colon \R_+ \to \R$. That function $q$  is related to  $q_1$ and $q_2$ in \eqref{potential} by $q_1=\cos(2\int_0^\tau q\,ds)$ and $q_2=-\sin(2\int_0^\tau q\,ds)$, see page~534 in \cite{HSch14}, and a discussion at the end of Section \ref{wfp}. In the theory of completely integrable systems, the Dirac equation appears as a linear self-adjoint problem that is used to solve 
the nonlinear Schr\"odinger equation. See, e.g., \cite{FTbook}, Chapter~1 for the inverse scattering approach in the theory of nonlinear Schr\"odinger equation, and \cite{Marchenko}, \cite{deift}, \cite{kupin} for more on inverse scattering problems.   

\medskip

For every $Q \in L^1_{\loc}(\R_+)$, there is a unique Borel measure $\mu_D$ on $\R$ such that $(x^2+1)^{-1} \in L^1(\mu_D)$ such that the generalized Fourier transform
\begin{equation}\label{ft-i}
\F_Q\colon Z \mapsto \frac{1}{\sqrt{\pi}}\int_{\R_+}\langle Z(\tau), \Psi(\tau, \bar z)\rangle_{\C^2}\,d\tau, \qquad z \in \C,
\end{equation} 
densely defined on smooth functions with compact support, can be extended to a unitary map from $L^2(\R_+,\C^2)$ onto $L^2(\mu_D)$. That measure is called the {\bf main spectral measure} of $\Dd_{Q}$. For example, the Lebesgue measure on $\R$ is the main spectral measure for the free Dirac operator $\Dd_0$. We refer the reader interested in the spectral theory of the one-dimensional Dirac operator to the classical monograph \cite{LSb}. It can also be considered as a part of a more general spectral theory of canonical Hamiltonian systems \cite{dbbook},\cite{Potapov},  \cite{Remlingb}, \cite{Romanov}.

\medskip

Given a potential $Q \in L^1_{\loc}(\R_+)$, we let  $U_t = e^{it\Dd_Q}$ and $U^0_t = e^{it\Dd_{0}}$ in the context of general problem discussed in Section \ref{sect11}. These are unitary operators on $L^2(\R_+,\C^2)$ parameterized by $t \in \R$. The unperturbed dynamics $U_t^0$ has an explicit form given by formula \eqref{eq32bisbis} below. 
Fix a measurable time-independent ``phase function'' $\gamma$ on $\R_+$ which takes its values on the unit circle $\T = \{z \in \C: |z| = 1\}$. Let 
$M_\gamma\colon Z \mapsto \gamma Z$ denote the multiplication operator on $L^2(\R_+,\C^2)$ with function $\gamma$. Set $U^{0}_{\gamma, t} = M_\gamma U^{0}_{t}$ for $t \ge 0$ and $U^{0}_{\gamma, t} = M_{\overline{\gamma}} U^{0}_{t}$ for $t < 0$.

\medskip

\noindent {\bf Definition.} If the following limits
\begin{equation}\label{wo-}
W_\pm(\Dd_Q, \Dd_{0},\gamma) =  \lim_{t \to \pm\infty} U_{-t} U_{\gamma,t}^{0} 
\end{equation}
exist in the strong operator topology, we will call them the modified {\bf wave operators} for evolution groups $U_t$ and $U^0_t$.  The standard M\"{o}ller wave operators correspond to the choice $\gamma = 1$ on $\R_+$. \medskip

Studying wave operators is a classical problem in scattering theory, theory of  partial differential equations, and mathematical physics (see, e.g., H\"{o}rmander \cite{Hormander}, Kato \cite{Kato}, Birman and M. Krein \cite{BirmanKrein}, Lax and Phillips \cite{LPbook}, Yafaev \cite{Yafbook}, etc.). In our paper, we show that wave operators for $\Dd_Q$, $\Dd_{0}$ exist if and only if the spectral measure $\mu_D$ belongs to Szeg\H{o} class \eqref{sz}. Specifically, we prove the following two theorems.

\begin{Thm}\label{t4-i}
Let $\Dd_{Q}$ be the Dirac operator \eqref{do-i} with potential $Q \in L^1_{\loc}(\R_+)$. Assume that for some measurable function $\gamma\colon \R_+ \to \T$ 
one of the wave operators $W_{\pm}(\Dd_{Q}, \Dd_{0}, \gamma)$ exists. Then, the main spectral measure $\mu_D$ of $\Dd_{Q}$ belongs to the Szeg\H{o} class $\sz$.
\end{Thm}

The wave operators $W_\pm(\Dd_Q, \Dd_0, \gamma)$ are called {\bf complete} if they are unitary operators from $L^2(\R_+,\C^2)$ onto the absolutely continuous subspace $H_{\bf ac}(\Dd_Q)$ of $\Dd_Q$. 

\begin{Thm}\label{t5-i}
Let $\Dd_{Q}$ be the Dirac operator \eqref{do-i} with potential $Q \in L^1_{\loc}(\R_+)$. Assume that
the spectral measure $\mu_D$ of $\Dd_Q$ belongs to the Szeg\H{o} class $\sz$. Then, the wave operators $W_{\pm}(\Dd_{Q}, \Dd_{0}, \gamma)$ exist and are complete for some measurable function $\gamma\colon \R_+ \to \T$. Moreover, one can take $\gamma = 1$ if $Q$ is anti-diagonal, i.e., $q_1 = 0$.  
\end{Thm}
In general, one cannot take $\gamma = 1$ in Theorem \ref{t5-i} as has been shown in \cite{Den06} (see discussion after Theorem 14.7 there and Teplyaev's work \cite{Tep05}). The wave operators $W_{\pm}(\Dd_{Q}, \Dd_{0}, \gamma)$ can be used to describe asymptotics of evolution in the remote future given  its behavior in the remote past.
In the setting of Theorem~\ref{t5-i}, the scattering operator 
$$
S = W_{+}^{-1}W_{-}, \qquad W_{\pm} = W_{\pm}(\Dd_{Q}, \Dd_{0}, \gamma),
$$
defines a unitary map on $L^2(\R_+,\C^2)$. We prove that it does not depend on $\gamma$. Our wave operators are complete and  $S$ describes the asymptotic dynamics of a state $Z \in H_{\bf ac}(\Dd_Q)$ under the evolution $U_t$ by
\begin{equation}\label{eq1-i}
S\colon Y_{Z,-} \mapsto Y_{Z,+}, \qquad \lim_{t \to \pm \infty}\|U_t Z - U^{0}_{\gamma, t}Y_{Z,\pm}\|_{L^2(\R_+,\C^2)} = 0.
\end{equation}
In other words, if $U_t Z$ is asymptotically close to $U^{0}_{\gamma, t}Y_{Z,-}$ in the remote past for some $Y_{Z,-}$, then $U_t Z$ is asymptotically close to $U^{0}_{\gamma, t}Y_{Z,+}$ in the remote future, hence $Y_{Z,+} = S Y_{Z,-}$. Moreover, both $Y_{Z,+}$ and $Y_{Z,-}$ are in one-to-one correspondence with $Z$ because we have $Y_{Z,\pm} = W^{-1}_{\pm}Z$ and the operators $W_{\pm}$ are unitary from   $L^2(\R_+,\C^2)$ onto $H_{\bf ac}(\Dd_Q)$. Notice that analogs of wave operator and scattering map can also be defined in the context of classical Hamiltonian mechanics (see, e.g., \cite{Buslaev}).

\medskip

Paraphrasing Theorem \ref{t4-i} and Theorem \ref{t5-i}, we can now say that the scattering phenomenon in the sense of \eqref{eq1-i} takes place for the Dirac evolution groups $U_t$, $U_t^{0}$ if and only if the main spectral measure of $\Dd_{Q}$ lies in the Szeg\H{o} class $\sz$. That motivates us to introduce the class of potentials
$$
\szd = \bigl\{Q \in L^1_{\loc}(\R_+)\colon \; \mbox{the main spectral measure of $\Dd_{Q}$ is  in }\sz\bigr\}.
$$
In previous works \cite{BD2017} and \cite{BD2019}, we characterized potentials $Q \in \szd$. That description is summarized in the following theorem.
\begin{Thm}\label{t3-i}
Let $Q$ be a potential in $L^{1}_{\loc}(\R_+)$ and $N_0$ be the solution to the Cauchy problem
$$
JN_0'(\tau) + Q(\tau)N_0(\tau) = 0, \qquad N_0(0) = \idm, \qquad \tau \in \R_+.
$$ 
Then, $Q \in \szd$ if and only if 
$$
\sum_{n \ge 0}\left(\det\int_{n}^{n+2}N_0^{*}(\tau)N_0(\tau)\,d\tau - 4\right) < +\infty.
$$
If, moreover, $Q$ has the form 
\begin{equation}\label{eq60-i}
 Q = \begin{pmatrix}q & 0\\ 0 & -q\end{pmatrix} \quad \mbox{or}\quad Q = \begin{pmatrix}0 & q\\ q & 0\end{pmatrix},
\end{equation}
then $Q \in \szd$  if and only if 
\begin{equation}\label{eq61-i}
\sum_{n \ge 0}\left(\int_{n}^{n+2}h(\tau)\,d\tau\int_{n}^{n+2}\frac{d\tau}{h(\tau)} - 4\right) < +\infty,
\end{equation}
where $h(\tau) = e^{2\int_{0}^{\tau}q(s)\,ds}$, $\tau \ge 0$.
\end{Thm}
Previous works  provide sufficient conditions for scattering in  Dirac evolution and we will discuss some of them now.
We write $Q \in L^p(\R_+)$ if the entries of a potential $Q$ belong to the space $L^p(\R_+)$. 
The case $Q \in L^1(\R_+)$ is classical, the existence and completeness of $W_{\pm}(\Dd_{Q}, \Dd_{0}, 1)$ follow from the general theorems on trace-class perturbations. Indeed, one can show that the operator $(\Dd_{Q} + i)^{-1} - (\Dd_{0} + i)^{-1}$ belongs to the trace class $S^1(L^2(\R_+,\C^2))$ and then use a result of Birman and Krein (see Theorem 2 in \cite{BirmanKrein} or Theorems XI.8, XI.9 in \cite{RSbook3}) which is an extension of the classical Kato-Rosenblum theorem.  In a  general setting, the trace class in Birman-Krein theorem cannot be replaced by any other Schatten class $S^p(L^2(\R_+,\C^2))$, $p > 1$. However,  for the Dirac equation, the existence of wave operators  was proved under  assumptions  much weaker than $(\Dd_{Q} + i)^{-1} - (\Dd_{0}+i)^{-1} \in S^1(L^2(\R_+,\C^2))$. 
For example,  Christ and Kiselev \cite{ChK} showed that the wave operators  $W_{\pm}(\Dd_{Q}, \Dd_{0}, 1)$ exist and are complete for $Q\in L^p(\R_+), 1\le p<2$. The second author covered the borderline case in \cite{SA_GAFA}: he proved that the wave operators  exist and are complete for $Q \in L^2(\R_+)$. On $L^p(\R_+)$-scale, the class $L^2(\R_+)$ is optimal: a well-known result by Pearson \cite{P78}, when stated  for  Dirac equation, says that there exists $Q \in \cap_{p > 2} L^p(\R_+)$ for which the Dirac operator $\Dd_{Q}$ has empty absolutely continuous spectrum. That implies $H_{\bf ac}(\Dd_{Q}) = \{0\}$ and the wave operators do not exist in a ``very strong sense''. Indeed, there are no isometric operators between $L^2(\R_+,\C^2)$ and $\{0\} = H_{\bf ac}(\Dd_{Q})$!

\medskip

Previous results indicated a connection between the convergence of logarithmic integral of the spectral measure $\mu_D$ and the existence of wave operators. The first author proved \cite{B2018} that Szeg\H{o} condition is sufficient for the existence of certain modified wave operators and  one can show that the spectral measure  belongs to the Szeg\H{o} class $\sz$ for $Q \in L^p(\R_+)$, $1 \le p \le 2$. In higher dimensions, similar results were obtained in \cite{Laptev}, \cite{Safronov}, and \cite{Den_M}.
The present paper provides the final answer in the form of a necessary and sufficient condition for Dirac scattering given both in terms of spectral data (via Theorem \ref{t4-i} and Theorem \ref{t5-i}) and an  explicit condition on potential $Q$ (via Theorem \ref{t3-i}). We rely on our previous works \cite{BD2017}, \cite{BD2019}, \cite{B2018}, \cite{Den06},  and, more broadly, on M.~Krein's idea (see \cite{Krein54} and \cite{Krein55}) to use the theory of polynomials orthogonal on the unit circle when studying problems of spectral theory. In particular, Lemma \ref{l2} that appears later in the text has the counterpart known as ``Khrushchev's theorem'' for orthogonal polynomials (see \cite{KH01} by Khrushchev).

\medskip

\begin{figure}
\centering
\definecolor{rvwvcq}{rgb}{0.08235294117647059,0.396078431372549,0.7529411764705882}
\begin{tikzpicture}[line cap=round,line join=round,>=triangle 45,x=2cm,y=1cm]
\begin{axis}[
x=1.6cm,y=1.6cm,
axis lines=middle,
xlabel = $\alpha$,
ylabel = $\beta$,
ymajorgrids=true,
xmajorgrids=true,
xmin=-3.8297866142200694,
xmax=4.876831972332251,
ymin=-1.360897428230408,
ymax=2.7475382172989677,
xtick={-3.5,-3,...,4.5},
ytick={-1,-0.5,...,2.5},]
\clip(-3.8297866142200694,-1.360897428230408) rectangle (4.876831972332251,2.7475382172989677);
\fill[line width=0pt,color=rvwvcq,fill=rvwvcq,fill opacity=0.1] (6,0) -- (6,6) -- (-6,6) -- (-5.977749617746582,-5.477749617746582) -- (0,0.5) -- (1,0.5) -- (6.016007630200937,-4.516007630200937) -- cycle;
\draw [line width=1pt,dash pattern=on 4pt off 4pt] (0,0.5)-- (1,0.5);
\draw [line width=1pt,dash pattern=on 4pt off 4pt,domain=-3.8297866142200694:0] plot(\x,{(-0.25-0.5*\x)/-0.5});
\draw [line width=1pt,dash pattern=on 4pt off 4pt,domain=1:4.876831972332251] plot(\x,{(--0.75-0.5*\x)/0.5});
\begin{scriptsize}
\draw [fill=black] (0,0.5) circle (1.5pt);
\draw [fill=black] (-0.5,0) circle (1.5pt);
\draw [fill=black] (1,0.5) circle (1.5pt);
\draw [fill=black] (1.5,0) circle (1.5pt);
\draw [fill=white] (0,0.5) circle (1pt);
\draw [fill=white] (-0.5,0) circle (1pt);
\draw [fill=white] (1,0.5) circle (1pt);
\draw [fill=white] (1.5,0) circle (1pt);
\end{scriptsize}
\end{axis}
\end{tikzpicture}
  \caption{The set $A$ in Theorem \ref{t6-i} (filled blue). Pairs $(\alpha, \beta) \in A$ correspond to potentials $Q \in \szd$. The set $A$ is open, which demonstrates the stability of the scattering phenomenon in parameters $\alpha$ and $\beta$.} 
  \label{fig:plot-i}
\end{figure}
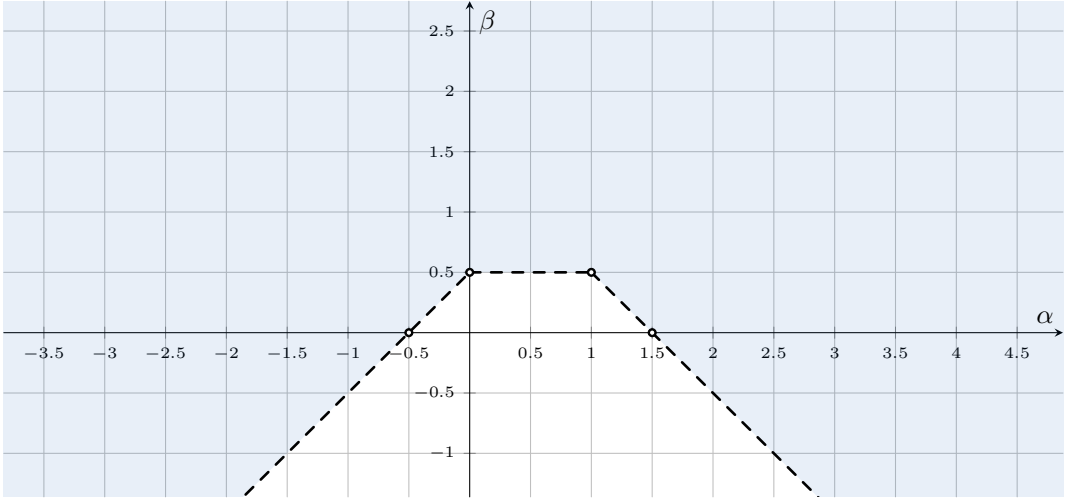

Consider $Q$ of the form \eqref{eq60-i} with $q = \frac{\sin \tau^\alpha}{\tau^\beta}$ for some $\alpha,\beta \in \R$. That class of potentials, introduced by von Neumann and Wigner in a different context, was extensively studied  in the literature (see, e.g., \cite{ChK} for some references). Different values of parameters $\alpha, \beta$ give rise to highly oscillating, slowly oscillating, periodic, decaying or growing potentials. The scattering problem for such potentials has been actively studied (see \cite{BA80}, \cite{BAD79},  \cite{D80},\cite{MS72},  \cite{White83}) in the setting of Schr\"odinger equation.
\begin{Thm}\label{t6-i}
Suppose that a potential $Q \in L^1_{\loc}(\R_+)$ has the form \eqref{eq60-i}, with $q = \frac{\sin \tau^\alpha}{\tau^\beta}$ 
on $[\tau_0, +\infty)$ for some $\tau_0 > 0$ and $\alpha, \beta \in \R$. Then, $Q \in \szd$ if and only $(\alpha, \beta) \in A$, where 
$$
A= \{\alpha \le 0,\; \beta-\alpha > \frac 12\} \cup  \{\alpha \in (0,1),\; \beta > \frac 12\} \cup \{\alpha \ge 1, \;\alpha+\beta > \frac 32\}
$$ 
is the open set depicted on Figure~\ref{fig:plot-i}.
\end{Thm}

\medskip

Most of the results in this paper have a dynamical interpretation. It also appears that some difficult  questions of spectral theory have precise answers in dynamical terms. Here is one example: as the following theorem shows, the propagation of just one nontrivial state $Z$ under the Dirac evolution characterizes the absence of the singular continuous spectrum for the corresponding Dirac operator. 

\begin{Thm}\label{t05-i}
Let $\Dd_{Q}$ be the Dirac operator with potential $Q \in \szd$.  The singular continuous spectrum of $\Dd_Q$ is empty if and only if 
\begin{equation}\label{eq0001-i}
\lim_{b \to +\infty}\lim_{\tb \to +\infty}\frac{1}{\tb}\int_{0}^{\tb}\|U_t Z\|_{L^2(\C^2, [b, t-b])}^{2}\,dt = 0 
\end{equation}
for some (and then for every) compactly supported nonzero state $Z\in L^2(\R_+,\C^2)$.  
\end{Thm}

The next theorem shows that in the Szeg\H{o} case the long-time Dirac evolution of any state $Z \in L^2(\R_+,\C^2)$ decomposes into three parts: ``the bound states part'' localized near the origin, the ``scattering part'' propagating to infinity with constant velocity, and the ``singular continuous part'' between them. These parts correspond to the orthogonal projections $P_{\bf pp}Z$, $P_{\bf ac}Z$, $P_{\bf sc}Z$ of $Z$ onto the pure point, the absolutely continuous, and the singular continuous subspaces of $\Dd_{Q}$. We provide the dynamical description of the sizes of these projections (see  \cite{Damanik}, \cite{Kiselev}, and \cite{Ruelle} for related results) which makes the connection between spectral type and evolution of wavepacket transparent.

\begin{Thm}\label{t7-i}
Let $\Dd_{Q}$ be the Dirac operator with potential $Q \in \szd$. Then, for every $Z \in L^2(\R_+,\C^2)$ we have
\begin{align*}
\lim_{b \to +\infty}\lim_{\tb \to +\infty}\frac{1}{\tb}\int_{0}^{\tb}\|U_{t}Z\|_{L^2(\C^2, [0, b])}^{2}\,dt 
&= \|P_{\bf pp}Z\|_{L^2(\R_+,\C^2)}^{2},\\
\lim_{b \to +\infty}\lim_{\tb \to +\infty}\frac{1}{\tb}\int_{0}^{\tb}\|U_t Z\|_{L^2(\C^2, [b, t-b])}^{2}\,dt &= \|P_{\bf sc} Z\|_{L^2(\R_+,\C^2)}^{2},\\
\lim_{b \to +\infty}\lim_{t \to +\infty}\|U_t Z\|_{L^2(\C^2, [t-b, t+b])}^{2} &= \|P_{\bf ac} Z\|_{L^2(\R_+,\C^2)}^{2},\\
\lim_{b \to +\infty}\lim_{t \to +\infty}\|U_t Z\|_{L^2(\C^2, [t+b, +\infty))} & = 0.
\end{align*} 
\end{Thm}
There is a large amount of literature in which the spectral types and the corresponding subspaces of an operator were studied in connection to the dynamics it generates. That can be done for very general setting (see, e.g., the celebrated RAGE theorem \cite{RSbook3}, Theorem XI.115) or for some standard operators of mathematical physics (see  \cite{DT2010} or \cite{KSL2003}, where the Schr\"odinger evolution was considered). In Theorems \ref{t05-i} and \ref{t7-i}, we focus on Dirac equation and give a complete dynamical description of classical spectral subspaces in the spirit of RAGE theorem. Again, the Szeg\H{o} condition on the spectral measure is central to our analysis as it provides the sharp asymptotics of  $U_tP_{\bf ac}Z$ which simply does not hold in the general case. For example, the Dirac equation with constant positive mass has different dispersion relation (see, e.g., \cite{EG}) and  the third equality in Theorem \ref{t7-i} does not hold for such a model.

\smallskip

We want to make a few remarks about other existing methods. In \cite{PeBook}, p. 406, the  scattering for regular stationary Gaussian sequences has been studied in the case when the spectral measure of the process is  purely absolutely continuous and its density satisfies Szeg\H{o} condition on the unit circle. That is one example of the general Lax-Phillips approach, in which the so-called representation theorem (see, e.g., Chapter~II in \cite{LPbook}) can be applied to the general unitary groups with both discrete and continuous time. It provides the abstract scattering operator under the assumption that the so-called ``outgoing''  and ``incoming'' subspaces exist.   To define these subspaces for concrete evolution equations, one usually works  with compactly supported 
perturbations of the canonical operators (e.g., Laplacian, Dirac, etc.). Our methods give asymptotics of evolution $U_t{Z}$ on the physical side for a large class of ${Z}$ avoiding such strong assumptions on perturbation $Q$ in \eqref{cp-i}.  Also, when viewed on the spectral side, our technique allows the spectral measures to have essentially arbitrary nontrivial singular parts as long as these measures are in Szeg\H{o} class. 

\smallskip

Another approach to scattering theory is based on proving the large-$\tau$ asymptotics of  solutions to \eqref{cp-i} for Lebesgue almost every spectral parameter $z\in \R$. This is an area of active research (see, e.g.,  \cite{ChK},\cite{Mohamed},\cite{Liu}, and more recent work \cite{Polt21}). 

\medskip

\subsection{String equation} To define the mathematical model of a vibrating string, one starts with prescribing its length  $L \in (0, \infty]$ and the non-decreasing right-continuous function $M\colon [0, L) \to \R_+$. Given $\xi\in [0,L)$, the number $M(\xi)$ is interpreted as the mass of the  $[0, \xi]$ piece.  Define the Lebesgue-Stieltjes measure $\mf$ by $\mf[0,\xi]=M(\xi)$ and write its decomposition into the absolutely continuous and singular parts: $\mf=\mf_{\bf ac}+\mf_{\bf s}=\rho\, d\xi+\mf_{\bf s}$. Usually, the function $\rho$ is referred to as the density of the string. Denote $M(L-)=\lim_{\xi\uparrow L}M(\xi)$. We will call the $[M,L]$ pair {\bf proper} if $M$ and $L$ satisfy the following conditions
\begin{align}
&L+M(L-)=\infty, \label{0sa17}\\
&0<M(\xi)<M(L-), \quad \forall\xi\in (0,L).\label{0sa15}
\end{align}
The second condition can be interpreted as the left and the right ends of the string being ``heavy''. The free motion of the vibrating string $[M, L]$ with a given initial displacement $u_0\colon [0,L) \to \R$ is described by the solution $u = u(\xi, t)$ of the string equation
\begin{align}
&\mf(\xi) u_{tt}(\xi, t) = u_{\xi\xi}(\xi, t), \label{seq}\\
&u(\xi, 0) = u_0(\xi), \qquad\qquad\qquad \xi \in [0, L), \quad t \in \R_+{},\label{seq2} \\
&u_t(\xi, 0) = u_\xi(0,t) = 0{.} \label{seq3}
\end{align}
Under mild assumptions \eqref{0sa17}--\eqref{0sa15}, $\mf$ is essentially an arbitrary non-negative Borel measure on $[0, L)$ and one needs to explain how to understand equation \eqref{seq}. 
In Section \ref{ks}, we will define the self-adjoint non-negative operator $\cal{S}_M$, the Krein string operator, which corresponds to the pair $[M,L]$. Then, the spectral theorem for $\cal{S}_M$ and operator calculus can be used  to define  the solution $u$ as follows:
\begin{equation}\label{e44-b}
u(\cdot,t)=\cos (t\sqrt{\cal{S}_M})u_0,  \quad u_0\in L^2(\mf).
\end{equation}  
With this general picture in mind, we mention that if $L=+\infty$, $\mf = \rho \,d\xi$, and $\rho$ and $u_0$ satisfy additional regularity assumptions, then  our solution $u$  coincides with the unique classical solution to the problem
$$
\rho(\xi)u_{tt}(\xi, t) = u_{\xi\xi}(\xi, t),
\qquad \xi \in \R_+, \qquad t \in \R_+,
$$
that satisfies \eqref{seq2} and \eqref{seq3}.  In that case, the value of $u(\xi, t)$ gives the displacement of the string  at the point $\xi \in [0, L)$ at the moment $t \in \R_+$ where $u_0$ is the initial real-valued displacement. Assumption $u_t(\xi, 0) = 0$ indicates that the initial velocity is equal to zero, and the Neumann boundary condition $u_\xi(0,t) = 0$ says that the left end is ``loose'' (see \cite{Belishev1} for some background).

\medskip

In our setup, function $u$, the solution to \eqref{seq}, will be considered as an element of  $L^2(\mf)$ for all $t\ge 0$. For $\mf$-measurable function $f$, we introduce its  {\bf front} as   
\[
\sff[f]=\inf \{a\in \R_+: f=0 \,\, \mf\text{-almost everywhere on\,\,} [a,L)\}
\]
and we will call $\sff_t=\sff[u]$ the {\bf wavefront} of solution $u(\cdot, t)$ at time $t$.
Using a classical result by Krein, one can explicitly compute the wavefront of a wave $u$ with compactly supported initial profile~$u_0$. Given a string $[M,L]$, we  define two functions:
\[
\Ts(\xi)=\int_0^\xi \sqrt{\rho(s)}\,ds, \qquad\Ls(\eta)=\inf\{\xi\in [0,L): \Ts(\xi)\ge \eta\}
\]
for $\xi\in [0,L), \eta \in \R_+$. In physics literature, the former function is sometimes referred to as eikonal or optical metric. The subscript ``$M$'' above refers to the mass distribution $M$ of a string. Later in Section \ref{s22} we use similar functions $\Tc$, $\Lc$ for a canonical Hamiltonian system generated by a Hamiltonian $\Hh$.
\begin{Thm}\label{tswf}
Let $[M,L]$ be a proper string   and let $u_0 \in L^2(\mf)$ be a  nonzero compactly supported initial profile.  Assume that $t>0$ is such that  $\Ls(\Ts(\sff_0) + t+\eps) < \infty$ for some $\eps>0$. Then, the wavefront of the solution $u$ of \eqref{seq} can be found by the formula
\begin{equation}\label{ffwf}
\sff_t =\Ls(\Ts(\sff_0) + t).
\end{equation}
\end{Thm}
We say that a Borel measure $\sigma = \upsilon\,dx + \sigma_{\bf s}$ on $\R_+$ with the density $\upsilon$ and the singular part $\sigma_{\bf s}$  belongs to the {\bf Szeg\H{o} class $\szp$} if $(x+1)^{-1} \in L^1(\sigma)$ and 
$$
\int_{\R_+}\frac{\log \upsilon(x)}{\sqrt x(x+1)}\,dx > -\infty.
$$
To each proper string $[M, L]$, one can associate the unique non-negative Borel measure $\sigma$ on $\R_+$ with $(x+1)^{-1} \in L^1(\sigma)$ called the {\bf main spectral measure} of the string $[M, L]$. 
The theorem below provides a dynamical characterization of the Szeg\H{o} class $\szp$. 
\begin{Thm}\label{ts2bis}
Let $[M,L]$ be a proper string  and let $\sigma$ be its main spectral measure. Then, $\sigma \in \szp$ if and only if for some \textup{(}and then for every\textup{)}  nonzero compactly supported initial profile $u_0 \in L^2(\mf)$ and for some \textup{(}and then for every\textup{)} $\ell > 0$ we have  
\begin{equation}\label{eqs14}
\limsup_{t \to +\infty} \|u(\cdot, t)\|_{L^2(\mf,[\sff_{t-\ell}, \sff_t])} > 0.
\end{equation}
\end{Thm}
Put differently, the result says that the spectral measure of a string $[M, L]$ belongs to $\szp$ if and only if the part of wave $u$ near its wavefront does not vanish as $t \to +\infty$. 

\medskip

 For the homogeneous string with positive constant density, we have $L = +\infty$ and $M\colon \xi \mapsto  \rho_0\xi$ where $\rho_0>0$. In that case, the propagation of the wave with the initial profile $u_0$ has the well-known ``traveling wave''  form given by d'Alembert's formula:
\begin{equation}\label{dala}
u(\xi, t) =\frac{u_0(\xi + at) + u_0(\xi - at)}{2}, \qquad t \ge 0, \qquad \xi \ge 0, \qquad a = \rho_0^{-\frac 12}, 
\end{equation}
where we extended $u_0$ to the whole real line $\R$ as an even function.  Moreover, if $u_0 \in L^2(\R_+)$, then
\begin{equation}\label{zao}
u(\xi, t) = \frac{{u_0}(\xi - at)}{2} + o(1), \qquad t \to +\infty, 
\end{equation}
where the remainder ``$o(1)$'' is with respect to the  $L^2(\R_+)$--norm. Below, we prove a similar result for arbitrary strings with spectral measures in the Szeg\H{o} class. Define
$$
\szks = \{[M, L]: \mbox{the main spectral measure of $[M,L]$ is in }\szp\}.
$$
Consider any two measurable sets
$\Omega_{\bf s}(\mf), \Omega_{\bf ac}(\mf)\subseteq \R_+$ that satisfy
\begin{equation}
\mf_{\bf s}(\Omega_{\bf ac}(\mf))=|\Omega_{\bf s}(\mf)|=0, \quad \Omega_{\bf ac}(\mf)=\R_+\backslash \Omega_{\bf s}(\mf),
\end{equation}
where $|E|$ refers to Lebesgue measure of a set $E$. In the theorem below, we denote $\Delta_{a,t} = [\Ls(t-a), \Ls(t+a)]$ for $t \ge a >0$.
\begin{Thm}\label{aps2bis}
Suppose $[M,L]$ is a proper string in the class $\szks$. Then, for each  $u_0\in L^2(\mf)$, there exists $G_{u_0} \in L^2(\R)$ such that  for every $a>0$ we have
\begin{equation}\label{eq71}
u(\xi,t) 
=\chi_{\Omega_{\bf ac}(\mf)}(\xi)\cdot \rho^{-\frac 14}(\xi)G_{u_0}(\Ts(\xi) - t) + o(1), \qquad t \to +\infty,
\end{equation}
with $o(1)$ in $L^2(\mf, \Delta_{a,t})$. 
If, moreover, $u_0$ belongs to the absolutely continuous subspace $H_{\bf ac}(\cal{S}_M)$ of the Krein string operator $\cal{S}_M$, then  \eqref{eq71} can be strengthened: $o(1)$ now holds with respect to $L^2(\mf)$--norm.
\end{Thm} 
The ``traveling wave'' $G_{u_0}$ in Theorem~\ref{aps2bis} can be explicitly written in terms of $u_0$ and Szeg\H{o} function of the spectral measure of $[M, L]$, see details in Section \ref{ks}. 
Similar results were obtained recently in a different setting (see \cite{Den_B1} and \cite{Den_M}).
\smallskip

The class $\szks$ can be described purely in terms of string's length $L$ and mass distribution  $M$. That characterization was obtained in \cite{BD2017}. Below, we give  somewhat more general version which has already been applied in the theory of quantum graphs (see \cite{KoNi21}). 
\begin{Thm}\label{ch1sz}
Let $[M, L]$ be a proper string, and let $\{\eta_n\}$ be an increasing sequence of positive numbers such that $c_1 \le \eta_{n+1} - \eta_{n} \le c_2$ for all $n \ge 0$ and some positive $c_{1}$, $c_2$. Then, we have $[M, L] \in \szks$ if and only if $\sqrt \rho \notin L^1[0,L)$ and 
\begin{equation}\label{sdk7}
\sum_{n=0}^{+\infty}\Bigl((\xi_{n+2} - \xi_n)(M(\xi_{n+2}) - M(\xi_n)) - (\eta_{n+2} - \eta_{n})^2\Bigr) < +\infty,
\end{equation}
where $\xi_n=\Ls(\eta_{n})$.
\end{Thm}
\medskip

The Lax-Phillips scattering theory for vibrating strings (see \cite{LPbook} and  \cite{Kh81},\cite{NikPavKh} for connections with the theory of $K_\Theta$-spaces, basis property of exponents, and Regge's problem) usually assumes that the string is homogeneous on a half-line $\xi>\xi_0$ for some $\xi_0$, which places it in Szeg\H{o} class, and then the solution to problem \eqref{seq}-\eqref{seq3} is studied in the energy norm $\|\cdot\|_E$ defined by (see, e.g., p. 73 in \cite{Kh81})
\[
\|u\|^2_{E}=\frac 12\left(  \int_0^\infty u_\xi^2\,d\xi+\int_0^\infty u_t^2\,d\mf\right).
\]
In contrast to that setup, we make no such assumptions on the string. We also find it more suitable to work in the original space ${L^2(\mf)}$. 

\medskip

\definecolor{ccqqqq}{rgb}{0.8,0,0}
\definecolor{cqcqcq}{rgb}{0.7529411764705882,0.7529411764705882,0.7529411764705882}

\begin{figure}
\centering
\resizebox{\textwidth}{!}{
\begin{tikzpicture}[line cap=round,line join=round,>=triangle 45,x=2cm,y=1cm]
\draw [color=cqcqcq,, xstep=0.5cm,ystep=0.5cm] (-0.25,-1) grid (6.25,1);
\clip(-1,-1) rectangle (7,1);
\draw [line width=2.8pt,color=ccqqqq] (0,0.5)-- (0.23853714038534002,0.5);
\draw [line width=2.8pt] (0.23853714038534002,0.5)-- (1,0.5);
\draw [line width=2.8pt,color=ccqqqq] (1,0.5)-- (1.6777197437266214,0.5);
\draw [line width=2.8pt] (1.6777197437266214,0.5)-- (2,0.5);
\draw [line width=2.8pt,color=ccqqqq] (2,0.5)-- (2.8088321512068997,0.5);
\draw [line width=2.8pt] (2.8088321512068997,0.5)-- (3,0.5);
\draw [line width=2.8pt,color=ccqqqq] (3,0.5)-- (3.8687548602458772,0.5);
\draw [line width=2.8pt] (3.8687548602458772,0.5)-- (4,0.5);
\draw [line width=2.8pt,color=ccqqqq] (4,0.5)-- (4.902133096181497,0.5);
\draw [line width=2.8pt] (4.902133096181497,0.5)-- (5,0.5);
\draw [line width=2.8pt,color=ccqqqq] (5,0.5)-- (5.923032714743279,0.5);
\draw [line width=2.8pt] (5.923032714743279,0.5)-- (6,0.5);
\draw [line width=2.8pt,color=ccqqqq] (0,-0.5)-- (0.4201743134274646,-0.5);
\draw [line width=2.8pt] (0.4201743134274646,-0.5)-- (1,-0.5);
\draw [line width=2.8pt,color=ccqqqq] (1,-0.5)-- (1.7922708728327308,-0.5);
\draw [line width=2.8pt] (1.7922708728327308,-0.5)-- (2,-0.5);
\draw [line width=2.8pt,color=ccqqqq] (2,-0.5)-- (2.890364560763455,-0.5);
\draw [line width=2.8pt] (2.890364560763455,-0.5)-- (3,-0.5);
\draw [line width=2.8pt,color=ccqqqq] (3,-0.5)-- (3.9310988531636832,-0.5);
\draw [line width=2.8pt] (3.9310988531636832,-0.5)-- (4,-0.5);
\draw [line width=2.8pt,color=ccqqqq] (4,-0.5)-- (4.952110345684899,-0.5);
\draw [line width=2.8pt] (4.952110345684899,-0.5)-- (5,-0.5);
\draw [line width=2.8pt,color=ccqqqq] (5,-0.5)-- (5.964456222001263,-0.5);
\draw [line width=2.8pt] (5.964456222001263,-0.5)-- (6,-0.5);
\begin{scriptsize}
\draw [fill=red] (0,0.5) circle (2pt);
\draw [fill=black] (1,0.5) circle (1.2pt);
\draw [fill=black] (2,0.5) circle (1.2pt);
\draw [fill=black] (3,0.5) circle (1.2pt);
\draw [fill=black] (4,0.5) circle (1.2pt);
\draw [fill=black] (5,0.5) circle (1.2pt);
\draw [fill=black] (6,0.5) circle (1.2pt);
\draw [fill=black] (0,0.5) circle (1pt);
\draw[color=black] (0.08765063241805225,0.6688941712622664) node {0};
\draw [fill=black] (0.23853714038534002,0.5) circle (1pt);
\draw [fill=black] (1.6777197437266214,0.5) circle (1pt);
\draw [fill=black] (2.8088321512068997,0.5) circle (1pt);
\draw [fill=black] (3.8687548602458772,0.5) circle (1pt);
\draw [fill=black] (4.902133096181497,0.5) circle (1pt);
\draw [fill=black] (5.923032714743279,0.5) circle (1pt);
\draw [fill=black] (0.23853714038534002,0.5) circle (1pt);
\draw [fill=black] (1,0.5) circle (1pt);
\draw[color=black] (1.081495134191589,0.6688941712622664) node {1};
\draw [fill=black] (1.6777197437266214,0.5) circle (1pt);
\draw [fill=black] (2,0.5) circle (1pt);
\draw[color=black] (2.0861422935930993,0.6688941712622664) node {2};
\draw [fill=black] (2.8088321512068997,0.5) circle (1pt);
\draw [fill=black] (3,0.5) circle (1pt);
\draw[color=black] (3.0907894529946094,0.6688941712622664) node {3};
\draw [fill=black] (3.8687548602458772,0.5) circle (1pt);
\draw [fill=black] (4,0.5) circle (1pt);
\draw[color=black] (4.084633954768146,0.6688941712622664) node {4};
\draw [fill=black] (4.902133096181497,0.5) circle (1pt);
\draw [fill=black] (5,0.5) circle (1pt);
\draw[color=black] (5.089281114169656,0.6688941712622664) node {5};
\draw [fill=black] (5.923032714743279,0.5) circle (1pt);
\draw [fill=black] (6,0.5) circle (1pt);
\draw[color=black] (6.083125615943193,0.6688941712622664) node {6};
\draw [fill=red] (0,-0.5) circle (2pt);
\draw [fill=black] (1,-0.5) circle (1.2pt);
\draw [fill=black] (2,-0.5) circle (1.2pt);
\draw [fill=black] (3,-0.5) circle (1.2pt);
\draw [fill=black] (4,-0.5) circle (1.2pt);
\draw [fill=black] (5,-0.5) circle (1.2pt);
\draw [fill=black] (6,-0.5) circle (1.2pt);
\draw [fill=black] (0,-0.5) circle (1pt);
\draw[color=black] (0.08765063241805225,-0.3357529881392607) node {0};
\draw [fill=black] (0.4201743134274646,-0.5) circle (1pt);
\draw [fill=black] (1.7922708728327308,-0.5) circle (1pt);
\draw [fill=black] (2.890364560763455,-0.5) circle (1pt);
\draw [fill=black] (3.9310988531636832,-0.5) circle (1pt);
\draw [fill=black] (4.952110345684899,-0.5) circle (1pt);
\draw [fill=black] (5.964456222001263,-0.5) circle (1pt);
\draw [fill=black] (0.4201743134274646,-0.5) circle (1pt);
\draw [fill=black] (1,-0.5) circle (1pt);
\draw[color=black] (1.081495134191589,-0.3357529881392607) node {1};
\draw [fill=black] (1.7922708728327308,-0.5) circle (1pt);
\draw [fill=black] (2,-0.5) circle (1pt);
\draw[color=black] (2.0861422935930993,-0.3357529881392607) node {2};
\draw [fill=black] (2.890364560763455,-0.5) circle (1pt);
\draw [fill=black] (3,-0.5) circle (1pt);
\draw[color=black] (3.0907894529946094,-0.3357529881392607) node {3};
\draw [fill=black] (3.9310988531636832,-0.5) circle (1pt);
\draw [fill=black] (4,-0.5) circle (1pt);
\draw[color=black] (4.084633954768146,-0.3357529881392607) node {4};
\draw [fill=black] (4.952110345684899,-0.5) circle (1pt);
\draw [fill=black] (5,-0.5) circle (1pt);
\draw[color=black] (5.089281114169656,-0.3357529881392607) node {5};
\draw [fill=black] (5.964456222001263,-0.5) circle (1pt);
\draw [fill=black] (6,-0.5) circle (1pt);
\draw[color=black] (6.083125615943193,-0.3357529881392607) node {6};
\end{scriptsize}
\end{tikzpicture}
}
\caption{The upper string has black pieces of length $\frac{1}{n\log(e+n)}$. That string does not belong to the class $\szks$. The bottom string has black pieces of length $\frac{1}{n\log^2(e+n)}$. It does belong to the class $\szks$.} \label{fig:plot2_bis}
\end{figure}
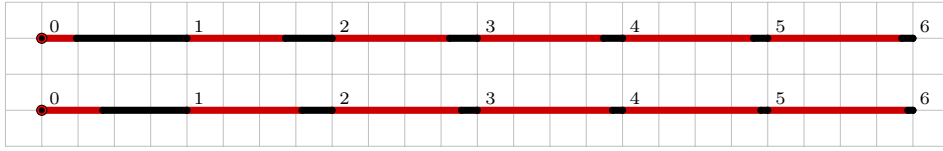

We end this section with two examples which are discussed in detail in Section \ref{ssexamples}. In the first one, the density $\rho = 1$ and we study how the properties of the string depend on $\mf_{\bf s}$, the singular component  of  $\mf$.  Measure $\mf_{\bf s}$ describes the  ``impurities'' in the material.
\begin{Ex}\label{ex2bis}
Let $[M, \infty]$ be the string with $\mf = d\xi + \mf_{\bf s}$ on $\R_+$, and let $u_0 \in L^2(\mf)$ have compact support, $u_0 \neq 0$. Then, 
$$
\sff_t=\sff_0+t, \qquad t \ge 0.
$$
If $\mf_{\bf s}(\R_+) = \infty$, we have $[M, \infty] \notin \szks$ and  
$$
\lim_{t\to\infty}\|u(\cdot,t)\|_{L^2(\mf,[\sff_t-a,\sff_t])}=0,
$$
for every $ a>0$. In the case $\mf_{\bf s}(\R_+) < \infty$, we have $[M, \infty] \in \szks$ and
\begin{align*}
\lim_{t\to+\infty}\|u(\cdot,t)\|_{L^2(\mf_{\bf s}, [\sff_t-a,\sff_t])} = 0,\\
\lim_{t\to+\infty}\|u(\cdot,t)-G_{u_0}(\cdot-t)\|_{L^2[\sff_t-a,\sff_t]}= 0,
\end{align*}
for some  $G_{u_0} \in L^2(\R), G_{u_0}\neq 0$, and all $ a>0$.
\end{Ex}
Example \ref{ex2bis} shows that the propagation of the wave depends solely on whether $\mf_{\bf s}(\R_+)$ is finite or not. 

\medskip

In our second example, $\mf_{\bf s}(\R_+)=0$ and the density $\rho$ takes two  positive values: $a$ and $b$. 
 So, we have 
\begin{equation}\label{eq68b}
\rho(\tau) 
= \begin{cases}
a, & \tau \in E,\\
b, & \tau \in \R_+ \setminus E,
\end{cases}
\end{equation}
for some Lebesgue-measurable set $E\subseteq \R_+$.
We interpret such strings as those made of two types of material (see Figure \ref{fig:plot2_bis}).
\begin{Ex}\label{ex19}
The  string with $\mf_{\bf s}(\R_+)=0$ and density $\rho$ of the form \eqref{eq68b} lies in $\szks$ if and only if either $a = b$ (the string is homogeneous) or one of the sets $E$, $\R_+\backslash E$ has finite Lebesgue measure.
\end{Ex}
Given that statement, Theorem \ref{ts2bis} and Theorem \ref{aps2bis} show that the propagation of the wave depends only on whether $\min\{|E|, |\R_+\backslash E|\}$ is finite or not.

\medskip

The paper has four parts. This introduction is followed by the second section which is focused on canonical systems. The general results obtained in the second part are applied to Krein strings in the third section and to Dirac operators in the fourth section. Finally, the Appendix contains  some auxiliary statements and proofs.  

\bigskip

\subsection{Notation}\label{notation} 
We use the following notation:

\smallskip

\noindent $\bullet$ $\R_+=[0,\infty)$, $\C_+ = \{z \in \C: \Im z > 0\}$, $\C_- = \{z \in \C: \Im z < 0\}$, $\T=\{z\in \C: |z|=1\}$.

\noindent $\bullet$ We define the direct and inverse Fourier transforms by:
\[
\widehat{f}(\xi)=\frac{1}{\sqrt{2\pi}}\int_{\R} f(x)e^{-ix\xi}\,dx, \quad \widecheck{f}(\xi)=\frac{1}{\sqrt{2\pi}}\int_{\R} f(x)e^{ix\xi}\,dx, \qquad \xi \in \R.
\]

\noindent $\bullet$ The symbol $C(\R)$ denotes the set of continuous functions on $\R$.

\noindent $\bullet$ The symbol $C^\infty_c(\R)$ denotes the set of compactly supported infinitely smooth functions on $\R$. Similarly, the symbol $L^2_c(\Hh)$ stands for compactly supported functions in $L^2(\Hh)$, and $H_c$ denotes the set of compactly supported elements in the space $H$ (see Section \ref{cs-basics} for definition of these spaces).

\noindent $\bullet$ We write $f\in L^1_{\rm loc}(\R_+)$ is $\int_0^r|f(x)|\,dx<\infty$ for every $r\in (0,\infty)$.

\noindent $\bullet$ Given a function $f\in C(\R)$, its support is defined as $\supp f=\ov{\{x: f(x)\neq 0\}}$.

\noindent $\bullet$ The symbol $C$ denotes the absolute constant which can change the value from formula to formula. 

\noindent $\bullet$  For two non-negative functions $f_1$ and $f_2$, we write $f_1\lesssim f_2$ if  there is a
constant $C$ such that $f_1\le Cf_2$ for all values of the arguments of $f_1$ and $f_2$. We define $\gtrsim$
similarly and say that $f_1\sim f_2$ if $f_1\lesssim f_2$ and
$f_2\lesssim f_1$ simultaneously. If $|f_3|\lesssim f_4$, we will write $f_3=O(f_4)$. Given $f_1$ and $f_2$,  two real-valued functions defined on $\R_+$, we write $f_1=o(f_2)$ when $x\to +\infty$ if $f_1=\alpha f_2$ for some function $\alpha$ that satisfies $\lim_{x\to+\infty}\alpha=0$.

\noindent $\bullet$ If $\mu$ is a non-negative Borel measure on the real line and $f,g\in L^2(\mu)$, we denote 
$(f,g)_\mu=\int_{\R} f\ov g\,d\mu$.

\noindent $\bullet$ We denote $f(x^*-) =\lim_{x\uparrow x^*}f(x)$. Similarly, $f'(x^*-)$ is the left derivative at point $x^*$.

\noindent $\bullet$ Given any measurable set $E\subseteq \R$, the symbol $|E|$ will denote its Lebesgue measure.

\noindent $\bullet$ Given a set $E\subseteq \R$, the symbol $\chi_E$ stands for the characteristic function of $E$.

\noindent $\bullet$ Suppose $\sigma$ is a non-negative Borel measure on $\R$ and $\sigma=\sigma_{\bf ac}+\sigma_{\bf s}$ is its decomposition into the sum of absolutely continuous and singular parts. In this paper,  $\Omega_{\bf s}(\sigma)$ and $\Omega_{\bf ac}(\sigma)$ will denote any sets that satisfy
\begin{equation}
\sigma_{\bf s}(\Omega_{\bf ac}(\sigma))=|\Omega_{\bf s}(\sigma)|=0, \quad \Omega_{\bf ac}(\sigma)=\R\backslash \Omega_{\bf s}(\sigma). \label{sa53}
\end{equation}
Analogous notation is used for measures defined on $\R_+$:
\begin{equation}
\sigma_{\bf s}(\Omega_{\bf ac}(\sigma))=|\Omega_{\bf s}(\sigma)|=0, \quad \Omega_{\bf ac}(\sigma)=\R_+\backslash \Omega_{\bf s}(\sigma). \label{sa53bis}
\end{equation}

\noindent $\bullet$ For $z\in \ov{\C_+}$, the symbol $\sqrt[n]z$ always defines the branch of the root such that $\sqrt[n]z>0$, $z\in \R_+$.

\noindent $\bullet$ Given an interval $I\subset \R$ and $f\in L^1(I)$, we write
$
\langle f\rangle_I=\frac{1}{|I|}\int_I f\,dx
$
for an average of $f$ over $I$.

\noindent $\bullet$ $\langle a, b \rangle_{\C^2} = a_1 \bar b_1 + a_2 \bar b_2$ for vectors $a, b \in \C^2$ with coordinates $a_1$, $a_2$ and $b_1$, $b_2$, correspondingly. 

\noindent $\bullet$ For $x>0$, we denote $\log^+ x=\max(0, \log x)$.

\noindent $\bullet$ The entire function $f$ has finite exponential type if 
$\type f:=\limsup_{|z|\to +\infty}\frac{\log |f(z)|}{|z|}<+\infty$.

\noindent $\bullet$ We denote $\diag(a,b)=\left(\begin{smallmatrix}a & 0\\0&b\end{smallmatrix}\right)$, $a, b\in \C.$

\noindent $\bullet$ Given a matrix $A$, we denote its transpose by $A^{t}$.

\noindent $\bullet$ Given a function $g$, we define the corresponding multiplication operator by $g$ as $M_g\colon f\mapsto gf$.

\noindent $\bullet$ 
The composition of two functions $F$ and $G$ will be denoted by $F\circ G$.

\noindent $\bullet$ 
Given a non-negative non-decreasing function $F$ defined on $\R_+$, we denote its {\bf generalized inverse} as
\[
F^{(-1)}(x)=\inf \{y\ge 0: F(y)\ge x\}
\]
and let $F^{(-1)}(x)=+\infty$ if $\{y\ge 0: F(y)\ge x\}=\emptyset$.  It can be shown that $F^{(-1)}$ is non-decreasing and left-continuous on $\R_+$. If $F$ is continuous on $\R_+$ and $F(0)=0$, then 
$F^{(-1)}(x)=\min \{y\ge 0: F(y)=x\}$,
and $F(F^{(-1)}(x))=x$ provided that $F^{(-1)}(x)<+\infty$.

\newpage

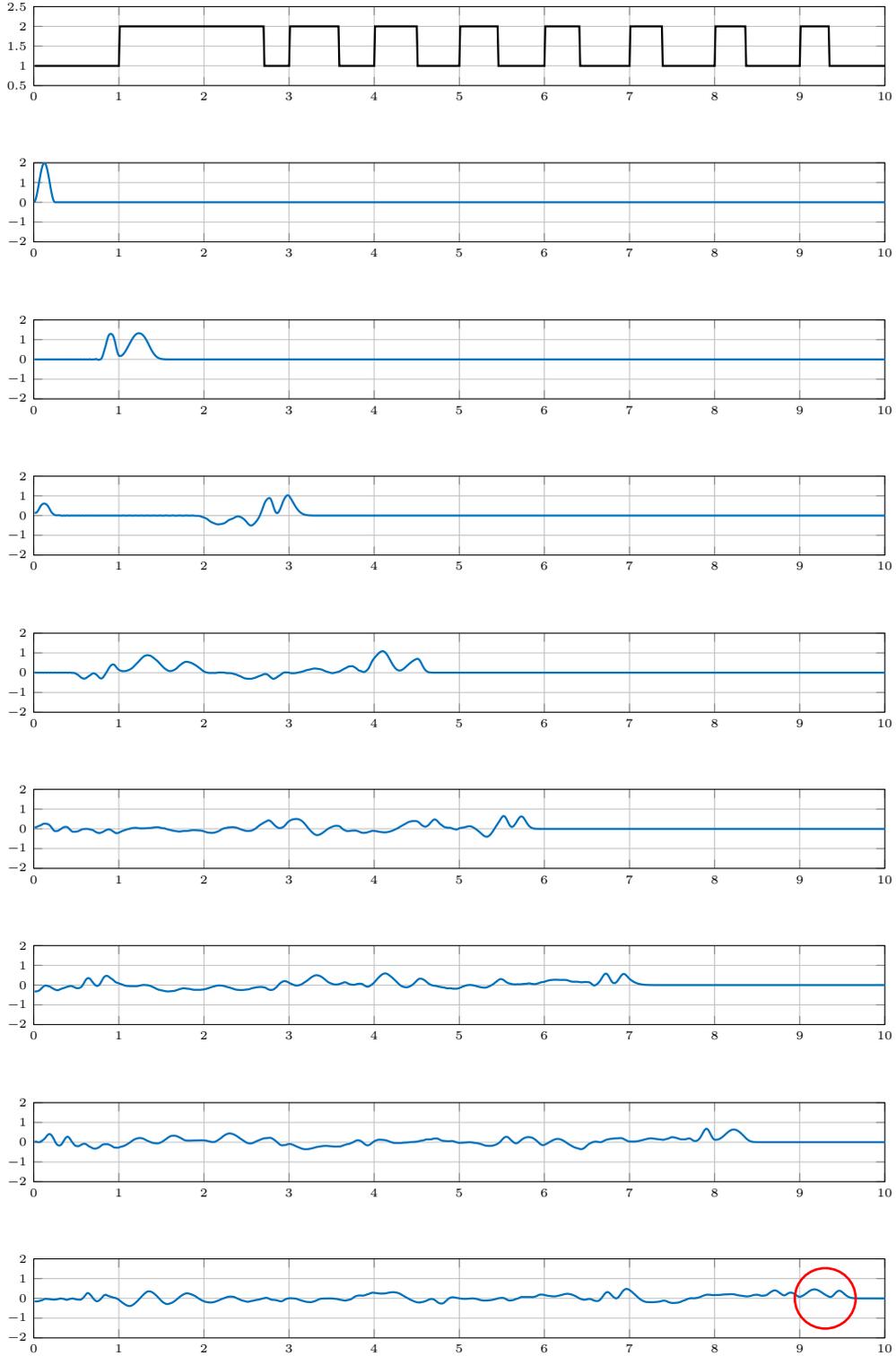
\subsection{Figure: wave propagation for a string in the non-Szeg\H{o} case}\label{s1p4} 
The first graph shows the density of the string. For each interval $[n, n+1]=E_n\cup F_n$, $E_n$ carries the density $1$, $F_n$ carries the density $2$,  and $|F_n| \sim 1/{\sqrt{n+1}}$. As time increases, only a vanishing portion of the wave (shown in the red circle) propagates with the maximal speed.

\input{FIG1.tex}

\newpage

\subsection{Figure: wave propagation for a string in the Szeg\H{o} case}\label{s1p5} The first graph shows the density of the string.  For each interval $[n, n+1]=E_n\cup F_n$, $E_n$ carries the density $1$, $F_n$ carries the density $2$. This time, $|F_n| \sim 1/{(n+1)^2}$. As time increases, a non-vanishing portion of the wave (shown in the red circle) propagates with the maximal speed.

\input{FIG2.tex}
 
\newpage

\section{Canonical Hamiltonian systems}\label{css}

\subsection{Some definitions and known results}\label{cs-basics}
We first recall some basics of the theory of canonical Hamiltonian systems \cite{Remlingb}, \cite{Romanov}. As we have seen in the Introduction, a Hamiltonian $\Hh$ on the positive half-axis $\R_+ = [0,+\infty)$ is a matrix-valued mapping of the form
$$
\Hh = \begin{pmatrix}h_1 & h\\h&h_2\end{pmatrix}, \quad \Hh(\tau) \ge 0, \quad \trace\Hh(\tau) >0, \quad \mbox{for a.e.\ }\tau \in\R_+.
$$
The functions $h_1$, $h_2$, $h$ are  real-valued and belong to $L_{\loc}^1(\R_+)$. 
If $h = 0$ almost everywhere on $\R_+$, we say that $\Hh$ is diagonal. A  Hamiltonian $\Hh$ on $\R_+$ is called {\bf  singular} if
\begin{equation}\label{eq33}
\int_{0}^{\infty}\trace\Hh(\tau)\, d\tau = + \infty.
\end{equation}
A Hamiltonian $\Hh$ is called {\bf nontrivial} if it  is {\it not} of the form $\Hh = \kappa A$ where a non-negative  function $\kappa$ is in $ L^1_{\rm loc}(\R_+)$ and a constant matrix $A \ge 0$ has rank one.  The Hilbert space $L^2(\Hh)$ is the set of (equivalence classes of) measurable vector-functions  
\begin{align}\label{eq31}
L^2(\Hh) &= \Bigl\{X \colon \R_+ \to \C^2: \int_{0}^{\infty}\bigr\langle\Hh(\tau)X(\tau), X(\tau)\bigr\rangle_{\C^2} \,d\tau < +\infty\Bigr\}\Big/{\mathcal Ker}\,\Hh,\\
{\mathcal Ker}\,\Hh &= \Bigl\{X\colon \;\Hh(\tau) X(\tau) = 0 \mbox{ for almost all } \tau \in \R_+\Bigl\}, \notag
\end{align} 
equipped with the inner product
$$
(X, Y)_{L^2(\Hh)} = \int_{0}^{\infty}\bigr\langle\Hh(\tau)X(\tau), Y(\tau)\bigr\rangle_{\C^2}\,d\tau.
$$
An open interval $I \subseteq \R_+$ is called {\bf indivisible} for $\Hh$ if there exists a function $\kappa$ and a nonzero vector $e \in \R^2$ such that $\Hh$ coincides with the operator $f \mapsto \kappa \langle f,e\rangle_{\C^2}e$ almost everywhere on $I$, and $I$ is the maximal open interval (with respect to inclusion) having this property. Let $\mathfrak{I}(\Hh)$ denote the set of all indivisible intervals of $\Hh$, and let 
\begin{equation}\label{sa1}
H = \bigl\{X \in L^2(\Hh): X = x_I \mbox{ on }  I\in\mathfrak{I}(\Hh),\;x_I \in \C^2\bigr\}.
\end{equation}
Since $L^2(\Hh)$ is a set of equivalence classes of functions, we say that $X = x_I$ on an interval $I$ if $\Hh X=\Hh x_I$ almost everywhere on $I$. We write $L^2_c(\Hh)$ and $H_c$ for compactly supported elements in $L^2(\Hh)$ and $H$, respectively.

  \medskip

\noindent In this paper, we will only work with Hamiltonians $\Hh$ that satisfy the following three conditions:
\begin{itemize}
\item[$(a)$] $\Hh$ is singular;
\item[$(b)$] for every $r\ge 0$, we have $(r, +\infty) \notin \mathfrak{I}(\Hh)$; 
\item[$(c)$] there is no $\eps>0$ such that $\Hh = \kappa\excl$ almost everywhere on $[0, \eps]$ for some function $\kappa$.
\end{itemize}
Later it will be clear that these assumptions are both convenient and natural for the kind of problems we consider in this section.  We refer to Hamiltonias satisfying $(a)$-$(c)$ as {\bf proper} Hamiltonians.

\smallskip
Fix $J = \jm$. With each Hamiltonian $\Hh$, one can associate a self-adjoint differential operator 
\begin{equation}\label{eq42}
\Di_\Hh\colon X \mapsto Y, \qquad Y: \; JX'=\Hh Y,
\end{equation}
defined on a certain dense linear subset of the Hilbert space $H$ we will introduce shortly. Note that if there are two functions $Y_1$, $Y_2$ such that  $JX' = \Hh Y_1 = \Hh Y_2$ almost everywhere on $\R_+$, then $Y_1 = Y_2$ as elements of $H$ according to \eqref{eq31} and \eqref{sa1}.  Under our assumptions on $\Hh$,  the domain of $\Di_{\Hh}$ is given by 
$$
\dom \Di_{\Hh} = \left\{X \in H: \; 
\begin{aligned}
&X \mbox{ is locally absolutely continuous on  }\R_+,\\ 
&JX' = \Hh Y \mbox{ for some } Y \in H,\\
&\langle X(0),\zo \rangle_{\C^2} = 0.
\end{aligned}
\right\}.
$$
In the first line above, ``$X$ is locally absolutely continuous''  means that there is a locally absolutely continuous representative of $X$ and then the boundary value $X(0) \in \C^2$ is defined for this representative. When considered on $\dom \Di_{\Hh}$, the operator $\Di_\Hh$ is in fact a self-adjoint operator densely defined on $H$ (e.g., check Section 2 of \cite{Romanov}).  Our assumption $(c)$ is related to the choice of the boundary condition $\langle X(0),\zo \rangle_{\C^2} = 0$ (see Theorem 3 in \cite{Romanov}). 

\medskip

Alternatively, the spectral theory of canonical Hamiltonian systems can be presented in the language of symmetric linear relations defined on the whole space $L^2(\Hh)$, not just on its subspace $H$. That approach was pioneered by I.~Kats \cite{Kats}. More details, including historical remarks, can be found in \cite{LeMa}, \cite{Remlingb}. 

\medskip

\noindent A Hamiltonian $\Hh$ on $\R_+$ generates a canonical system -- the differential equation of the form
\begin{equation}\label{cs}
J \Theta'(\tau,z) = z \Hh(\tau) \Theta(\tau,z), \quad \Theta(0, z) = \oz, \quad
\quad \tau \in \R_+, \quad z \in \C.
\end{equation}
As we mentioned in the Introduction, it can be considered as the eigenvalue problem for $\Di_{\Hh}$. Indeed, if $\Theta(\cdot, z) \in \dom \Di_{\Hh}$, then $z$ is a eigenvalue of $\Di_{\Hh}$ and $\Theta(\cdot, z)$ is an eigenfunction. Since $\Hh\in L^1_{\rm loc}(\R_+)$,   the Cauchy problem \eqref{cs} has the locally absolutely continuous (with respect to $\tau$) solution $\Theta$ for each $z \in \C$. It is also easy to see that for fixed $\tau \ge 0$, this solution is an entire $\C^{2}$-valued function with respect to $z$. We will use notation $\Theta^{+}$ and $\Theta^-$ for its entries: 
\begin{equation}\label{r6}
\Theta(\tau, z) = \begin{pmatrix}\Theta^+(\tau,z)\\ \Theta^{-}(\tau,z)\end{pmatrix}.
\end{equation}
The Titchmarsh-Weyl transform (or the ``generalized Fourier transform'') associated with $\Hh$ is densely defined by 
\begin{equation}\label{wt}
\W_{\Hh}\colon X \mapsto \frac{1}{\sqrt{\pi}}\int_{0}^{\infty}\langle\Hh(\tau)X(\tau), \Theta(\tau, \bar z)\rangle_{\C^2}\,d\tau, \qquad z \in \C,  
\end{equation}
on the set of elements $X \in H_{c}$. For such $X$,  $\W_{\Hh}X$ is an entire function with respect to $z$. It is known (see Section 9 in \cite{Romanov} or \cite{Winkler95}) that for every singular Hamiltonian there exists a unique  measure $\mu$ on $\R$ such that $(1 + x^2)^{-1} \in L^1(\mu)$ and the mapping $\W_{\Hh}$ is the unitary operator from $H$ onto $L^2(\mu)$. Note that for $z \in \C$ and $X \in \dom \Di_{\Hh} \cap H_c$ we have 
\begin{align}
z \W_{\Hh}X
&= \frac{1}{\sqrt{\pi}}\int_{0}^{\infty}\langle\Hh(\tau) X(\tau), \bar z\Theta(\tau, \bar z)\rangle_{\C^2} \, d\tau = \frac{1}{\sqrt{\pi}}\int_{0}^{\infty}\langle X(\tau), J\Theta'(\tau, \bar z)\rangle_{\C^2} \, d\tau \notag\\
&=\frac{1}{\sqrt{\pi}}\int_{0}^{\infty}\langle JX'(\tau), \Theta(\tau, \bar z)\rangle_{\C^2} \, d\tau
=\W_{\Hh}(\Di_{\Hh} X). \label{eq20}
\end{align}
Therefore, the operator $\W_{\Hh}\colon H \to L^2(\mu)$ diagonalizes $\Di_{\Hh}$ and $\mu$ is the spectral measure for $\Di_{\Hh}$ (see Section 8 in \cite{Romanov}). That measure $\mu$ is often called the main spectral measure of $\Di_{\Hh}$ or simply the spectral measure of the Hamiltonian $\Hh$. 

\medskip

\noindent Take $r$ that satisfies
\begin{equation}\label{sa02}
 r\in \R_+ \setminus \bigcup_{I \in \mathfrak{I}(\Hh)} I.
 \end{equation}
 For such $r$ the multiplication operator $Y \mapsto \chi_{[0,r]}Y$ by the characteristic function of $[0,r]$ acts from $H$ to $H$ since the multiplication with such $\chi_{[0, r]}$ is consistent with the condition that $X = x_I$ on $I \in \mathfrak{I}(\Hh)$ in \eqref{sa1},  the definition of $H$. Consider the space 
\begin{equation}\label{eq15}
\B_{r} = \W_{\Hh} H_r, \qquad H_r = \{Y \in H: \esssupp Y \subseteq [0, r]\}.  
\end{equation} 
Since $\W_{\Hh}$ is a unitary map, the set $\B_r$ is a Hilbert space of entire functions with respect to the inner product
\[
(f_1,f_2)_\mu=\int_{\R} f_1\bar{f}_2\,d\mu
\]
inherited from $L^2(\mu)$. It is called the {\bf de Branges space} generated by the restriction of $\Hh$ to $[0,r]$. Given an entire function $f$, we let $f^\sharp$ denote the function $z \mapsto \overline{f(\bar z)}$. Then, $f = f^\sharp$ if and only if $f$ is real on the real line $\R$. We can define
\begin{equation}\label{eq49}
E_{r}(z)= \Theta^+(r,z) + i\Theta^{-}(r, z), \qquad E_{r}^\sharp(r, z) = \Theta^+(r,z) - i\Theta^{-}(r, z), \qquad z \in \C,
\end{equation}
{where} the second formula follows from the fact that $\Theta^+(r,z)$ and $\Theta^-(r,z)$ are real for real $z$.  It is known (see, e.g., Section 4.3 in \cite{Remlingb}) that $E_r$ has no roots in the upper half-plane $\C_+$ and that $\B_r$ admits the following description in terms of the Hardy space $H^2(\C_+)$: 
\begin{equation}\label{eq50} 
\B_r = \left\{\mbox{entire }f\colon \; \frac{f}{E_r} \in H^2(\C_+), \; \frac{f^\sharp}{E_r} \in H^2(\C_+)\right\}.  
\end{equation}
Moreover, 
\begin{equation}\label{eq51} 
\int_{\R}|f|^2\,d\mu = \int_{\R}\left|\frac{f}{E_r}\right|^2\,dx, \qquad f \in \B_r.
\end{equation}
An immediate corollary of \eqref{eq15} and \eqref{eq51} is the following nesting property of subspaces $\B_r$: if $0 \le r_1 < r_2$ and both $r_1$ and $r_2$ satisfy \eqref{sa02}, then we have the isometric inclusion $\B_{r_1} \subsetneq \B_{r_2}$. In particular,
\begin{equation}\label{eq52}
\int_{\R}\left|\frac{f}{E_{r_1}}\right|^2\,dx =  \int_{\R}\left|\frac{f}{E_{r_2}}\right|^2\,dx=\int_{\R}|f|^2\,d\mu, \qquad f \in \B_{r_1}.
\end{equation}
We say that an entire function $E$  is a {\bf Hermite-Biehler function}  if
\begin{equation}\label{eq55}
|E(z)| > |E(\bar z)|, \qquad z \in \C_+.
\end{equation}
A Hermite-Biehler function $E$ is called {\bf regular}  if 
\begin{equation}\label{eq54}
\frac{1}{(z+i)E} \in H^2(\C_+).
\end{equation}
For each $r \ge 0$, the function $E_r$ in \eqref{eq49} is known to be a regular Hermite-Biehler function, see  Proposition~6 in~\cite{Romanov}.

\medskip

Using description \eqref{eq50}   of the space $\B_r$ and formula \eqref{eq31}, it is easy to check (see Theorem 4.4 in \cite{Remlingb}) that the Hilbert space $\B_r$ has a reproducing kernel at each point $\lambda \in \C$:
\begin{equation}\label{eq11}
k_{\B_r, \lambda}\colon z \mapsto - \frac{1}{2\pi i} \frac{E_r(z)\ov{E_r(\lambda)} - E_r^\sharp(z)\ov{E_r^\sharp(\lambda)}}{z - \bar \lambda}, \qquad z \in \C. 
\end{equation}
The latter means that $k_{\B_r, \lambda} \in \B_r$ and $(f, k_{\B_r, \lambda})_\mu = f(\lambda)$ for every $f \in \B_r$.

\medskip

\subsection{The Krein-de Branges theorem and  front of the wave}\label{s22} Let a proper Hamiltonian $\Hh$ be given. In the context of the general problem considered in Section \ref{sect11}, we
 define the unitary group $U_t=e^{it \Di_{\Hh}}$ for all $t\in \R$ using the spectral theorem. In the current subsection, we study evolution $U_tX$ for $X\in H_c$. We define the {\bf front} of $X$ as follows
\[
\sff[X]=\inf\{\ell \ge 0: \Hh(\tau)X(\tau) = 0\, \, \,\text{for a.e. } \tau>\ell\}.
\]
For $t \in \R$, we will refer to the number $\sff[U_t X]$ as the {\bf wavefront} of wave function $U_t X$. 

The next theorem gives the formula for the wavefront in terms of two auxiliary functions $\Tc(\tau)$ and $\Lc(\eta)$. For $\tau, \eta \in \R_+$, they are defined as follows
\begin{equation}\label{sa21}
\Tc(\tau) = \int_{0}^{\tau}\sqrt{\det\Hh(s)}\,ds,  
\qquad
\Lc(\eta) = \Tc^{(-1)}(\eta)=\inf \left\{\tau \ge 0:\, \Tc(\tau) \ge  \eta\right\}.
\end{equation}
$\Lc$ is the generalized inverse of $\Tc$ and, if $\eta > 0$ and the set $\{\tau\ge 0: \Tc(\tau) \ge \eta\}$ is empty, we let $\Lc(\eta) = +\infty$. The latter can happen only if $\sqrt{\det\Hh} \in L^1(\R_+)$. Notice that $\Tc(\Lc(\eta)) = \eta$
provided that $\eta\ge 0$ and $\Lc(\eta)<+\infty$. Moreover, when $\det\Hh> 0$ almost everywhere on $\R_+$, we have $\Lc(\Tc(\tau)) = \tau$ for each $\tau \ge 0$. 
Informally, for every $\tau_0>0$, the quantity $\Tc(\tau_0)$ is equal to the time it takes  for a wave  $U_t X$ to travel from $0$ to the point $\tau_0$. 

\medskip

Later in the text,  an element $X\in H$ is called real if both components of $X$ are real-valued. 
\begin{Thm}\label{p0bb}
Let $X \in H_c$ be real, $\Tc(\sff[X]) = a$, $a>0$. Assume that $t \in \R\setminus\{0\}$ is such that there is $\eps > 0$ such that $\Lc(|t| + a + \eps) < \infty$. Then,
\begin{equation}
\sff[U_t X] = \Lc(|t|+a), \quad \sff[(U_t+U_{-t})X] = \Lc(|t|+a). \label{add1}
\end{equation}
In particular, we have 
 \[
 \Tc(\sff[U_t X]) = |t|+a, \,\quad \Tc(\sff[(U_t+U_{-t}) X]) = |t|+a,
 \] for every such $t$.
\end{Thm}  
Recall that $\widecheck g$ denotes the inverse Fourier transform of a function $g$ as defined in Section \ref{notation}. Consider the sets
\begin{equation}\label{eq8}
\E_{(0, \infty)} = \left\{\widecheck g \colon g\in C_c^\infty(\R), \, \supp g\subseteq (0,\infty)\right\}, \qquad
\E_{b} = \left\{\widecheck g: \, g\in C^\infty_c(\R), \,\supp g \subseteq (-b,b)\right\}.
\end{equation} 
To obtain Theorem \ref{p0bb}, we will need a few results from complex analysis. For the proof of the following theorem, see  Section 4.2 in \cite{DymMcKean} or Theorem A.6 in \cite{Den06}.
\begin{Thm}[Krein-Wiener theorem]\label{nmt0}
Let $\mu$ be a measure on $\R$ such that $(1 + x^2)^{-1} \in L^1(\mu)$. Then, $\mu \in \sz$ if and only if  $\E_{(0, \infty)}$ is not dense in $L^2(\mu)$.
\end{Thm}
For the proof of the following theorem, see Section 6.4 in \cite{DymMcKean}, p.~241. 
\begin{Thm}[Krein's alternative]\label{nmt4}
Let $\mu$ be a measure on $\R$ such that $(1 + x^2)^{-1} \in L^1(\mu)$. Take any $b>0$. Then, either the set 
$\clos_{L^2(\mu)}\bigcap_{\eps>0}\E_{b+\eps}$ coincides with $L^2(\mu)$ or it is equal \textup{(}in $L^2(\mu)$\textup{)} to the set of all entire functions of type at most $b$ that belong to $L^2(\mu)$ when restricted to the real line.
\end{Thm}
The proof of Theorem~\ref{nmt4} in \cite{DymMcKean} contains a step ``$\textbf{Z}^T \subset \textbf{I}^T$ because $\textbf{I}^T$ is closed provided it is not dense in $L^2(\mu)$'' (we use the notation on p.~110 of that book).
The reader can find a more detailed proof to that claim
 in Appendix II of \cite{B2018} (see formula $(65)$ in the proof of Proposition~2.5 on p.~300 therein). In turn, Appendix II of \cite{B2018} uses ideas from Section 5.2 of \cite{BS11}. 

\medskip

The next result was announced by M.~Krein in \cite{Krein51} and was proved independently by de Branges \cite{dbii}. A short proof by Romanov can be found in  Section 6 of \cite{Romanov}, see also Section 5 in \cite{BLY1}. 
\begin{Thm}[Direct Krein-de Branges theorem on exponential type]\label{nmt1}
For every $r$ that satisfies \eqref{sa02}, the entire functions in  the space $\B_r$ defined by \eqref{eq15} have the first order and a finite exponential type. Moreover, 
\begin{equation}\label{eq53}
\Tc(r) = \max\{\type f, f\in \B_r\} = \type E_r = \lim_{y \to +\infty}\frac{\log |E_{r}(iy)|}{y},
\end{equation}
for the Hermite-Biehler function $E_r$ in \eqref{eq49} generating $\B_r$.
\end{Thm}

\medskip

The following result is folklore. See Appendix II in \cite{B2018} for its proof.
\begin{Thm}[Inverse Krein-de Branges theorem on exponential type]\label{nmt2}
Let $\Hh$ be a proper Hamiltonian on $\R_+$ and let $\mu$ be its spectral measure. If the set
$\E_{b}$ is not dense in $L^2(\mu)$ for some $b>0$, then the completion of $\E_b$
with respect to the inner product of $L^2(\mu)$ coincides with $\B_{r}$ for $r = \Lc(b)$.
\end{Thm}

For Krein strings, the next statement was proved in Section 6.4 of \cite{DymMcKean}.
We give a sketch of a similar argument in the case of canonical systems and their de Branges spaces.
\begin{Prop}\label{prop1} 
Let $\Hh$ be a proper Hamiltonian on $\R_+$ and let $\mu$ be its spectral measure. Assume that for some positive $b$ and $\eps$, we have $\Lc(b+\eps) < \infty$. Let $\B$ denote the  set of all entire functions of exponential type at most $b$ that belong to $L^2(|E_{\Lc(b+\eps)}(x)|^{-2})$ when restricted to {the} real line. Then, $\B = \B_{r_0}$ for $r_0 = \lim_{\delta \downarrow 0}\Lc(b+\delta)$. In particular, we have $\B \subseteq \B_{\Lc(b +\eps)}$.
\end{Prop}
\beginpf 
By \eqref{eq54}, the measure $\mu = |E_{\Lc(b+\eps)}|^{-2}\,dx$ belongs to the class $\sz$. Using Theorem~\ref{nmt0} and Theorem~\ref{nmt4} for this choice of $\mu$, we see that the set $\B$ is a Hilbert space of entire functions with respect to the inner product inherited from $L^2(|E_{\Lc(b+\eps)}|^{-2}\,dx)$.  Moreover, $\B$  satisfies the ``axiomatic'' description of de Branges spaces summarized in the following properties:
\begin{itemize}
\item[$(A_1)$] whenever $f$ is in the space and has a non-real zero $w$, the function $\frac{z - \bar w}{z - w}f$ is in
the space and has the same norm as $f$;
\item[$(A_2)$] for every non-real number $w$, the evaluation functional $f \mapsto f(w)$ is continuous;
\item[$(A_3)$] the function $f^\sharp$ belongs to the space whenever $f$ belongs to the space and it always has the same norm as $f$.
\end{itemize}  
For the proof of $(A_2)$, see formula $(65)$, p.~300, in \cite{B2018}. Notice that $\B$ satisfies an additional property: the function 
$z \mapsto \frac{f(z) - f(\lambda)}{z-\lambda}$ belongs to $\B$ for every $f \in \B$ and $\lambda \in \C$. In other words, $\B$ is the so-called {\bf regular de Branges space} isometrically embedded in $L^2(|E_{\Lc(b+\eps)}|^{-2}\,dx)$. The same is true for every space $\B_{r}$, $r \in \R_+ \setminus \bigcup_{I \in \mathfrak{I}(\Hh)} I$. Then, the de Branges ordering theorem for regular spaces states that for every $r\ge 0$ we have either $\B \subseteq \B_{r}$ or $\B_r \subseteq \B$. {Take $r_0$ as in the statement of the proposition.} Comparing the maximal exponential types of functions in $\B$ and $\B_{\Lc(b+\eps')}$ for positive $\eps'$ and using Theorem \ref{nmt1}, we get 
$$
\B \subseteq \bigcap_{\eps' > 0} \B_{\Lc(b+\eps')} = \W_{\Hh}\left(\bigcap_{\eps' > 0} H_{\Lc(b+\eps')}\right) = 
\W_{\Hh}H_{r_0} = \B_{r_0}.
$$
On the other hand, for every $f \in \B_{r_0}$, we have $f \in L^2(|E_{\Lc(b+\eps)}|^{-2}\,dx)$ by \eqref{eq52} and $\type f \le \Tc(r_0) = b$ by \eqref{eq53}. Hence, $\B_{r_0} \subseteq \B$ and the result follows. \qed

\medskip

\begin{Prop}\label{prop2}
For $X \in H_c$, $\W_{\Hh} X$ is an entire function of finite exponential type which can be computed by the formula
\begin{equation}\label{eq10}
\type \W_{\Hh} X = \Tc(\sff[X]) = \int_{0}^{\sff[X]}\sqrt{\det \Hh(\tau)}\,d\tau.
\end{equation}
\end{Prop}  
\beginpf Take an element $X \in H_c$. The definition of $\W_{\Hh}$ shows that $\W_{\Hh} X$ is an entire function. By Theorem \ref{nmt1}, $\W_{\Hh} X$ has first order and is of finite exponential type. {Moreover, by} Theorem \ref{nmt1}, we have $\type \W_{\Hh} X \le \Tc(\sff[X])$. To prove that this inequality is in fact equality, assume that $\type \W_{\Hh} X < \Tc(\sff[X])$. Then, there is a number $r < \sff[X]$ such that $\type \W_{\Hh} X < \Tc(r) < \Tc(\sff[X])$ and $\Lc(\Tc(r)) = r$. Consider the space $\B$ of all entire functions $f$ of exponential type at most $\type \W_{\Hh}X$ such that $f$ belongs to $L^2(|E_{\Lc(\Tc(r))}|^{-2}\,dx) = L^2(|E_{r}|^{-2}\,dx)$. Proposition~\ref{prop1} shows that $\B \subseteq \B_{\Lc(\Tc(r))} = \B_r$.  Since $\W_{\Hh} X \in \B$ by construction, it follows that $X$ belongs to $\W_{\Hh}^{-1}\B_r =H_{r}$. The latter contradicts that $r < \sff[X]$ and so \eqref{eq10} holds. \qed

\medskip

Now, we are ready to prove Theorem \ref{p0bb}.

\medskip

\noindent {\bf Proof of Theorem \ref{p0bb}.} We will do the proof for $U_tX$, the argument for $(U_t+U_{-t})X$ is identical. Consider a real  element $X \in H_c$ and set $f = \W_{\Hh} X$.  Since $X$ is real,  $f$ is an entire function taking real values on $\R$ and so $f = f^\sharp$. From \eqref{eq10}, one has $\type f = \Tc(\sff[X]) = a$. We claim that function $f$ is of a bounded type both in the lower and upper half-planes $\C_{\pm}$. Indeed, if we put $E(z) = \Theta^+(\sff[X], z) + i\Theta^-(\sff[X], z)$, then $E$ is an entire function of bounded type in $\C_{+}$ and it has no zeroes  there (see Theorem 4.19 in \cite{Remlingb}). 
Similarly, $E^\sharp$ is an entire functions of bounded type in $\C_{-}$ without zeroes in $\C_{-}$.
From \eqref{eq50}, we get $f/E\in H^2(\C_+)$ and $f/E^\sharp\in H^2(\C_-)$.  Since functions in $H^2(\C_{\pm})$ have bounded type in $\C_{\pm}$, the product $f = E\cdot(f/E) = E^\sharp\cdot(f/E^\sharp)$ has bounded type in $\C_\pm$ as well and the claim is proved. For every entire function which is of bounded type in both  $\C_+$ and $\C_-$, its exponential type can be computed by the formula
\begin{equation}\label{eq21}
\type f = \limsup_{y \to +\infty}\frac{\log\max(|f(iy)|, |f(-iy)|)}{y},
\end{equation}
(we sketch the proof of that known identity in Section \ref{app3}).
In our case, $|f(iy)| = |f(-iy)|$ and the same formula gives $\type(e^{itz} f) = |t| + a$. By our assumption, there exists $\eps > 0$ such that $\Lc(|t| + a + \eps) < \infty$. Then, by Proposition \ref{prop1}, the set of all entire functions of exponential type at most $|t|+a$ that belong to $L^2(|E_{\Lc(|t|+a+\eps)}|^{-2}\,dx)$ coincides with $\B_{r}$ where $r = \lim\limits_{\delta \to 0, \delta > 0}\Lc(|t|+a+\delta)$. Note that by \eqref{eq52} we have
$$
\int_{\R}\Bigl|\frac{e^{itx}f(x)}{E_{\Lc(|t|+a+\eps)}(x)}\Bigr|^2\,dx = \int_{\R}\Bigl|\frac{f(x)}{E_{\Lc(|t|+a+\eps)}(x)}\Big|^2\,dx< \infty.
$$
It follows that $e^{itz}f \in \B_{r}$. Let $Y \in H_r$ be such that $\W_{\Hh}Y = e^{itz}f$. By the spectral theorem, we have $\W_{\Hh} (U_t X) = e^{itx} f$ and that function is an element of $L^2(\mu)$. Therefore, we have $\W_{\Hh} (U_t X) = \W_{\Hh}Y$ in $L^2(\mu)$, and hence $U_t X = Y$ belongs to $H_r$. In particular, there exists a representative of $U_t X$ in $L_c^2(\Hh)$ and we have $\W_{\Hh} (U_t X)(z) = e^{itz}f(z)$ everywhere in~$\C$. Then, 
$$
\type \W_{\Hh} (U_t X) = \type(e^{itz} f) = |t| + a.
$$
The formula \eqref{eq10} gives $\Tc(\sff[U_t X])= |t| + a$. Hence, $\Lc(|t| + a) \le \sff[U_t X]$ by the definition of function~$\Lc(\eta)$. If the solution $\tau$ to the equation $\Tc(\tau) = |t|+a$ is unique, then we immediately have $\sff[U_t X] = \tau = \Lc(|t| + a)$. In the general case, for every $t \neq 0$ we can find a sequence $t_n$ such that $t_n \to t$, $|t_n| < |t|$ and the equation $\Tc(\tau) = |t_n|+a$ has unique solution for each $n$ (here we use the fact that $t \neq 0$). Notice that,  by the spectral theorem,
\[
\lim_{t_n\to t}\|U_{t_n}X-U_tX\|_{L^2(\Hh)}=\lim_{t_n\to t}\|(e^{it_nx}-e^{itx})f\|_{L^2(\mu)}=0.
\]
Since $\sff[U_{t_n} X] = \Lc(|t_n| + a) \le \Lc(|t| + a)$, we obtain $\sff[U_{t} X] \le \Lc(|t| + a)$. Hence, $\sff[U_{t} X] = \Lc(|t| + a)$ and the proof is finished.
\qed

\medskip

\subsection{Spectral measures in Szeg\H{o} class and their dynamical characterization}  Recall that a measure $\mu = w\,dx + \mus$ on $\R$ with the absolutely continuous part $w\,dx$ and the singular part $\mus$  belongs to the Szeg\H{o} class $\sz$ if $(x^2+1)^{-1} \in L^1(\mu)$ and 
$$
\int_{\R}\frac{\log w(x)}{x^2+1}\,dx > -\infty.
$$
Since $(x^2+1)^{-1} \in L^1(\mu)$, the last condition is in fact equivalent to $\frac{\log w}{1+x^2} \in L^1(\R)$.
We now define a class of Hamiltonians as follows
$$
\szcs = \bigl\{\,\Hh:   \; \Hh\, \mbox{\, is proper and its main spectral measure  is  in }\sz\bigr\}.
$$
The class $\szcs$ was characterized in \cite{BD2019} (for Dirac and Schr\"odinger operators, similar results were obtained in \cite{Den06} and \cite{KS-S}). Assuming that $\sqrt{\det\Hh} \notin L^1(\R_+)$, we define
\begin{equation}\label{eq02}
\widetilde \K(\Hh) = \sum_{n=0}^\infty \left(\det \int_{\Lc(n)}^{\Lc(n+2)}\Hh(\tau)\,d\tau-4\right).
\end{equation}
It can be shown that all terms in this series are non-negative.
In particular, $\widetilde \K(\Hh) \in \R_+ \cup\{+\infty\}$ is well-defined but could be $+\infty$, in general. 
In \cite{BD2019}, we proved that
\begin{equation}\label{s_a5}
\Hh\in \szcs\Longleftrightarrow \sqrt{\det\Hh} \notin L^1(\R_+)\quad {\rm and}\quad \widetilde\K(\Hh)<+\infty.
\end{equation} 
The partition in \eqref{eq02} does not have to be done over  the integer lattice $\{0,1,2,\ldots\}$. In fact, we have
 the following result.
\begin{Prop} \label{str1}
 Consider any monotonically increasing sequence $0=\alpha_0<\alpha_1<\alpha_2<\ldots$ of real numbers $\{\alpha_n\}$ such that
$
0<C_1<\alpha_{n+1}-\alpha_n<C_2
$
for all $n$. Then, 
\begin{equation}\label{fix1}
\Hh\in \szcs\Longleftrightarrow \sqrt{\det\Hh} \notin L^1(\R_+)\quad {\rm and}\quad \sum_{n=0}^\infty \left(\det \int_{\Lc(\alpha_n)}^{\Lc(\alpha_{n+2})}\Hh(\tau)\,d\tau-(\alpha_{n+2}-\alpha_n)^2\right)<+\infty.
\end{equation}
\end{Prop}
For completeness, we give the proof of this result in Section \ref{app1}.

\medskip

Recall our convention to write $\Lc(\eta) = +\infty$ for some $\eta > 0$ if there is no $\tau > 0$ such that $\Tc(\tau) \ge \eta$. In the following theorem, we regard $[+\infty, +\infty]$ as the empty set. We also use notation $\|Y\|_{L^2(\Hh, S)}$ for the norm of a function $\chi_{S} Y$ in $L^2(\Hh)$  on $\R_+$, where $\chi_{S}$ denotes the characteristic function of a measurable set $S \subseteq \R_+$. Our next result gives a dynamical characterization of  $\szcs$. In particular, it says that the property $\Hh\in \szcs$ can be established by observing the dynamics of $U_tX$ near its wavefront for any real nonzero $X \in H_c$.

\smallskip

We will need the following notation: given three parameters $t,s,\ell$ that satisfy $s,\ell\ge 0$ and $|t|+s-\ell\ge 0$, we define $\Delta_{\ell,s,t}=\bigl[\Lc(|t| + s - \ell), \Lc(|t| + s)\bigr]$.
\begin{Thm}\label{t2}
Let $\Hh$ be a proper Hamiltonian. Suppose $\Hh \in \szcs$ and $X$ is any real nonzero element in $H_c$.  Define  $a=\Tc(\sff[X])$. Then, we have 
\begin{equation}\label{eq14}
\liminf_{t \to \pm \infty} \|U_t X\|_{L^2(\Hh, \dlto)} > 0
\end{equation}
for all $\ell>0$.  Conversely, suppose there is $X \in H_c$ such that one of the following two conditions holds
\begin{eqnarray}\label{eq141}
\limsup_{t\to +\infty} \|U_tX\|_{L^2(\Hh, \dlto)} > 0,\qquad
\limsup_{t \to +\infty} \|U_{-t}X\|_{L^2(\Hh, \dlto)} > 0,
\end{eqnarray}
for some $\ell > 0$ and $a=\Tc(\sff[X])<+\infty$.  Then, $\Hh \in \szcs$. 
\end{Thm}
\begin{Rema} Combining this result with Theorem \ref{t1} below, one can conclude that condition $\Hh\in \szcs$ actually implies that the limits
$
\lim_{t \to \pm\infty} \|U_{t}X\|_{L^2(\Hh, \dlto)}
$
exist and are positive
for every $\ell>0$ as long as $X$ is real-valued. 
\end{Rema}

\beginpf Suppose that the first bound in \eqref{eq141} holds for some $X \in H_c$. Each component of $X$ can be written as a sum of its real and imaginary parts $X=X_R+iX_I$, so
\[
0<\limsup_{t\to +\infty} \|U_tX\|_{L^2(\Hh, \dlto)}\le \limsup_{t\to +\infty} \|U_tX_R\|_{L^2(\Hh, \dlto)}+\limsup_{t\to +\infty} \|U_tX_I\|_{L^2(\Hh, \dlto)}.
\]
Thus, we can assume that $X$ is, e.g., real without loss of generality. Then, $\sqrt{\det\Hh} \notin L^1(\R_+)$ since otherwise $\dlto = \emptyset$ for large $t$. Theorem \ref{p0bb} implies 
$$
\inf\{\|U_t X - Z\|_{L^2(\Hh)}: \supp Z \subseteq [0, \Lc(t + a - \ell)]\} =\|U_t X\|_{L^2(\Hh, \dlto)}.
$$
By the spectral theorem, it means that the function $f = \W_{\Hh} X$ satisfies 
$$
\limsup_{t\to +\infty}
\left(\inf_{u \in \B_{r}}
\|e^{itx}f-u\|_{L^2(\mu)}\right) > 0, \quad r = \Lc(t + a - \ell),
$$
where $r>0$ and $\B_r$ is defined by \eqref{eq15}.  Theorem \ref{nmt2} claims that $\B_r$ coincides with the closure in $L^2(\mu)$ of the linear manifold $\E_{t+a-\ell}$ defined by \eqref{eq8}. It follows that
$$
\limsup_{t\to+\infty}\left(\inf\{
\|e^{-i(a-\ell)x}f-h\|_{L^2(\mu)} : h \in L^1(\R), \; \widehat h\in C_c^\infty(\R), \; \supp \widehat h\subseteq (-2(t+a-\ell), 0) 
\}\right) > 0.
$$
The infimum above is a non-increasing function in $t$, hence
$$
\inf\{
\|e^{-i(a-\ell)x}f-h\|_{L^2(\mu)}\colon h \in L^1(\R), \; \widehat h\in C_c^\infty(\R), \,\supp \widehat  h\subseteq
(-\infty, 0) 
\} > 0.
$$
Now, Theorem \ref{nmt0} implies $\mu \in \sz$ and, therefore, $\Hh \in \szcs$. Similarly, $\limsup\limits_{t \to -\infty} \|U_{t}X\|_{L^2(\Hh, \dlto)}> 0$ gives $\Hh\in \szcs$.

\medskip

Conversely, suppose that $\Hh \in \szcs$. Then, $\mu \in \sz$. Arguing by contradiction, let $X$ be a real element in $H_c$ such that $\Tc(\sff[X]) = a$, $a> 0$, and either
\begin{equation}\label{ch1}
\liminf_{t \to +\infty} \|U_t X\|_{L^2(\Hh, \dlto)} = 0 \,\,\quad {\rm or}\, \quad
\liminf_{t \to -\infty} \|U_t X\|_{L^2(\Hh, \dlto)} = 0
\end{equation}
for some $\ell > 0$. 
Assume that the first limit is zero, the other case can be handled similarly.
Consider the function $f = \W_{\Hh} X$ and denote $g=e^{-i(a-\ell)z}f$. Using Theorem~\ref{p0bb} as in the first part of the proof, we obtain 
$$
\eqref{ch1} \Rightarrow
\liminf_{t\to +\infty}
\left(\inf_{u \in \B_{r}}
\|e^{itx}f-u\|_{L^2(\mu)}\right) = 0, \quad r = \Lc(t + a - \ell),
$$
and
\begin{equation}\label{eq16}
\inf\{
\|g - h\|_{L^2(\mu)}\colon h \in L^1(\R),\; \widehat{h}\in C_c^\infty(\R), \; \supp \widehat{h}\subseteq 
(-\infty, 0) 
\} = 0. 
\end{equation}
Recall the decomposition of $\mu = w\,dx + \mus$ of $\mu$ into the absolutely continuous and singular parts. Now we use assumption $\mu \in \sz$. Let $\mathcal{O}$ be an outer function in  $\C_-$ such that $|\mathcal{O}|^2 = w$ almost everywhere on $\R$ in the sense of non-tangential boundary values. Then, \eqref{eq16} gives
$$
\inf\{\|g \mathcal{O} - h \mathcal{O}\|_{L^2(\R)}: h \in L^1(\R),\; \widehat{h}\in C_c^\infty(\R), \; \supp \widehat{h}\subseteq 
(-\infty, 0)   
\} = 0. 
$$ 
That implies, in particular, that $g \mathcal{O}$, when restricted to the real line, is $L^2(\R)$ function whose Fourier transform is supported on the negative half-line. In the proof of Theorem \ref{p0bb}, we showed  that $f$ is of a bounded type in $\C_-$ and $\C_+$ and that $f=f^\sharp$. Then, $g\mathcal{O}$ is also of bounded type there and, therefore, it is in fact an element of  $H^2(\C_-)$. By the Lebesgue dominated convergence theorem,  the function $\mathcal{O}$ satisfies
$$
\lim_{y \to +\infty}\frac{\log|\mathcal{O}(-iy)|}{y} = -\lim_{y \to +\infty}\frac{1}{\pi y}\int_{\R}\log|\mathcal{O}(x)|\frac{y}{x^2 + y^2}\,dx = 0\,.
$$
The same argument applies to the outer factor of $g\mathcal{O}$ in its Smirnov-Nevanlinna factorization. As a consequence, if $g\mathcal{O} = I\cdot \mathcal{O}_1$ is the inner-outer factorization of $g\mathcal{O}$ in $\C_-$, we have
\begin{equation}\label{eq43}
\limsup_{y \to +\infty}\frac{\log|g(-iy)|}{y} = \limsup_{y \to +\infty}\frac{\log|g(-iy)\mathcal{O}(-iy)|}{y}=\limsup_{y \to +\infty}\frac{\log|I(-iy) \cdot\mathcal{O}_1(-iy)|}{y} \le 0,
\end{equation}
due to the fact that $|I| \le 1$ on $\C_-$. On the other hand, we have
\begin{equation}\label{eq44}
\limsup_{y \to +\infty}\frac{\log|g(-iy)|}{y} =
-a+\ell + \limsup_{y \to +\infty}\frac{\log|f(-iy)|}{y} = -a+\ell+\type f,
\end{equation}
because $f$ is an entire function of bounded type in both $\C_{+}$ and $\C_-$, $f=f^\sharp$, and hence we can use Lemma \ref{appe1}. However $\type f = \Tc(\sff[X]) = a$ by \eqref{eq10}, and we have got a contradiction of \eqref{eq43} and \eqref{eq44} if $\ell > 0$. 
Thus, for $\mu \in \sz$ we always have $\liminf_{t \to +\infty} \|U_t X\|_{L^2(\Hh, \dlto)} > 0$. Analogously, one can show that   
$\liminf_{t \to +\infty} \|U_{-t} X\|_{L^2(\Hh, \dlto)}$
is equal to 
$$
\inf\{
\|e^{i(a-\ell)x}f - h\|_{L^2(\mu)}\colon h \in L^1(\R), \;
 \widehat h\in C_c^\infty(\R), \,\supp \widehat h\subseteq 
(0,\infty) 
\},
$$
for $f = \W_{\Hh}X$, which implies  $\liminf_{t \to +\infty} \|U_{-t} X\|_{L^2(\Hh, \dlto)} > 0$ for every $\mu \in \sz$. \qed

\medskip
We have the following two corollaries. 
\begin{Cor}\label{t2bis1}
Let $\Hh$ be a proper Hamiltonian. Suppose $\Hh$ is not in the class $\szcs$. Then, we have 
\begin{equation}\label{eq14bis1}
 \lim_{t \to \pm\infty} \|U_t X\|_{L^2(\Hh, \dlt)}=0, \quad 
\dlt = \bigl[\Lc(|t| -b), \Lc(|t| + b)\bigr]
\end{equation}
for all $b>0$ and all $X\in H$.
\end{Cor}
\begin{Rema} If $\Lc(|t|-b)=+\infty$ for some $t$ and $b$ in the formula for $\dlt$, then $\|U_t X\|_{L^2(\Hh, \dlt)}=0$ by definition.  \end{Rema}

\noindent{\bf Proof of Corollary \ref{t2bis1}.}
Arguing by contradiction, suppose there is some $X\in H$ and $b>0$ such that, e.g.,
\begin{equation}\label{eq14bis2}
 \limsup_{t \to +\infty} \|U_t X\|_{L^2(\Hh, \dlt)}>0.
\end{equation}
Given an arbitrary $\eps>0$, there is $X_\eps\in H$ such that 
\[
\|X-X_\eps\|_{L^2(\Hh)}\le \eps, \quad \Tc(\sff[X_\eps])=a_\eps<+\infty.
\]
Choosing $\eps$ such that $\eps<\frac{1}{2}\limsup_{t \to +\infty} \|U_t X\|_{L^2(\Hh, \dlt)}$, we get 
\begin{equation}\label{eq14bis3}
 \limsup_{t \to +\infty} \|U_t X_\eps\|_{L^2(\Hh, \dlt)}>0,
\end{equation}
since $U_t$ {preserves} the norm $L^2(\Hh)$. By Theorem \ref{p0bb}, $\sff[U_tX_\eps]\le \Lc(|t|+a_\eps)$. Now, we  apply the second part of the Theorem \ref{t2} to $X_\eps$. Given $a_\eps$ and $b$, we can find $l$ so large that \eqref{eq14bis3} yields
\[
 \limsup_{t \to +\infty} \|U_t X_\eps\|_{L^2(\Hh, \Delta_{l,a_\eps,t})}>0
\]
and, therefore, $\mu \in \sz$ which  gives a contradiction. The case when $t\to -\infty$ can be handled similarly. 
\qed

\begin{Cor}\label{t2bis}
Let $\Hh$ be a proper Hamiltonian. Suppose $\Hh \in \szcs$ and $X$ is any real nonzero element in $H_c$.  Define  $a=\Tc(\sff[X])$. Then, we have 
\begin{equation}\label{eq14bis}
 \liminf_{t \to +\infty} \|(U_t+U_{-t}) X\|_{L^2(\Hh, \dlto)} > 0,\,\,
\end{equation}
for all $\ell>0$.  Conversely, suppose there is real $X \in H_c$ such that 
\begin{eqnarray}\label{eq143bis}
\limsup_{t \to +\infty} \|(U_t+U_{-t})X\|_{L^2(\Hh, \dlto)} > 0,
\end{eqnarray}
for some $\ell > 0$ and, again,   $a=\Tc(\sff[X])$.  Then, $\Hh \in \szcs$. 
\end{Cor}
\beginpf Suppose \eqref{eq143bis} holds.  If $\Hh\notin \szcs$, then 
$$ \lim_{t\to +\infty} \|U_tX\|_{L^2(\Hh, \dlto)}=\lim_{t\to +\infty} \|U_{-t}X\|_{L^2(\Hh, \dlto)}=0$$ by the previous theorem.
Since $\|(U_t+U_{-t})X\|_{L^2(\Hh, \dlto)}\le \|U_tX\|_{L^2(\Hh, \dlto)}+\|U_{-t}X\|_{L^2(\Hh, \dlto)}$, we get a contradiction. Conversely, suppose \begin{equation}\label{ch2}
\liminf_{t \to +\infty} \|(U_t+U_{-t})X\|_{L^2(\Hh, \dlto)} = 0
\end{equation}
for some $\ell > 0$. Again,  consider the function $f = \W_{\Hh} X$ and denote $g=e^{-i(a-\ell)z}f$. As in the proof of Theorem \ref{t2}, we have
$$
\liminf_{t\to +\infty}
\left(\inf_{u \in \B_{r}}
\|e^{itx}f-(u-e^{-itx}f)\|_{L^2(\mu)}\right) = 0.
$$
That gives \eqref{eq16} and the rest of the argument repeats the proof of the Theorem \ref{t2}.
\qed

\medskip

\subsection{Long-time asymptotics of the evolution in Szeg\H{o} case: preliminaries} 

In the rest of the section, we are going to  study the long-time behavior of the  group 
$U_t = e^{it\Di_\Hh}$.  To this end, we need to do some additional work first.
In this subsection, we collect all necessary definitions and auxiliary results. In many places, the presentation follows  \cite{BD2017},  \cite{BD2019} and \cite{B2018}, where one can find more details and references. Let $\Hh$ be a singular Hamiltonian on $\R_+$ and let $\Phi = \phs$ be the solution of the Cauchy problem $J\Phi'(\tau,z) = z\Hh \Phi(\tau,z)$, $\Phi(0, z) = \zo$, $\tau \in \R_+$, $z \in \C$. Recall that $\Theta = \ths$ solves the same differential equation but  satisfies different boundary condition: $\Theta(0, z) = \oz$. The Titchmarsh-Weyl function of any singular Hamiltonian $\Hh$, which is not equal to $\kappa\excl$ a.e. on $\R_+$ for some function $\kappa$, is defined by
\begin{equation}\label{mf}
m(z) = \lim_{\tau \to +\infty}\frac{\Phi^-(\tau,z)}{\Theta^-(\tau,z)}, \qquad z \in \C_+.
\end{equation}
That function is analytic and takes $\C_+$ into its closure $\ov{\C_+}$. The Herglotz representation of $m$ has the form
\begin{equation}\label{hg}
m(z) = \frac{1}{\pi}\int_{\R}\left(\frac{1}{x - z} - \frac{x}{x^2+1}\right) d\mu(x) + b_0z + a_0, \quad z \in\C_+,
\end{equation}
where $\mu$ is called the main spectral measure of $\Hh$, $b_0 \ge 0$, and $a_0 \in \R$. Given proper $\Hh$,   define $\Hh_r$ by $\tau \mapsto \Hh(\tau+r)$ for every $r>0$.  Let $m_r$, $\mu_r$, $b_r$, $a_r$ denote the Titchmarsh-Weyl function of $\Hh_r$, its spectral measure, and the coefficients in the Herglotz representation \eqref{hg} for $m_r$. Define 
\begin{align*}
\I_{\Hh}(r) &= \Im m_r(i) = \frac{1}{\pi}\int_{\R}\frac{d\mu_r(x)}{x^2+1} + b_r, \\
\Rr_{\Hh}(r) &= \Re m_r(i) = a_r,\\
\J_{\Hh}(r) &= \frac{1}{\pi}\int_{\R}\frac{\log w_r(x)}{x^2+1}\,dx, 
\end{align*}  
where $\mu_r = w_r \,dx + \mu_{r,\mathbf{s}}$. It is well-known that $\Rr_{\Hh}$ is identically zero if the Hamiltonian $\Hh$ is diagonal (see, e.g., Lemma 2.2 in \cite{BD2017}). The quantity $\K_\Hh(r) = \log \I_{\Hh}(r) - \J_{\Hh}(r)$, $r \ge 0$ is called the entropy function of $\Hh$. 
Jensen's inequality gives $\K_\Hh(r)\ge 0$ for all $r \ge 0$ and we have $\K_{\Hh}(0) < \infty$ if and only if $\mu \in \sz$. In \cite{B2018}, it was proved that $\mu\in \sz$ implies  $\K_{\Hh}(r) < \infty$ for every $r \in \R_+$. Moreover, the function $\K_{\Hh}$ is absolutely continuous, non-increasing, $\lim_{r\to\infty}\K_{\Hh}=0$, and
\begin{equation}\label{sa4}
\|\K'_{\Hh}\|_{L^1(\R_+)} = \K_{\Hh}(0),
\end{equation}
see Lemma 2.3 and Lemma 2.4 in \cite{B2018}. 
We will need an auxiliary matrix-function $G$ 
\begin{equation}\label{eq6}
G := 
\begin{pmatrix}
1/\sqrt{\I_{\Hh}} & \Rr_{\Hh}/\sqrt{\I_{\Hh}} \\
0 & \sqrt{\I_{\Hh}}
\end{pmatrix},
\qquad 
G^{-1} = 
\begin{pmatrix}
\sqrt{\I_{\Hh}} & -\Rr_{\Hh}/\sqrt{\I_{\Hh}} \\
0 & 1/\sqrt{\I_{\Hh}}
\end{pmatrix}.
\end{equation}
For $\Hh = \sth$, the formula
$$
(G^{-1})^*\Hh G^{-1} = 
\begin{pmatrix}
\I_{\Hh}h_1 & -\Rr_{\Hh}h_1 + h \\
-\Rr_{\Hh}h_1 + h & (\Rr_{\Hh}^2 h_1 - 2\Rr_{\Hh}h + h_2)/\I_{\Hh}
\end{pmatrix}
$$
holds.
It was proved in Lemma 2.4 of \cite{B2018}, that for every $\Hh$ whose spectral measure is in the Szeg\H{o} class $\sz$, the function $\K_{\Hh}$ satisfies
\begin{equation}\label{eq22}
\K_{\Hh}' = 2\sqrt{\det\Hh} - \trace (G^{-1})^*\Hh G^{-1}.
\end{equation}
Recall that $\Theta = \Theta(\tau, z)$ is the solution of the Cauchy problem \eqref{cs}. Define 
$\widetilde\Theta(\tau, z) = G(\tau) \Theta(\tau, z)$, 
$\Et_{\tau}( z) = \widetilde \Theta^+(\tau, z) + i \widetilde\Theta^-(\tau, z)$, and $
\Et^\sharp_{\tau}( z) = \widetilde \Theta^+(\tau, z) - i \widetilde\Theta^-(\tau, z)$.

\begin{Lem}\label{l0}
We have \[
E_\tau(z)\ov{E_\tau(\lambda)} - E_\tau^\sharp(z)\ov{E_\tau^\sharp(\lambda)} = \Et_\tau(z)\ov{\Et_\tau(\lambda)} - \Et_\tau^\sharp(z)\ov{\Et_\tau^\sharp(\lambda)}\]
for all $z, \lambda \in \C$, $\tau \ge 0$.
\end{Lem}
\beginpf Take $z, \lambda \in \C$, $\tau \ge 0$. We have
\begin{align*}
\Et_\tau(z)\ov{\Et_\tau(\lambda)} - \Et_\tau^\sharp(z)\ov{\Et_\tau^\sharp(\lambda)} 
&=\left\langle \begin{pmatrix}
1 & 0 \\ 0 & -1
\end{pmatrix}\begin{pmatrix}
\Et_{\tau}(z) \\ \Et_{\tau}^{\sharp}(z)
\end{pmatrix}, \begin{pmatrix}
\Et_{\tau}(\lambda) \\ \Et_{\tau}^{\sharp}(\lambda)
\end{pmatrix}\right\rangle_{\C^2} 
\\
&=-2i\left\langle 
J
\begin{pmatrix}
\widetilde \Theta^+(\tau, z) \\ \widetilde \Theta^-(\tau, z)
\end{pmatrix},
\begin{pmatrix}
\widetilde \Theta^+(\tau, \lambda) \\ \widetilde \Theta^-(\tau, \lambda)
\end{pmatrix}
\right\rangle_{\C^2}\\
&=-2i\left\langle 
G^*(\tau)JG(\tau)
\begin{pmatrix}
\Theta^+(\tau, z) \\ \Theta^-(\tau, z)
\end{pmatrix},
\begin{pmatrix}
\Theta^+(\tau, \lambda) \\ \Theta^-(\tau, \lambda)
\end{pmatrix}
\right\rangle_{\C^2}.
\end{align*}
For every $\tau \ge 0$, the matrix $G(\tau)$ has real entries and unit determinant which gives $G^*(\tau)JG(\tau) = J$. Thus, we have
\begin{align*}
\Et_\tau(z)\ov{\Et_\tau(\lambda)} - \Et_\tau^\sharp(z)\ov{\Et_\tau^\sharp(\lambda)} 
&=-2i\left\langle 
J
\begin{pmatrix}
\Theta^+(\tau, z) \\ \Theta^-(\tau, z)
\end{pmatrix},
\begin{pmatrix}
\Theta^+(\tau, \lambda) \\ \Theta^-(\tau, \lambda)
\end{pmatrix}
\right\rangle_{\C^2}=E_\tau(z)\ov{E_\tau(\lambda)} - E_\tau^\sharp(z)\ov{E_\tau^\sharp(\lambda)}
\end{align*}
and that proves the statement. \qed\medskip

Since $E_\tau$ is Hermite-Biehler function, taking $z=\lambda$ in the last lemma implies that $\Et_{\tau}$ is  Hermite-Biehler function as well, and hence it has no zeroes in $\C_+$.
Define $\alpha_\tau \in \T$ such that  $\alpha_\tau\Et_{\tau}( i)> 0$ and put
\begin{equation}\label{eq7}
\Pt_{2\tau}(z) = \bar\alpha_\tau e^{i\Tc(\tau) z}\Et^\sharp_{\tau}( z), \qquad \widetilde P_{2\tau}^{*}(z) = \alpha_\tau e^{i\Tc(\tau) z}\Et_{\tau}( z), \qquad z \in \C, \quad \tau \ge 0.
\end{equation}
As in Section 4 of \cite{B2018}, we call $\Pt_{\tau}$ and $\Pt_{\tau}^{*}$ the  regularized Krein's orthogonal entire  functions generated by $\mu$.  To some extent, the  introduction of these functions will allow us to use ideas of the theory of polynomials orthogonal on the unit circle, see, e.g., Lemma \ref{l2} below. Both $\Phi^+,\Phi^-,\Theta^+$, and $\Theta^-$ are entire functions of finite exponential type (see, e.g., Lemma 17 in \cite{Romanov}) so $\Pt_{2\tau}$ and $\Pt_{2\tau}^{*}$ have finite exponential type as well. 
Their basic properties  were studied in the papers  \cite{BD2019}, \cite{B2018} and we discuss some of them now. 
 From the definition, it is immediate that $\widetilde P_{2\tau}^{*}$ satisfies relation $|\widetilde P_{2\tau}^{*}(x)| = |\Et_\tau(x)|$ for $x \in \R$. Therefore, by Theorem 1.3 in \cite{B2018} we have
\begin{equation}\label{eq17}
\lim_{\tau \to \infty}\int_{\R}\frac{|\log|\Pt_{\tau}^{*}(x)|^{-2} - \log w(x)|}{x^2+1}\,dx = 0.  
\end{equation}
The formula \eqref{eq11} and the Lemma \ref{l0} yield
\begin{equation}\label{eq9}
k_{\B_\tau,\lambda}(z) = -\frac{1}{2\pi i} \frac{\Et_{\tau}(z)\ov{\Et_{\tau}(\lambda)} - \Et_{\tau}^\sharp(z)\ov{\Et_{\tau}^\sharp(\lambda)}}{z - \bar \lambda}, \qquad z \in \C. 
\end{equation} 
Consider now the ``shifted'' Hilbert space $e^{i\Tc(\tau) z}\B_{\tau}$. From \eqref{eq9}, we conclude that its reproducing kernel at $\lambda \in \C$ is given by 
\begin{align}
K_{\tau, \lambda}(z) \label{fi1}
&=-\frac{e^{i\Tc(\tau)(z - \bar\lambda)}}{2\pi i} \frac{\Et_{\tau}(z)\ov{\Et_{\tau}(\lambda)} - \Et_{\tau}^\sharp(z)\ov{\Et_{\tau}^\sharp(\lambda)}}{z - \bar \lambda},\\
&= - \frac{1}{2\pi i} \frac{\Pt_{2\tau}^*(z)\ov{\Pt_{2\tau}^*(\lambda)} - \Pt_{2\tau}(z)\ov{\Pt_{2\tau}(\lambda)}}{z - \bar \lambda}. \nonumber
\end{align}
In the case when $\det \Hh=1$, the following result was obtained in Lemma 4.1 of \cite{B2018} where the expression for $\alpha_\tau$ was  found in terms of $\I_{\Hh}$ and $\Rr_{\Hh}$. 
\begin{Lem}\label{l24}
For every $\tau > 0$, 
the function $\widetilde P_{2\tau}^{*}$ is outer in $\C_+$.
\end{Lem}
\beginpf 
 Recall that $\Et_\tau$ is a Hermite-Biehler function. By definition, it can be written as
\begin{eqnarray}\label{dop1}
\Et_\tau=\frac{E_\tau(\Rr_{\Hh}+i(\I_{\Hh}+1))-E_\tau^\sharp(\Rr_{\Hh}+i(\I_{\Hh}-1))}{2i\sqrt{\I_{\Hh}}}\\=
\frac{E_\tau}{2i\sqrt{\I_{\Hh}}}\Bigl(\Rr_{\Hh}+i(\I_{\Hh}+1)\Bigr)\Bigl(1-\frac{E^\sharp_\tau}{E_\tau} \frac{\Rr_{\Hh}+i(\I_{\Hh}-1)}{\Rr_{\Hh}+i(\I_{\Hh}+1)} \Bigr).\nonumber
\end{eqnarray}
The formula \eqref{eq53} says
\[
\type E_\tau=\limsup_{y\to +\infty}\frac{\log|E_\tau(iy)|}{y}=\Tc(\tau).
\]
Since {$\I_{\Hh}> 0$} and $|E^\sharp_\tau/E_\tau|< 1$ in $\C_+$, we get  $|(\Rr_{\Hh}+i(\I_{\Hh}-1))/(\Rr_{\Hh}+i(\I_{\Hh}+1))|<1$ and
\begin{equation}\label{simj}
\limsup_{y\to +\infty}\frac{\log|\Et_\tau(iy)|}{y}=\Tc(\tau).
\end{equation}
The formula \eqref{dop1} shows that $\Et_\tau$ is a linear combination of $E_\tau$ and $E^\sharp_\tau$, two functions of exponential type $\Tc(\tau)$, and so its exponential type is at most $\Tc(\tau)$. Thus, identity \eqref{simj} gives $\type \Et_\tau=\Tc(\tau)$.\smallskip

{Recall} again that $\Et_\tau$ is Hermite-Biehler function. If it has no roots in $\C_-$, then $\Et_\tau=Ce^{-i\Tc(\tau)z}$ with some nonzero constant $C$. Hence, $\widetilde P_{2\tau}^{*}$ is a positive constant and we are done. If $\Et_\tau$ does have a root in $\C_-$, we call it $\lambda$ and argue as follows. Since $\Et_\tau(\lambda)=0$, we also have $\Pt_{2\tau}(\bar\lambda)=0$ by definition. Formula \eqref{fi1} takes the form
\begin{equation}\label{simk}
K_{\tau, \bar\lambda}(z)= - \frac{1}{2\pi i} \frac{\Pt_{2\tau}^*(z)\ov{\Pt_{2\tau}^*(\bar\lambda)}}{z - \lambda}.
\end{equation}
Function $K_{\tau, \bar\lambda}(z)$ is a reproducing kernel of $e^{i\Tc(\tau) z}\B_{\tau}$ at point $\bar\lambda$. Thus, $K_{\tau, \bar\lambda}$ belongs to this space and, by \eqref{eq50}, 
$K_{\tau, \bar\lambda}=e^{i\Tc(\tau) z}E_\tau g$ with some $g\in H_2(\C_+)$. The function $E_\tau$ is of bounded type in $\C_+$. Hence, $K_{\tau, \bar\lambda}$ is also of bounded type there. We know that $
\Pt_{2\tau}^*$ is entire and has no roots in $\C_+$ because $\Et_{\tau}$ is Hermite-Biehler and has no roots there.
Hence,  Smirnov-Nevanlinna factorization of $K_{\tau, \bar\lambda}$ can be written as
$
K_{\tau, \bar\lambda}=\xi e^{icz}\mathcal{O},
$
where $\xi\in \T$ is a constant, $c\in \R$,  and $\mathcal{O}$ is outer. Since
\[
\lim_{y\to +\infty}\frac{\log|\mathcal{O}(iy)|}{y}=0,
\]
we also have
\[
-c=\limsup_{y\to +\infty}\frac{\log|K_{\tau, \bar\lambda}(iy)|}{y}=\limsup_{y\to +\infty}\frac{\log |\Pt_{2\tau}^*(iy)|}{y}=0,
\]
where \eqref{simj}, \eqref{simk}, and the definition of $\Pt_{2\tau}^*$ have been used. Since $c=0$, we get $K_{\tau, \bar\lambda}=\xi \mathcal{O}$ and the formula \eqref{simk} along with normalization $\Pt_{2\tau}^*(i)>0$ prove the lemma.
\qed
\medskip

Given a measure $\mu = w\,dx + \mus$ in $\sz$, we   denote by $D_{\mu}$ its Szeg\H{o} function: 
\begin{equation}\label{eq94}
D_\mu(z) =  \exp\left(\frac{1}{\pi i}\int_{\R}\log \sqrt{w(x)}\left(\frac{1}{x-z} - \frac{x}{x^2+1}\right)dx\right), \qquad z \in \C_+.
\end{equation}
In other words, $D_{\mu}$ is the outer function in $\C_+$ such that $D_{\mu}(i) > 0$ and $|D_{\mu}|^2 = w$ almost everywhere on the real line $\R$ in the sense of non-tangential boundary values.

\medskip

The following lemma will play a key role later on.
\begin{Lem}\label{l2}
Let $\mu = w\,dx + \mus$ be a measure in $\sz$, and let $\Pt_{\tau}$, $\Pt_{\tau}^{*}$ be its regularized Krein's orthogonal entire functions. Then,
\begin{align}
&\lim_{\tau \to \infty}\Pt_{\tau}^{*}(z) = D_{\mu}^{-1}(z), \label{eq18}\\ 
&\lim_{\tau \to \infty} \Pt_{\tau}(z) = 0, \label{eq19}
\end{align}
uniformly on compacts in $\C_+$, and
\begin{align}
&\lim_{\tau \to \infty}\int_{\R}\frac{|\Pt_{\tau}^{*}(x)|^2}{x^2+1}d\mus(x) = 0, \label{eq4}\\
&\lim_{\tau \to \infty}\int_{\R}\frac{|\Pt_{\tau}^{*}(x) - D^{-1}_{\mu}(x)|^2}{x^2+1}w(x)\,dx = 0. \label{eq1}
\end{align}
\end{Lem}
\beginpf  Formula \eqref{eq17} gives \eqref{eq18} after comparing the multiplicative representations for outer functions $\Pt_{\tau}^{*}$ and $D_{\mu}^{-1}$. Let $\B_{\tau}$ be defined by \eqref{eq15}.  The standard variational property of the reproducing kernel yields
\begin{equation}\label{dg0}
{\|K_{\tau, \lambda}\|_{L^2(\mu)} = \sup\{|f(\lambda)|\colon  f \in e^{i\Tc(\tau) z}\B_{\tau},\; \|f\|_{L^2(\mu)} \le 1\}.}
\end{equation}
We claim that $\|K_{\tau, \lambda}\|_{L^2(\mu)}$ is non-decreasing in $\tau \in \R_+$. To prove it, we first notice that the space $e^{ib z}\B_{\Lc(b)}$ coincides with the completion in $L^2(\mu)$ of the set $e^{ib z}\E_{b}$ for every fixed $b>0$ according to Theorem~\ref{nmt2}. It follows that
\begin{equation}\label{dg2}
\|K_{\Lc(b), \lambda}\|_{L^2(\mu)} = \sup\{|f(\lambda)|\colon f \in e^{ib z}\E_{b},\; \|f\|_{L^2(\mu)} \le 1\}.
\end{equation}
Since $e^{ib_1 z}\E_{b_1} \subseteq e^{ib_2 z}\E_{b_2}$ if $b_1 \le b_2$, we have 
\begin{equation}\label{dg1}
\|K_{\Lc(b_1), \lambda}\|_{L^2(\mu)} \le \|K_{\Lc(b_2), \lambda}\|_{L^2(\mu)},
\end{equation}
provided that $0<b_1 \le b_2$.
Now, take  arbitrary positive $\tau > 0$ and let $b = \Tc(\tau)$. We have  
$\B_{\tau} \subseteq \bigcap_{\eps > 0}   \B_{\Lc(b+\eps)} \subseteq \bigcap_{\eps > 0}\clos_{L^2(\mu)}\E_{b+\eps}$ with the last inclusion 
following from Theorem \ref{nmt2}.
 Hence,   $e^{i\Tc(\tau)z}\B_{\tau} \subseteq \bigcap_{\eps > 0}\clos_{L^2(\mu)}(e^{i(b+\eps)z}\E_{b+\eps})$.
From \eqref{dg0}, \eqref{dg2}, and \eqref{dg1}, we get $\|K_{\tau, \lambda}\|_{L^2(\mu)} \le \inf_{\eps > 0}\|K_{\Lc(b+\eps), \lambda}\|_{L^2(\mu)}$. Finally, if we have $\tau_1<\tau_2$ for which $\Tc(\tau_1)=\Tc(\tau_2)=b$, then $\B_{\tau_1}\subsetneq\B_{\tau_2}$ and \eqref{dg0} yields
\[
\|K_{\tau_1, \lambda}\|_{L^2(\mu)} \le \|K_{\tau_2, \lambda}\|_{L^2(\mu)}
\]
in that situation too. Putting together all cases, we get our claim.
Therefore, $\|K_{\tau, \lambda}\|_{L^2(\mu)}$ is non-decreasing in $\tau \in \R_+$. In particular, we have
$$
\lim_{\tau \to \infty} \|K_{\tau, \lambda}\|_{L^2(\mu)} = \lim_{\eta \to \infty} \|K_{\Lc(\eta), \lambda}\|_{L^2(\mu)}.
$$
By the Krein-Wiener theorem (combine formulas $(9.9)$, $(9.13)$, and $(9.14)$ in \cite{Den06}), we have
$$
\lim_{b \to \infty}\|K_{\Lc(b), \lambda}\|_{L^2(\mu)}^{2} = \frac{1}{4\pi} \frac{|D_{\mu}(\lambda)|^{-2}}{\Im\lambda},
$$
where the convergence is uniform on compact sets in $\C_+$. Since $K_{\tau, \lambda}$ is a reproducing kernel, one has 
\begin{equation}\label{eq67}
\|K_{\tau, \lambda}(\cdot)\|_{L^2(\mu)}^{2} = K_{\tau, \lambda}(\lambda) =  \frac{1}{4\pi} \frac{|\Pt_{2\tau}^*(\lambda)|^2 - |\Pt_{2\tau}(\lambda)|^2}{\Im\lambda}
\end{equation}
for every $\tau \ge 0$. It follows that  
\begin{equation}
\lim_{\tau \to \infty} \frac{1}{4\pi} \frac{|\Pt_{2\tau}^{*}(\lambda)|^2 - |\Pt_{2\tau}(\lambda)|^2}{\Im\lambda} =
\lim_{\tau \to \infty} \|K_{\tau, \lambda}\|_{L^2(\mu)}^{2}
=
 \lim_{b \to \infty} \|K_{\Lc(b), \lambda}\|_{L^2(\mu)}^{2}
=
 \frac{1}{4\pi} \frac{|D_{\mu}(\lambda)|^{-2}}{\Im\lambda}\label{fi2}
\end{equation}
holds locally uniformly in $\C_+$.
Combined with \eqref{eq18}, that implies \eqref{eq19}. From \eqref{eq17} and Jensen's inequality, we get
\begin{align}
1 &= \lim_{\tau \to \infty}\exp\left(\frac{1}{\pi}\int_{\R}\frac{\log(|\Pt^*_{\tau}(x)|^2 w(x))}{x^2+1}\,dx\right) \le \limsup_{\tau \to \infty}\frac{1}{\pi}\int_{\R}\frac{|\Pt^*_{\tau}(x)|^2}{x^2+1}\,d\mu(x). \label{eq46}  
\end{align}
On the other hand, one can write
\begin{align}
\left\|\frac{\Pt^*_{\tau}}{x+i}\right\|_{L^2(\mu)}	
&=
\frac{1}{|\Pt_{\tau}^{*}(i)|}\left\|\frac{\Pt_{\tau}^{*}\ov{\Pt_{\tau}^{*}(i)}}{x + i}\right\|_{L^2(\mu)} \notag\\
&\le 
\frac{1}{|\Pt_{\tau}^{*}(i)|}\left\|\frac{\Pt_{\tau}^{*}\ov{\Pt_{\tau}^{*}(i)} - \Pt_{\tau}\ov{\Pt_{\tau}(i)}}{x + i}\right\|_{L^2(\mu)} + \left|\frac{\Pt_{\tau}(i)}{\Pt_{\tau}^{*}(i)}\right|\cdot \left\|\frac{ \Pt_{\tau}}{x + i}\right\|_{L^2(\mu)}  \notag\\
&= 
\frac{1}{|\Pt_{\tau}^{*}(i)|}\left\|\frac{\Pt_{\tau}^{*}\ov{\Pt_{\tau}^{*}(i)} - \Pt_{\tau}\ov{\Pt_{\tau}(i)}}{x + i}\right\|_{L^2(\mu)} + \left|\frac{\Pt_{\tau}(i)}{\Pt_{\tau}^{*}(i)}\right|\cdot \left\|\frac{ \Pt^*_{\tau}}{x + i}\right\|_{L^2(\mu)}.\label{27yan}
\end{align}
Relations \eqref{eq18} and \eqref{eq19} yield $\lim_{\tau \to +\infty}|\Pt_{\tau}(i)/\Pt_{\tau}^{*}(i)| = 0$. {We also have} 
\begin{align*}
\limsup_{\tau \to \infty}\frac{1}{|\Pt_{\tau}^{*}(i)|^2}\left\|\frac{\Pt_{\tau}^{*}\ov{\Pt_{\tau}^{*}(i)} - \Pt_{\tau}\ov{\Pt_{\tau}(i)}}{x + i}\right\|_{L^2(\mu)}^{2} 
&\stackrel{\eqref{fi1}}{=} \limsup_{\tau \to \infty}\frac{4\pi^2}{|\Pt_\tau^*(i)|^2}\|K_{\tau/2, i}\|_{L^2(\mu)}^{2} \\
&\stackrel{\eqref{eq67}}{=} \pi\limsup_{\tau \to \infty}\frac{|\Pt_{\tau}^{*}(i)|^2 - |\Pt_{\tau}(i)|^2}{|\Pt_{\tau}^{*}(i)|^2} = \pi. 
\end{align*}
{That identity, along with  \eqref{27yan}, yields $\Pt_{\tau}^{*}/(x+i) \in L^2(\mu)$ and $\limsup_{\tau\to\infty}\|\Pt_{\tau}^{*}/(x+i)\|^2_{L^2(\mu)}\le\pi$.} Moreover, the inequality in \eqref{eq46} is, in fact, equality, and we get
\begin{equation}\label{sa2}
\lim_{\tau \to \infty}\frac{1}{\pi}\int_{\R}\frac{|\Pt^*_{\tau}(x)|^2}{x^2+1}\,d\mu(x)=1.
\end{equation}
Next, we claim that 
\begin{equation}\label{out2}
\frac{1}{\pi}\int_{\R}\frac{\Pt_\tau^*(x)D_{\mu}(x)}{x^2+1}\,dx = \Pt^*_\tau(i)D_{\mu}(i).
\end{equation}
Given the properties of $\Pt_\tau^*$, this is nearly obvious. However, in the next few lines,  we give the proof of \eqref{out2}.
Indeed, as {showed in Lemma \ref{l24} above}, $\Pt_\tau^*$ is outer in $\C_+$. Then, the function $\Pt_\tau^* D_{\mu}/(z+i)$ lies in $N_+(\C_+)$ and has non-tangential boundary values in $L^2(\R)$ thanks to the following bound
$$
\int_{\R}\left|\frac{\Pt_\tau^*(x) D_{\mu}(x)}{x+i}\right|^2\,dx \le \int_{\R}\frac{|\Pt_\tau^*(x)|^2}{x^2+1}\,d\mu(x) \stackrel{\eqref{sa2}}{<} \infty.
$$
Hence, $\Pt_\tau^* D_{\mu}/(z+i)$ is in $N_+(\C_+) \cap L^2(\R) = H^2(\C_+)$ (see the discussion after Theorem 5.4 in \cite{Garnett} concerning the last equality of sets). Therefore, the function $f$ defined by $$f(\xi) = \Pt_\tau^*(z) D_{\mu}(z), \qquad z \in \C_+,\qquad \xi=\frac{i-z}{i+z}\in \D$$
belongs to the Hardy space $H^2(\D)$ in the open unit disk $\D$ as established in Chapter VI.C in \cite{Koos98}. The mean-value formula for functions in $H^2(\D)$ yields
\[
\frac{1}{2\pi}\int_0^{2\pi} f(e^{i\theta})\,d\theta= f(0).
\]
That gives \eqref{out2}, when written in terms of $z$.

Having proved \eqref{out2}, we can  write
\begin{equation}\label{sa3}
\frac{1}{\pi}\int_{\R}\frac{|\Pt_{\tau}^*(x) - D^{-1}_\mu(x)|^2}{x^2 + 1} |D_\mu(x)|^2 \, dx
=\frac{1}{\pi}\int_{\R}\frac{|\Pt_{\tau}^*(x) D_{\mu}(x)|^2}{x^2 + 1} \, dx + 1 - 2\Re(\Pt_{\tau}^*(i)D_{\mu}(i)).
\end{equation}
Notice now that
\[
\frac{1}{\pi}\int_{\R}\frac{|\Pt_{\tau}^*(x) D_{\mu}(x)|^2}{x^2 + 1} \, dx+\frac{1}{\pi}\int_{\R}\frac{|\Pt_{\tau}^*(x)|^2}{x^2 + 1} \, d\mus(x)=\frac{1}{\pi}\int_{\R}\frac{|\Pt_{\tau}^*(x)|^2}{x^2 + 1} \, d\mu(x)
\]
and the right-hand side converges to $1$ when $\tau\to\infty$ by \eqref{sa2}. Then, 
\[
\lim_{\tau\to\infty}\left(1 - 2\Re(\Pt_{\tau}^*(i)D_{\mu}(i))\right)=-1
\]
by \eqref{eq18}. Thus, \eqref{sa3} yields
\[
0\le \limsup_{\tau\to\infty}\left(\frac{1}{\pi}\int_{\R}\frac{|\Pt_{\tau}^*(x) - D^{-1}_\mu(x)|^2}{x^2 + 1} |D_\mu(x)|^2 \, dx+\frac{1}{\pi}\int_{\R}\frac{|\Pt_{\tau}^*(x)|^2}{x^2 + 1} \, d\mus(x)\right)=0.
\]
Therefore, \eqref{eq4} and \eqref{eq1} follow. \qed

\medskip



The Lemmas \ref{l4} and \ref{l5} below are not new. We give their proofs for the reader's convenience.

\begin{Lem}\label{l4}
Let $\Hh$ be a proper Hamiltonian on $\R_+$. Then, the set $\dom \Di_{\Hh} \cap H_c$ is dense in $H$.
\end{Lem}
\beginpf Consider the linear manifold of functions $X\in H$ for which there is $L \in \R_{+} \setminus \cup_{I \in \I(\Hh)}I$ such that $X$ can be written as follows
$$
X(\tau) = 
\begin{cases}J\int_{\tau}^{L}\Hh(s) Y(s)\,ds, & \tau \le L\\
0, & \tau > L 
\end{cases},
\qquad 
Y \in H: \quad \supp Y \subseteq [0, L], \quad ( Y, \oz )_{L^2(\Hh)} = 0.
$$
In the proof of Theorem 3 in \cite{Romanov}, it was showed that if $X_0 \in H$ is orthogonal to  that linear manifold, then $\Hh X_0 = 0$ almost everywhere on $\R_+$ (we apply Theorem 3 of \cite{Romanov} to Hamiltonians  $\Hh$ that do not coincide with those that are equal to $\kappa\excl$ on some interval $[0, \eps]$, since we  study only such Hamiltonians in our paper). That implies this manifold is  dense in $H$. On the other hand,  every $X$ in that manifold has compact support  and 
$$
\langle X(0), \zo \rangle_{\C^2} = \langle J X(0), J\zo \rangle_{\C^2} = -(Y, \oz )_{L^2(\Hh)} = 0.
$$
Therefore, $X \in \dom\Di_{\Hh} \cap H_c$. The lemma follows. \qed

\medskip

Let $\Hh$ be a proper Hamiltonian on $\R_+$, and let $H$ be the Hilbert space generated by $\Hh$. On functions $X \in L_c^2(\Hh)$, define 
\begin{equation}\label{wtm}
\widetilde\W_{\Hh}\colon X \mapsto \frac{1}{\sqrt{\pi}}\int_{0}^{\infty}\langle\Hh(\tau)X(\tau), \Theta(\tau, \bar z)\rangle_{\C^2}\,d\tau, \qquad z \in \C.  
\end{equation}
Note that $\widetilde\W_{\Hh}$ coincides with $\W_{\Hh}$ on $H$. Denote by $P_H$ the orthogonal projector in $L^2(\Hh)$ to $H$. The orthogonal complement $L^2(\Hh) \ominus H$ consists of functions $X$ that satisfy the following conditions:
 $X = 0$ on $\R_+ \setminus \bigcup_{I \in \mathfrak{I}(\Hh)} I$ and 
$$
\quad \int_{I}\kappa(\tau)\langle X(\tau), e_I\rangle_{\C^2}\,d\tau = 0, \qquad \Hh(\tau) = \kappa(\tau)\langle\cdot, e_I\rangle_{\C^2}e_I, \qquad \tau \in I, \qquad I \in \mathfrak{I}(\Hh),
$$
where $e_I\in \R^2, \|e_I\|_{\R^2}=1$, and $\kappa>0$ a.e. on each $I$.
It follows that $P_H\colon L^2(\Hh) \to H$ coincides with the operator 
$$
X \mapsto \begin{cases}
X(\tau), &\tau \in \R_+ \setminus \bigcup_{I \in \mathfrak{I}(\Hh)} I, \\ 
\sum_{I \in \mathfrak{I}(\Hh)} c_I \chi_{I}e_{I}, & \tau \in \bigcup_{I \in \mathfrak{I}(\Hh)} I,
\end{cases} \qquad c_I = \frac{\int_{I}\kappa(\tau)\langle X(\tau), e_I\rangle_{\C^2}\,d\tau}{\int_{I}\kappa(\tau)\,d\tau }.
$$
Indeed, this operator is linear, vanishes on $L^2(\Hh) \ominus H$, and acts as an identity on $H$ because $\chi_I (e_I + e_I^\bot) = \chi_I e_I$ in $L^2(\Hh)$ for every vector $e_I^\bot \in \C^2$ orthogonal to $e_I$. As a consequence, if $r \in \R_+ \setminus \bigcup_{I \in \mathfrak{I}(\Hh)} I$ and $\supp X \subseteq[0, r]$, then $\supp P_H X \subseteq [0, r]$. We use this observation in the formula \eqref{eq02feb} below.   
\begin{Lem}\label{l5}
Let $\Hh$ be a proper Hamiltonian on $\R_+$ and let $\mu$ be its spectral measure. We have $\W_{\Hh}^{-1}\widetilde \W_{\Hh} = P_H$ and $\|\widetilde \W_{\Hh}X\|_{L^2(\mu)} \le \|X\|_{L^2(\Hh)}$ for every $X \in L^2(\Hh)$.
\end{Lem}
\beginpf Consider $X\in L^2(\Hh)$ such that $\supp X \subseteq [0, r]$, where $r$ is not in the interior of an indivisible interval (that is, $r \in \R_+ \setminus \bigcup_{I \in \mathfrak{I}(\Hh)} I)$. Then, taking any
 $z \in \C$, we have
$$
(\widetilde \W_{\Hh}X)(z) = (X, \Theta(\cdot, z))_{L^2(\Hh)} = (X, \chi_{[0, r]}\Theta(\cdot,z))_{L^2(\Hh)}. 
$$
If $I$ is an indivisible interval, we have $\Hh(\tau)=\kappa(\tau) \langle \cdot, e_I\rangle_{\C^2}e_I$ with some vector  $e_I\in \R^2, \, \|e_I\|_{\R^2}=1$ for $\tau\in I$. Equation 
$
J\Theta'(\tau,z)=z\kappa(\tau) \langle \Theta(\tau,z), e_I\rangle_{\C^2}e_I, \,\tau\in I
$
implies $
\Theta'(\tau,z)=-z\kappa(\tau) \langle \Theta(\tau,z), e_I\rangle_{\C^2}Je_I,$ $\tau\in I\,
$. Since $\langle Je_I,e_I\rangle_{\C^2}=0$, we have $\langle\Theta'(\tau,z),e_I\rangle_{\C^2}=0$, and hence $\langle\Theta(\tau,z),e_I\rangle_{\C^2}$ is constant in $\tau$ on $I$.
That gives $\chi_{[0, r]}\Theta(\cdot, z)\in H$ for every $z$ in a sense that $\Theta(\tau,z)$ is constant on each $I$ when considered as an element of $L^2(\Hh)$ defined in \eqref{eq31}. Thus, we have 
\begin{equation}\label{eq02feb}
(X, \chi_{[0, r]}\Theta(\cdot, z))_{L^2(\Hh)} = (P_H X, \chi_{[0, r]}\Theta(\cdot, z))_{L^2(\Hh)} = (\W_{\Hh}P_HX)(z).
\end{equation}
That gives $\|\widetilde \W_{\Hh}X\|_{L^2(\mu)} \le \|X\|_{L^2(\Hh)}$ and $\W_{\Hh}^{-1}\widetilde \W_{\Hh}X = P_HX$. The set of $X$ we considered is dense in $L^2(\Hh)$ and the operator $\W_{\Hh}$ is unitary. Therefore,  the  lemma is true  for all $X\in L^2(\Hh)$.\qed

\subsection{Long-time asymptotics of the evolution in Szeg\H{o} case: the main result } 

{Recall that we study the evolution $U_tX$ when $X\in H$ and $t\to\pm\infty$. We will describe $U_t X$ in terms of the ``free'' evolution $U^0_t Y_{X,\pm}$ of some states $Y_{X,\pm} \in L^2(\Hh_0)$ as $t \to \pm\infty$, where $U_t^0$ is generated by ``free'' Hamiltonian $\Hh_0 = \idm$ on $\R_+$. Note that $U_t$, $U^0_t$ act on different Hilbert spaces and an identification is needed to relate the ``perturbed'' and ``free'' dynamics governed by $U_t$ and $U^0_t$, respectively. First, we observe that given a pair of real states $X \in H_c$, $Y \in L_c^2(\Hh_0)$ such that $\Tc(\sff[X]) = a$ and $\sff[Y] = a$ for some $a>0$, the Theorem \ref{p0bb} yields
$$
\sff[U_t X] = \Lc(|t|+a) = \sff[(U_{t}^0 Y)(\Tc(\cdot))], \qquad t \in \R\backslash \{0\}.
$$
Thus, when $t$ varies in the interval $(0,t_0)$, the wavefronts of $U_t X$ and $(U^{0}_tY)(\Tc(\cdot))$ simultaneously propagate from $\Lc(a)$ to $\Lc(t_0+a)$. That provides an intuition how to map  $U_{t}^0 Y$ into $L^2(\Hh)$. First, we  introduce the  non-negative matrix-function $\Hh_{\mathfrak{n}}$:
$$
\Hh_{\mathfrak{n}}(\tau) 
= 
\begin{cases}
\Bigl(\det\Hh(\tau)\Bigr)^{-\frac 12}\Hh(\tau), & \tau: \det\Hh(\tau) > 0, \\
\infty, & \tau: \det\Hh(\tau) = 0.
\end{cases}
$$
Then, 
\begin{equation}\label{refe7}
\Hh_{\mathfrak{n}}^{-\frac 12}(\tau) 
= 
\begin{cases}
{\Bigl(\det\Hh(\tau)\Bigr)}^{\frac14}\Hh^{-\frac 12}(\tau), & \tau: \det\Hh(\tau) > 0, \\
0, & \tau: \det\Hh(\tau) = 0.
\end{cases}
\end{equation}
It is instructive to note that $\det\Hh_{\mathfrak{n}}(\tau) = 1$ for every $\tau$ that satisfies $\det\Hh(\tau) > 0$. 
Second, we fix a measurable function $\gamma\colon \R_+ \to \T$. Finally, for every $t \in \R$ and $Y \in L^2(\Hh_0)$, define
\begin{equation}\label{eq12}
\widetilde U_{\gamma, t}^0 Y\colon \tau \mapsto 
\begin{cases}
\gamma(\tau)\Hh_{\mathfrak{n}}^{-\frac 12}(\tau)(U_{t}^0 Y)(\Tc(\tau)), \qquad t \ge 0,\\
\overline{\gamma(\tau)}\Hh_{\mathfrak{n}}^{-\frac 12}(\tau)(U_{t}^0 Y)(\Tc(\tau)), \qquad t<0.
\end{cases}
\end{equation}
The role of the ``phase function'' $\gamma$ will become clear in Theorem \ref{t1} below. Given definition \eqref{eq12}, we get several important properties of the evolution $\widetilde U_{\gamma, t}^0$:
\begin{itemize}
\item[$(A)$] the dynamics $\widetilde U_{\gamma, t}^0 Y$ has an explicit expression in terms of $Y$, $\Hh$, $\gamma$;
\end{itemize}
Indeed, that follows from the explicit formula for  $U_{t}^0 Y$ which we obtain in  Lemma \ref{l1} below.
\begin{itemize}
\item[$(B)$] if $\sff[X] = \sff[Y(\Tc(\cdot))]$ for some real states $X \in H_c$, $Y \in L_c^2(\Hh_0)$, then
$$
\sff[\widetilde U_{\gamma, t}^0 Y] = \sff[U_t X], \qquad t \in \R\backslash \{0\}.
$$
\end{itemize}
That is the direct consequence of Theorem \ref{p0bb}.
\begin{itemize}
\item[$(C)$] The map $\widetilde U_{\gamma, t}^0$ is an isometry: 
$$
\|\widetilde U_{\gamma,t}^0 Y\|_{L^2(\Hh)}^{2} = \|Y\|_{L^2(\Hh_0)}^{2}. 
$$
\end{itemize}
The last relation comes from a change of variables:  
$$
\|\widetilde U_{\gamma,t}^0 Y\|_{L^2(\Hh)}^{2} = \int_{\R_+}\Tc'(\tau)\bigl\langle(U_t^0 Y)(\Tc(\tau)),(U_t^0 Y)(\Tc(\tau))\bigr\rangle_{\C^2}\,d\tau = \|U_{t}^0 Y\|_{L^2(\Hh_0)}^{2} = \|Y\|_{L^2(\Hh_0)}^{2}, 
$$
where $\Tc' = \sqrt{\det\Hh}$ is the derivative of the locally absolutely continuous function $\Tc$.
\begin{itemize}
\item[$(D)$] The map $\widetilde U_{\gamma, t}^0$ sends $L^2(\Hh_0)$ into $H$. 
\end{itemize}
Indeed, for every $Y \in L^2(\Hh_0)$ we have $\widetilde U_{\gamma, t}^0 Y = 0$ on each indivisible interval. Therefore, $\Ran \widetilde U_{\gamma, t}^0 \subseteq H$. 

\medskip

We aim to prove the following result.

\begin{Thm}\label{t1}
Let $\Hh \in \szcs$. Then, there exists a function $\gamma\colon \R_+ \to \T$ such that the following assertion holds. For every $X \in H$, there are unique $Y_{X,\pm} \in L^2(\Hh_0)$ such that for every  $b> 0$, we have 
\begin{equation}\label{eq0}
\lim_{t \to \pm\infty}\|U_t X - \widetilde U^0_{\gamma,t} Y_{X,\pm}\|_{L^2(\Hh, \dlt)} =0, 
\end{equation}
where $
\dlt = \bigl[\Lc(|t| -b), \Lc(|t| + b)\bigr]$.
These $Y_{X,\pm}$ can be computed by the formulas 
\begin{equation}\label{sa9}
Y_{X, +}=\W_{\Hh_{0}}^{-1}(f_{\bf ac} \ov{D_{\mu}}), \quad Y_{X, -}=\W_{\Hh_{0}}^{-1}(f_{\bf ac} {D_{\mu}}), 
\end{equation}
where $f=\W_{\Hh}X, f_{\bf ac}:=f\cdot \chi_{\Omega_{\bf ac}(\mu)}$.  Moreover, if $\Hh$ is diagonal, then one can take $\gamma = 1$ on $\R_+$.  
\end{Thm} 

\begin{Rema}
In Theorem \ref{t1}, we do not assume that $X$ has compact support or belongs to the absolutely continuous subspace of $\Di_{\Hh}$. 
\end{Rema}
We start with providing an explicit formula for the evolution $U_t^0$. Let $\Di_{\Hh_0}$ be the self-adjoint operator on $L^2(\Hh_0) = L^2(\R_+,\C^2)$ corresponding to $\Hh_0$. We have $U_t^0= e^{it\Di_{\Hh_0}}$. The main spectral measure of $\Di_{\Hh_0}$ is equal to the  Lebesgue measure on $\R$ and 
\[
(\W_{\Hh_0}Y)(z)=\frac{1}{\sqrt{\pi}}\int_0^\infty \Bigl\langle
\left( \begin{smallmatrix}
y_1(\tau)\\
y_2(\tau)\\
\end{smallmatrix}\right),
\left( \begin{smallmatrix}
\cos (\tau \bar z)\\
-\sin (\tau \bar z)\\
\end{smallmatrix}\right)
\Bigr\rangle_{\C^2}d\tau, \qquad z \in \C,
 \]
for every function $Y=\begin{pmatrix}y_1 \\ y_2\end{pmatrix}$  in $L_c^2(\R_+,\C^2)=L_c^2(\Hh_0)$. Recall that $\widehat h$ denotes the Fourier transform of a function $h\in L^2(\R)$.

\begin{Lem}\label{l1}
Let $h \in  L^2(\R)$ and $Y =\begin{pmatrix}y_1 \\ y_2\end{pmatrix}\in L^2(\Hh_0)$ be defined by $\W_{\Hh_0}Y = h$. Extend $y_1$ to all of $\R$ as an even function and $y_2$ as an odd function. Then, for every $t \in \R$ and $\tau\in \R_+$, we have
\begin{align}
(U_{t}^{0}Y)(\tau) 
&=\frac{1}{\sqrt 2} \widehat h(\tau-t) \begin{pmatrix}1\\-i\end{pmatrix}+\frac{1}{\sqrt 2} \widehat h(-(\tau+t)) \begin{pmatrix}1\\i\end{pmatrix} \label{sa10}\\
&= 
\frac{1}{2}
\begin{pmatrix}
y_1(\tau-t) + y_1(\tau+t)\\
-i(y_1(\tau-t) - y_1(\tau+t))
\end{pmatrix}
+
\frac{1}{2}
\begin{pmatrix}
i(y_2(\tau-t) - y_2(\tau+t))\\
y_2(\tau-t) + y_2(\tau+t)
\end{pmatrix}, \label{sa12}
\end{align}
where the integrals are understood in $L^2(\R)$-sense. In particular, we have
\begin{equation}\label{sa11}
(U_{t}^{0}Y)(\tau)
=\frac{1}{\sqrt 2} \widehat h(\tau-t)\cdot \omi+o(1),\quad   (U_{-t}^{0}Y)(\tau)
=\frac{1}{\sqrt 2} \widehat h(-\tau+t)\cdot \oi+o(1),   \\
\end{equation}
when $t\to +\infty$ and $o(1)$ is understood in $L^2(\Hh_0)$-sense.
\end{Lem}
\beginpf 
We claim that 
\begin{equation}
\label{eq25}
(U_{t}^{0}Y)(\tau)
=
\frac{1}{2\sqrt{\pi}}\int_{\R}e^{i(t - \tau)x}h(x)\,dx \cdot \omi
+
\frac{1}{2\sqrt{\pi}}\int_{\R}e^{i(t + \tau)x}h(x)\,dx \cdot \oi,
\end{equation}
where the integrals are understood in $L^2(\R)$-sense. To prove \eqref{eq25}, we first assume that $h\in L^1(\R)\cap L^2(\R)$. Denote the right-hand side of \eqref{eq25} by $Z_t$. Notice that the integrals in the right-hand side of \eqref{eq25} converge absolutely. We only need to check that the images of $U^0_t Y$ and $Z_t$ under $\W_{\Hh_0}$  coincide. Indeed, 
\begin{align*}
\W_{\Hh_0}(U_t^0 Y)(s) 
=& e^{its}h(s) = \frac{1}{2\pi}\lim_{r \to \infty}\int_{-r}^{r}e^{i\tau s} \int_{\R}e^{i(t - \tau)x}h(x)\,dx\,d\tau\\
=&\frac{1}{2\pi}\lim_{r \to \infty}\left(\int_{0}^{r}\Bigl\langle\omi,\left(\!
\begin{smallmatrix}
\cos (\tau s)\\
-\sin (\tau s)\\
\end{smallmatrix}
\!\right)\Bigr\rangle_{\C^2} \int_{\R}e^{i(t - \tau)x}h(x)\,dx\,d\tau\right)\\
&+\frac{1}{2\pi}\lim_{r \to \infty}\left(\int_{0}^{r}\Bigl\langle\oi,\left(\!
\begin{smallmatrix}
\cos (\tau s)\\
-\sin (\tau s)\\
\end{smallmatrix}
\!\right)\Bigr\rangle_{\C^2} \int_{\R}e^{i(t + \tau)x}h(x)\,dx\,d\tau\right)\\
=&(\W_{\Hh_0}Z_t)(s),
\end{align*}
for almost all $s \in \R$ where the limits are understood in the $L^2(\R)$-sense. Since $L^1(\R)\cap L^2(\R)$ is dense in $L^2(\R)$, we can extend \eqref{eq25} to all of $L^2(\R)$ by continuity. The formula \eqref{sa10} is immediate from \eqref{eq25} if we use the notation for Fourier transform. Finally, taking $t=0$ in \eqref{sa10} gives $y_{1}(\tau) = (\widehat{h}(\tau) + \widehat{h}(-\tau))/\sqrt{2}$, $y_{2}(\tau) = -i(\widehat{h}(\tau) - \widehat{h}(-\tau))/\sqrt{2}$ and  the formula \eqref{sa12} follows. Since $\lim_{t\to +\infty}\|\widehat h(-(\tau+t))\|_{L^2(\R_+)}=0$, we get the first limit in \eqref{sa11}. The second one can be proved similarly.
\qed

\medskip

The free dynamics $U_t^0$ is known to be reducible to the shift operator on the real line. We recall that construction in Section \ref{app2}.

\medskip

\noindent{\bf Proof of Theorem \ref{t1}.} {\it Existence of $Y_{X,\pm}$.}\smallskip

From Lemma \ref{l4}, we know that the set  $\dom \Di_{\Hh} \cap H_c$ is dense in $H$. We start the proof by  considering $X \in \dom \Di_{\Hh} \cap H_c$. Recall that $b> 0$ is a number and $\chi_{\dlt}$ denotes the characteristic function of the interval $\dlt.$ Put $f = \W_{\Hh}X$ and let $f_{\bf ac}=f\cdot \chi_{\Omega_{\bf ac}(\mu)},f_{\bf s}=f\cdot \chi_{\Omega_{\bf s}(\mu)}$.  By the spectral theorem, for every $t \in \R$ and $\lambda \in \C$ we have
\begin{align}
\W_{\Hh}(\chi_{\dlt}U_t X)(\lambda) 
&=\frac{1}{\sqrt{\pi}}\int_{\dlt}\langle\Hh(\tau)U_t X(\tau), \Theta(\tau,  \bar\lambda)\rangle_{\C^2}\,d\tau\notag\\
&=\frac{1}{\sqrt{\pi}}(U_t X, \chi_{\dlt}\Theta(\cdot, \bar \lambda))_{L^2(\Hh)}\notag\\
&=(e^{itx}f, k_{\mu, \dlt, \lambda})_{L^2(\mu)},\label{eq23}
\end{align}
where
\begin{equation}\label{eq24}
k_{\mu, \dlt, \lambda}(x) = 
\frac{1}{\sqrt{\pi}}\W_{\Hh}(\chi_{\dlt}\Theta(\cdot, \bar \lambda)) = \frac{1}{\pi}\int_{\dlt}\langle\Hh(\tau)\Theta(\tau,\bar\lambda), \Theta(\tau,x)\rangle_{\C^2}\,d\tau, \qquad x \in \R.
\end{equation}
We are going to study the asymptotic behavior of \eqref{eq23} when $t \to \pm\infty$ using representation \eqref{eq24}. Let $G$ be the matrix function from \eqref{eq6}. Fix $\tau \in \dlt$ and set
\begin{equation}
G_1 = \left(\begin{smallmatrix}1 & 1\\-i&i \end{smallmatrix}\right), \quad\Psi_1 = \frac{1}{2} G^{-1}G_1\oz, \quad\Psi_2 = \frac{1}{2} G^{-1}G_1\zo. \label{sl2}
\end{equation} 
Recall that $\Pt_\tau$, $\Pt^*_\tau$, $\alpha_\tau$ are defined in \eqref{eq7}. On the real line $\R$, we have 
\begin{align}
\Theta 
&= G^{-1}\begin{pmatrix}\frac{\Et_\tau + \Et^\sharp_\tau}{2}\\ \frac{\Et_\tau - \Et^\sharp_\tau}{2i}\end{pmatrix} = \frac{1}{2}G^{-1}G_1\begin{pmatrix}\Et_\tau\notag\\ \Et^\sharp_\tau\end{pmatrix}
=\frac{e^{-i\Tc(\tau) x}}{2}G^{-1}G_1\begin{pmatrix}\overline{\alpha}_\tau\Pt_{2\tau}^{*}\\ \alpha_\tau\Pt_{2\tau}\end{pmatrix}, \notag\\
&=\overline{\alpha}_\tau e^{-i\Tc(\tau) x}\Pt_{2\tau}^{*} \Psi_1  + {\alpha_\tau}e^{i\Tc(\tau) x}\ov{\Pt_{2\tau}^{*}}\Psi_2. \label{eq47}
\end{align}\smallskip
We continue by getting the estimates on $\|\Psi_{1}\|_{L^2(\Hh, \dlt)}$ and $\|\Psi_{2}\|_{L^2(\Hh, \dlt)}$.
Note that
\begin{align*}
4\langle\Hh\Psi_{1}, \Psi_1\rangle_{\C^2} + 4\langle\Hh\Psi_{2}, \Psi_2\rangle_{\C^2} 
&= \langle \Hh G^{-1}G_1\oz, G^{-1}G_1\oz \rangle_{\C^2}
+ \langle \Hh G^{-1}G_1\zo, G^{-1}G_1\zo \rangle_{\C^2}\\
&= \trace G_1^* (G^{-1})^*\Hh G^{-1}G_1 = 
2\trace (G^{-1})^*\Hh G^{-1}
\end{align*}
due to the fact that $G_1G_1^* = 2\idm$. 
Using this calculation and relation \eqref{eq22}, we get
\begin{align*}
2\|\Psi_{1}\|_{L^2(\Hh,\dlt)}^{2} + 2\|\Psi_{2}\|_{L^2(\Hh,\dlt)}^{2} &= \int_{\dlt}\trace{(G^{-1})^* \Hh G^{-1}}\,d\tau\\
&= \|\K'_{\Hh}\|_{L^1(\dlt)} + 2\int_{\dlt}\sqrt{\det \Hh}\,d\tau,
\end{align*}
because $\K_{\Hh}$ is a non-increasing function.
Note that if $|t|  \ge b$, we have $\int_{\dlt}\sqrt{\det \Hh}\,d\tau = \Tc(\Lc(|t| + b)) - \Tc(\Lc(|t| -b)) = (|t| + b) - (|t| -b) = 2b$.
Thus, for such $t$ we have
\begin{equation}
2\|\Psi_{1}\|_{L^2(\Hh,\dlt)}^{2} + 2\|\Psi_{2}\|_{L^2(\Hh,\dlt)}^{2} \leqslant \|\K'_{\Hh}\|_{L^1(\R_+)} + 2b. \label{eq26}
\end{equation}
Next, we study the inner product \eqref{eq23} using \eqref{eq24}. We need some auxiliary bounds first. By applying the spectral theorem, we have $(x+i)f \in L^2(\mu)$. Relations \eqref{eq4} and \eqref{eq1}, along with
 Cauchy-Schwarz inequality, imply that
\[
 \left(\int_{\R}|f_{\bf s}(x)\Pt_{2\tau}^{*}(x)| d\mus\right)^2\le \left(\int_{\R}|f(x)|^2(x^2+1)\, d\mu\right) \left(\int_{\R}\frac{|\Pt_{2\tau}^{*}(x)|^2}{x^2+1}\,d\mus\right)\to 0,
\]
and that
\begin{eqnarray*}
\left(\int_{\R}|f_{\bf ac}(x)(\Pt_{2\tau}^{*}(x)w(x) - \ \ov{D_{\mu}}(x))|\,dx\right)^2\le \\ \left(\int_{\R}|f(x)|^2(x^2+1) d\mu\right) \left(\int_{\R}\frac{| \Pt_{2\tau}^{*}(x)-D_\mu^{-1}(x)|^2}{x^2+1}w(x)\,dx\right)    \to 0,
\end{eqnarray*}
}
when $\tau\to \infty$.  Since $\Hh\in \szcs$, we get $\lim_{s\to\infty} \Lc(s)=+\infty$. Thus, recalling that $\mu = |D_{\mu}|^2\,dx + \mus$, one has
\begin{equation}
\label{eq36}
\lim_{t\to\infty} \left(\int_\R e^{i(t-\Tc)x}f(x)\Pt_{2\tau}^{*}(x)\,d\mu(x)-\int_\R e^{i(t-\Tc)x}f_{\bf ac}(x)\ov{D_{\mu}(x)}\,dx\right)=0
\end{equation}
uniformly in $\tau\in \dlt$. Similarly,  $(x+i)f \in L^2(\mu)$, relations \eqref{eq4}, \eqref{eq1} and the Riemann-Lebesgue lemma imply 
\begin{equation}\label{eq37}
\lim_{t \to +\infty} \int_{\R}e^{i(t + \Tc(\tau))x}f(x) \overline{\Pt_{2\tau}^{*}(x)}\,d\mu(x) =  
\lim_{t \to +\infty} \int_{\R}e^{i(t + \Tc(\tau))x}f_{\bf ac}(x) \ov{D_{\mu}(x)} \,dx = 0,
\end{equation}
uniformly with respect to $\tau \in \dlt$.  That follows from the  inclusion $f_{\bf ac} \ov{D_{\mu}}\in L^1(\R)$ which is immediate from our assumptions on $f$ and the bound
$$
\left(\int_{\R}|f_{\bf ac}(x)D_{\mu}(x)|\,dx\right)^2 \le  \left(\int_{\R} \frac{dx}{x^2+1}\right) \cdot \left(\int_{\R} |f(x)|(x^2+1)\,d\mu(x)\right) < \infty.
$$
Taking into account \eqref{eq23} and \eqref{eq47}, we see that $\W_{\Hh}(\chi_{\dlt}U_t X)(\lambda)$ equals
\begin{align*} 
(e^{itx}f, k_{\mu, \dlt, \lambda})_{L^2(\mu)} 
&=
\frac{1}{\pi}\int_{\R}e^{itx}f(x)\overline{\int_{\dlt}\langle\Hh(\tau)\Theta(\tau,\bar\lambda), \Theta(\tau,x)\rangle_{\C^2}\,d\tau}\,d\mu(x)\\
&=
\frac{1}{\pi}\int_{\R}\int_{\dlt}e^{itx}f(x)\langle\Hh(\tau)\Theta(\tau,x), \Theta(\tau,\bar\lambda)\rangle_{\C^2}\,d\tau\,d\mu(x)\\
&=
\frac{1}{\pi}\int_{\R}\int_{\dlt}e^{itx}f(x)\langle\Hh(\tau)[\ov\alpha_\tau e^{-i\Tc(\tau) x}\Pt_{2\tau}^{*} \Psi_1], \Theta(\tau,\bar\lambda)\rangle_{\C^2}\,d\tau\,d\mu(x)\\
&\phantom{=}+\frac{1}{\pi}\int_{\R}\int_{\dlt}e^{itx}f(x)\langle\Hh(\tau)[\alpha_\tau e^{i\Tc(\tau) x}\ov{\Pt_{2\tau}^{*}}\Psi_2], \Theta(\tau,\bar\lambda)\rangle_{\C^2}\,d\tau\,d\mu(x)\\
&=\frac{1}{\pi}\int_{\dlt}\left(\ov\alpha_\tau\int_{\R}e^{i(t - \Tc(\tau))x}f_{\bf ac}(x) \ov{D_{\mu}}(x)\,dx\right)\left\langle \Hh(\tau)\Psi_1,\Theta(\tau, \bar \lambda)\right\rangle d\tau + R(t,\lambda),
\end{align*}
where $R(t,\lambda) = R_1(t,\lambda) + R_2(t, \lambda)$,
$$
R_1(t,\lambda) = \frac{1}{\pi}\int_{\dlt}\ov\alpha_\tau\left(\int_{\R}e^{i(t - \Tc(\tau))x}f\Pt_{2\tau}^{*}\,d\mu(x) - \int_{\R}e^{i(t - \Tc(\tau))x}f_{\bf ac} \ov{D_{\mu}}\,dx\right)\left\langle \Hh(\tau)\Psi_1,\Theta(\tau, \bar \lambda)\right\rangle d\tau,
$$
and
$$
R_2(t,\lambda) = \frac{1}{\pi}\int_{\dlt}\alpha_\tau\left(\int_{\R}e^{i(t + \Tc(\tau))x}f(x) \overline{\Pt_{2\tau}^{*}}(x)\,d\mu(x)\right) \left\langle \Hh(\tau)\Psi_2,\Theta(\tau, \bar \lambda)\right\rangle d\tau.
$$
Observe that by \eqref{eq36} and \eqref{eq37} both $R_1(t,\lambda)$ and $R_2(t,\lambda)$ can be represented in the form
$$
R_{1,2}(t,\lambda) = \frac{1}{\sqrt{\pi}}\int_{\dlt}\left\langle \Hh(\tau)\psi_{1,2}(t,\tau),\Theta(\tau, \bar \lambda)\right\rangle d\tau = \widetilde\W_{\Hh}(\Delta_{b,t}\psi_{1,2})(\lambda),
$$
where $\psi_{1,2}(t, \cdot) = \kappa_{1,2}(t, \cdot) \Psi_{1,2}$ for some functions $\kappa_{1,2}(t, \cdot)$ such that $\lim_{t \to +\infty}\|\kappa_{1,2}(t, \cdot)\|_{L^\infty(\dlt)} = 0$. Estimate \eqref{eq26} shows that the quantities  $\|\Psi_{1,2}\|_{L^2(\Hh, \dlt)}$ are uniformly bounded with respect to $t \in \R_+$, hence $\lim_{t \to +\infty}\|\psi_{1,2}(t, \cdot)\|_{L^2(\Hh, \dlt)} =0$.
Therefore, we have $\lim_{t \to +\infty} \|R(t, \cdot)\|_{L^2(\mu)} = 0$ by Lemma \ref{l5}. Summarizing, we see that
$$
\W_{\Hh}(\chi_{\dlt}U_t X)
=\widetilde\W_{\Hh}\left(\frac{\ov\alpha_\tau \chi_{\dlt}}{2\sqrt{\pi}}\int_{\R}e^{i(t - \Tc(\tau))x}f_{\bf ac}(x) \ov{D_{\mu}(x)}\,dx \cdot  G^{-1}(\tau)\omi\right) + o(1),
$$
as $t \to +\infty$ with $o(1)$ in $L^2(\mu)$ (the function under $\widetilde \W_{\Hh}$ depends on $\tau$).
Applying $\W_{\Hh}^{-1}$ and using Lemma \ref{l5}, we get  
\begin{equation}\label{eq27}
(\chi_{\dlt}U_t X)(\tau) 
=P_{H}\left(\frac{\ov\alpha_\tau\chi_{\dlt}}{2\sqrt{\pi}}\int_{\R}e^{i(t - \Tc(\tau))x}f_{\bf ac}(x) \ov{D_{\mu}(x)}\,dx \cdot G^{-1}(\tau)\omi\right) + o(1),
\end{equation}
with $o(1)$ in $H$. Similar reasoning gives
$$
(\chi_{\dlt}U_t X)(\tau) =P_{H}\left(\frac{\alpha_\tau\chi_{\dlt}}{2\sqrt{\pi}}\int_{\R}e^{i(t + \Tc(\tau))x}f_{\bf ac}(x) D_{\mu}(x)\,dx \cdot G^{-1}(\tau)\oi\right) + o(1),
$$
when $t \to -\infty$.

\smallskip

Having established this asymptotics, we want to relate the integral in the right-hand side of \eqref{eq27} to the free evolution one finds in \eqref{eq0}. To this end, we first define $Y_{X,+} \in L^2(\Hh_0)$ by the relations $\W_{\Hh_0} Y_{X,+} = f_{\bf ac}\ov{D_{\mu}}$.  By the spectral theorem, we have $(x+i)f \in L^2(\mu)$. So, the inclusion $f_{\bf ac} \ov{D_{\mu}} \in L^1(\R) \cap L^2(\R)$
and property $(D)$ shows that   $\chi_{\dlt}\widetilde U^0_{\gamma, t} Y_{X,+}$ can be considered as an element of $H$ for every choice of the phase function $\gamma$. 
To understand this function better, we notice that \eqref{eq25} implies
\begin{equation}\label{eq48}
(U_{t}^{0}Y_{X,+})(\Tc(\tau)) 
=
\frac{1}{2\sqrt{\pi}}\int_{\R}e^{i(t - \Tc(\tau))x}f_{\bf ac}\ov{D_{\mu}}\,dx \cdot \omi
+
\frac{1}{2\sqrt{\pi}}\int_{\R}e^{i(t + \Tc(\tau))x}f_{\bf ac}\ov{D_{\mu}}\,dx \cdot \oi,
\end{equation}
where the integrals converge absolutely. Riemann-Lebesgue lemma gives
\begin{equation}\label{sa5}
\lim_{t \to +\infty}\int_{\R}e^{i(t + \Tc(\tau))x}f_{\bf ac}(x)\ov{D_{\mu}(x)}\,dx = 0
\end{equation}
uniformly with respect to $\tau \in \R_+$.  Next, we indicate how the phase function $\gamma$ is chosen in \eqref{eq0}.
For a.e. $\tau \ge 0$,  we have $\Hh(\tau)\ge 0$ and $\det G(\tau)=1$ so Lemma \ref{l3} allows us to choose  $\phi(\tau) \in [0, 2\pi)$ such that 
\begin{equation}\label{choice}
\sqrt{\Hh(\tau)} G^{-1}(\tau)\Sigma_{\phi(\tau)}^{-1} \ge 0.
\end{equation}
That $\phi$ is Lebesgue-measurable.
Observe that for $u \in [0, 2\pi)$ we have
$$
e^{iu}\omi = (\cos u + i \sin u) \omi = 
\Sigma_{-u}\omi = \Sigma_{u}^{-1}\omi.
$$
We can then choose measurable $\gamma$ so that $|\gamma(\tau)| = 1$ on $\R_+$ and
\begin{equation}\label{sa8}
\gamma(\tau)\Sigma_{\phi(\tau)}^{-1}\omi =  \ov\alpha_\tau \omi,
\end{equation}
i.e., 
\begin{equation}\label{sa7}
\gamma(\tau) =\ov\alpha_\tau e^{-i\phi(\tau)}.
\end{equation}
Note also that for the rotation matrix  $\Sigma_{u}$ from Lemma \ref{l3} we can use the definition of $\Psi_1$ and $\Psi_2$ in \eqref{sl2} to get
\begin{equation}\label{eq38}
\sup_{u \in [0,2\pi]}\|\chi_{\dlt} G^{-1}\Sigma_{u}^{-1}\oi\|_{L^2(\Hh)} \le 4\|\Psi_1\|_{L^2(\Hh,\dlt)} + 4\|\Psi_2\|_{L^2(\Hh,\dlt)},
\end{equation}
where the right-hand side is uniformly bounded by \eqref{eq26}. Then,  \eqref{eq48}, \eqref{sa5}, and \eqref{eq38} imply  
\begin{align}
\gamma(\tau)\chi_{\dlt}G^{-1}\Sigma_{\phi(\tau)}^{-1}(U_{t}^{0}Y_{X,+})(\Tc(\tau)) 
&=
\frac{\gamma(\tau)\chi_{\dlt}}{2\sqrt{\pi}}\int_{\R}e^{i(t - \Tc(\tau))x}f_{\bf ac}\ov{D_{\mu}}\,dx \cdot G^{-1}(\tau)\Sigma_{\phi(\tau)}^{-1}\omi
+
o(1),\notag \\
&\!\!\!\stackrel{\eqref{sa8}}{=}
\frac{\ov\alpha_\tau\chi_{\dlt}}{2\sqrt{\pi}}\int_{\R}e^{i(t - \Tc(\tau))x}f_{\bf ac}\ov{D_{\mu}}\,dx \cdot G^{-1}(\tau)\omi
+
o(1), \label{eq40}
\end{align}
with $o(1)$ in $L^2(\Hh)$ as $t \to +\infty$. Combining this with \eqref{eq27}, we obtain
\begin{equation}\label{eq28}
\chi_{\dlt}(\tau)U_tX(\tau) - P_H \left(\gamma(\tau)\chi_{\dlt}(\tau)G^{-1}\Sigma_{\phi(\tau)}^{-1}(U_{t}^{0}Y_{X,+})(\Tc(\tau))\right)  \to 0, \qquad t \to +\infty, 
\end{equation}
in $H$. Similarly, for the same choice of $\gamma$ and $Y_{X,-}$ defined by $\W_{\Hh_0} Y_{X,-} = f_{\bf ac}D_{\mu}$ we have 
$$
\ov{\gamma(\tau)}\Sigma_{\phi(\tau)}^{-1}\oi =  \alpha_\tau \oi,
$$
by taking conjugation of \eqref{sa8} and 
\begin{equation*}
\chi_{\dlt}(\tau)U_tX(\tau) - P_H \left( \overline{\gamma(\tau)}\chi_{\dlt}(\tau)G^{-1}\Sigma_{\phi(\tau)}^{-1}(U_{t}^{0}Y_{X,-})(\Tc(\tau))\right)  \to 0, \qquad t \to -\infty, 
\end{equation*}
in $H$. Consider the set $\dlt^{0}=\dlt \cap \{\tau\colon \det\Hh(\tau) = 0\}$ and denote $\dlt'=\dlt \setminus \dlt^0$. Recall the formula \eqref{eq40} and note that
$$
\|\chi_{\dlt^0} G^{-1}\omi\|_{L^2(\Hh)}^{2} = 4\|\Psi_1\|_{L^2(\Hh,\dlt^0 )}^{2}.
$$
Similarly to \eqref{eq26}, we have
$$
2\|\Psi_1\|_{L^2(\Hh,\dlt^0 )}^{2} \le \|\K'_{\Hh}\|_{L^1(\dlt^0)} + 2\int_{\dlt^0}\sqrt{\det \Hh}\,d\tau = \|\K'_{\Hh}\|_{L^1(\dlt^0)},
$$
which tends to zero as $t \to +\infty$ thanks to \eqref{sa4}. Together with \eqref{eq27} this yields 
\begin{equation}\label{eq29}
\lim_{t \to +\infty}\|\chi_{\dlt^0}\cdot U_tX\|_{H} = 0.
\end{equation}
We also have $\|\chi_{\dlt^0}\cdot\widetilde U^0_{t, \gamma} Y_{X,+}\|_{H} = 0$ 
by the definition of $\widetilde U^0_{t, \gamma}$. From \eqref{eq28}, it is now clear  that the relation 
\begin{equation}\label{sa60}
\|U_t X - \widetilde U^0_{t, \gamma} Y_{X,+}\|_{L^2(\Hh, \dlt)} \to 0, \qquad t \to +\infty,
\end{equation}
is equivalent to the relation 
$$
\|P_H\left(\gamma \chi_{\Delta_{b,t}}G^{-1}\Sigma_{\phi}^{-1}(U_{t}^{0}Y_{X,+})(\Tc(\cdot))\right) - \widetilde U^0_{t, \gamma} Y_{X,+}\|_{L^2(\Hh, \dlt')} \to 0, \qquad t \to +\infty.
$$
Since $\widetilde U^0_{t, \gamma} Y_{X,+}$ belongs to $H$ by the definition of $\widetilde U^0_{t, \gamma}$, we have $P_H \widetilde U^0_{t, \gamma} Y_{X,+} = \widetilde U^0_{t, \gamma} Y_{X,+}$. Moreover, one gets $\|P_{H} Y\|_{L^2(\Hh, \dlt')} = \|Y\|_{L^2(\Hh, \dlt')}$ for every $Y \in L^2(\Hh)$, because the operator $Y \mapsto \chi_{\dlt'} Y$ is the orthogonal projector in $L^2(\Hh)$ onto a subspace in $H$. Therefore, \eqref{sa60} will follow if we prove
$$
\|\gamma G^{-1}\Sigma_{\phi}^{-1}(U_{t}^{0}Y_{X,+})(\Tc(\cdot)) - \gamma\sqrt[4]{\det\Hh}\Hh^{-\frac 12}(U_{t}^{0}Y_{X,+})(\Tc(\cdot))\|_{L^2(\Hh, \dlt')} \to 0,
$$
as $t \to +\infty$. Similarly, 
$$
\|U_t X - \widetilde U^0_{t, \gamma} Y_{X,-}\|_{L^2(\Hh, \dlt)} \to 0, \qquad t \to -\infty,
$$
follows from 
$$
\|\bar\gamma G^{-1}\Sigma_{\phi}^{-1}(U_{t}^{0}Y_{X,-})(\Tc(\cdot)) - \bar\gamma\sqrt[4]{\det\Hh}\Hh^{-\frac 12}(U_{t}^{0}Y_{X,-})(\Tc(\cdot))\|_{L^2(\Hh, \dlt')} \to 0,
$$
as $t \to -\infty$. Since $|\gamma|=1$ on $\R_+$, we only need to prove
\begin{equation}\label{sa6}
\|G^{-1}\Sigma_{\phi}^{-1}(U_{t}^{0}Y_{X,\pm})(\Tc(\cdot)) - \sqrt[4]{\det\Hh}\Hh^{-\frac 12}(U_{t}^{0}Y_{X,\pm})(\Tc(\cdot))\|_{L^2(\Hh, \dlt')} \to 0,
\end{equation}
where $t \to \pm\infty$. Noting that $\|X\|_{L^2(\Hh, \dlt')} = \|\sqrt{\Hh} X\|_{L^2(\Hh_0, \dlt')}$, we see  that the norm in \eqref{sa6} is equal to 
\begin{align*}
\ldots &=\|(G^{-1}\Sigma_{\phi}^{-1}- \sqrt[4]{\det\Hh}\Hh^{-\frac 12})(U_{t}^{0}Y_{X,\pm})(\Tc(\cdot))\|_{L^2(\Hh, \dlt')}\\ 
&= \bigl\|\bigl[\sqrt{\Hh}G^{-1}\Sigma_{\phi}^{-1} - \sqrt[4]{\det\Hh}\bigr](U_{t}^{0}Y_{X,\pm})(\Tc(\cdot))\bigr\|_{L^2(\Hh_0, \dlt')}.
\end{align*}
Thus, \eqref{sa6} can be rewritten further in the form
\begin{equation}\label{eq30}
\bigl\|\bigl[V - \sqrt{\det V}\bigr](U^0_t Y_{X,\pm})(\Tc(\cdot))\bigr\|_{L^2(\Hh_0, \dlt')} \to 0,
\end{equation}
where the matrix-function $V$ is defined by $V = \sqrt{\Hh}G^{-1}\Sigma_\phi^{-1}$. Recall that $V \ge 0$ by the choice of $\phi$ we made in \eqref{choice}. 
For each $\tau \in \dlt'$, let $e_1(\tau)$ and $e_2(\tau)$ denote the orthonormal eigenvectors of $V(\tau)$ corresponding to the eigenvalues $\lambda_1(\tau)$ and $\lambda_2(\tau)$. Then, for every vector $e = c_1e_1(\tau) + c_2e_2(\tau)$ in $\C^2$, we have
\begin{align*}
\|(V - \sqrt{\det V})e\|_{\C^2}^{2} 
&= (\lambda_1(\tau)-\sqrt{\lambda_1(\tau)\lambda_2(\tau)})^2|c_1|^2 + (\lambda_2(\tau)-\sqrt{\lambda_1(\tau)\lambda_2(\tau)})^2|c_2|^2\\
&\le (\lambda_1(\tau) - \lambda_2(\tau))^2\|e\|_{\C^2}^{2},
\end{align*}
due to the fact that
$$
(a-\sqrt{ab})^2 + (b-\sqrt{ab})^2 = (a+b)(\sqrt{a}-\sqrt{b})^2 \le (\sqrt{a}+\sqrt{b})^2(\sqrt{a}-\sqrt{b})^2 = (a-b)^2.
$$
On the other hand, 
$$
(\lambda_1(\tau) - \lambda_2(\tau))^2 = \trace (V^2)  - 2\det V.
$$
Since $V^2=V^*V$, we can write
$
\trace V^2=\trace (  \Sigma_\phi  (G^{-1})^*{\Hh}G^{-1}\Sigma_\phi^{-1})=\trace ((G^{-1})^*\Hh G^{-1})$. 
So, 
$$
(\lambda_1(\tau) - \lambda_2(\tau))^2 =  \trace ((G^{-1})^*\Hh G^{-1}) - 2\sqrt{\det \Hh} = -\K_{\Hh}',
$$
as follows from \eqref{eq22}. Since $U^0_t Y_{X,\pm} \in L^{\infty}(\Hh_0)$, $\dist(0, \dlt') \to \infty$, and $\K'_{\Hh} \in L^1(\R_+)$, we see that \eqref{eq30} holds. Hence,  $Y_{X, \pm}$ satisfy \eqref{eq0}.

\medskip

Now, consider the case where $X \in H$ is an arbitrary element (that is, we do not assume now that $X \in \dom \Di_{\Hh} \cap H_c$). Lemma \ref{l4} allows us to find $X_n\in \dom \Di_{\Hh} \cap H_c$ such that
$$
X = \lim_{n\to \infty}X_n
$$ 
and this limit is in $L^2(\Hh)$-norm.
Let $Y_{X_n, \pm}$ be the corresponding elements of $L^2(\Hh_0)$:  if $f_n = \W_{\Hh}(X_n)$, then $Y_{X_n, +} := \W_{\Hh_{0}}^{-1}(f_n \ov{D_{\mu}}\cdot \chi_{\Omega_{\bf ac}(\mu)})$,   $Y_{X_n, -} := \W_{\Hh_{0}}^{-1}(f_n {D_{\mu}}\cdot \chi_{\Omega_{\bf ac}(\mu)})$, and
$$
\|Y_{X_n, +}\|_{L^2(\Hh_0)} = \|\W_{\Hh_{0}}(Y_{X_n, +})\|_{L^2(\R)} = \|f_n\ov{D_{\mu}}\cdot \chi_{\Omega_{\bf ac}(\mu)}\|_{L^2(\R)} \le \|f_n\|_{L^2(\mu)} = \|X_n\|_{L^2(\Hh)}.
$$
A similar relation holds for $Y_{X_n, -}$. Since $\{X_n\}$ converges to $X$, the sequence $f_n$ converges to $f=\W_{\Hh}X$ in $L^2(\mu)$. 
In particular, $(f_n-f)\ov D_\mu \cdot \chi_{\Omega_{\bf ac}(\mu)}\to 0$ and $(f_n-f) D_\mu \cdot \chi_{\Omega_{\bf ac}(\mu)}\to 0$  in $L^2(\R)$. The sequences $\{Y_{X_n,\pm}\}$ converge and we denote $Y_{X, \pm} = \lim_{n\to\infty}Y_{X_n, \pm}$. In fact,   $Y_{X, +}=\W_{\Hh_{0}}^{-1}(f \ov{D_{\mu}}\cdot \chi_{\Omega_{\bf ac}(\mu)})$ and $Y_{X, -}=\W_{\Hh_{0}}^{-1}(f {D_{\mu}}\cdot \chi_{\Omega_{\bf ac}(\mu)})$, which proves \eqref{sa9}.
Moreover, for each $t$ we have 
$$
\|\widetilde U^{0}_{\gamma, t} (Y_{X_n, \pm}) - \widetilde U^{0}_{\gamma, t} (Y_{X,\pm})\|_{L^2(\Hh)}^{2} =
\|Y_{X_n, \pm} -  Y_{X,\pm}\|_{L^2(\Hh_0)}^{2},
$$
by $(C)$. Now, given that  \eqref{eq0}  holds for every $X_n$, we can extend \eqref{eq0} to all $X$ by the standard approximation argument.\smallskip

{\it Uniqueness of $Y_{X,\pm}$.}\smallskip

 We will prove uniqueness of $Y_{X,+}$, the argument for $Y_{X,-}$ is similar. Suppose that $Y_{X,+}$ and $\widetilde Y_{X,+}$ both satisfy \eqref{eq0} for some $X \in H$. Denote $Y_0 = Y_{X,+} - \widetilde Y_{X,+}$ and let $h$ be such that $\W_{\Hh_0}Y_0 = h$. We have 
$$
\lim_{t\to+\infty}\|\widetilde U_{\gamma, t}^{0}Y_0\|_{L^2(\Hh, \dlt)} = 0, \qquad \dlt = \bigl[\Lc(|t| -b), \Lc(|t| + b)\bigr],
$$
for every $b>  0$. By Lemma \ref{l1}, we have
\begin{equation}\label{eq32}
U_{t}^{0}Y_0
=\frac{1}{\sqrt 2} \widehat h(\tau-t)\cdot \omi+\frac{1}{\sqrt 2} \widehat h(-(\tau+t))\cdot \oi.
\end{equation}
Then,
$$
\lim_{t\to+\infty}\int_{\Lc(t-b)}^{\Lc(t+b)}
\sqrt{\det\Hh(\tau)}    ( |\widehat h(\Tc(\tau)-t)+\widehat h(-(\Tc(\tau)+t))|^2+|\widehat h(\Tc(\tau)-t)-\widehat h(-(\Tc(\tau)+t))|^2) \,d\tau = 0,
$$
and, after changing variables,
$$
0=\lim_{t\to+\infty}\int_{t-b}^{t+b}
|\widehat h(s-t)|^2\,ds = \int_{-b}^{b}
|\widehat h(\alpha)|^2\,d\alpha.
$$
Since $b$ is arbitrary, we get $\widehat h= 0$ a.e. on $\R$. Hence, $h=0$ a.e. which gives $Y_0 = 0$ and so $Y_{X,+}$ is defined uniquely by $X$. 

\medskip

To complete the proof, it remains to check that $\alpha_\tau  = 1$ and $G(\tau) > 0$ for every $\tau \in \R_+$ in the case when $\Hh$ is diagonal. Then,  $\gamma(\tau) = 1$ as well by \eqref{sa7}. 
We have $\Rr_{\Hh}(\tau) = 0$, $\tau \in \R_+$, for any diagonal Hamiltonian $\Hh$, see Lemma 2.2 in \cite{BD2017}. Then, $G(\tau)$ is a diagonal matrix with positive entries for every $\tau \in \R_+$, in particular, $G(\tau) > 0$. Suppose for a moment that $\det\Hh = 1$ almost everywhere on $\R_+$. 
Then, formula $(42)$ in \cite{B2018} for $z = i$ together with the relation $\Rr_{\Hh} = 0$ says that $(e^{-\tau}\Et_\tau(i))' = -\K'_{\Hh}(\tau)e^{-\tau}\Et_\tau(i)$. Since $\Et_0(i) > 0$, this shows that for such Hamiltonians $\Hh$ we have $\Et_\tau(i) > 0$, $\tau \in \R_+$. 
That implies $\alpha_\tau = 1$ for all $\tau \in \R_+$. Now let $\Hh \in \szcs$ be a diagonal Hamiltonian such that $\det\Hh > 0$ almost everywhere on~$\R_+$. Define the new Hamiltonian $\widetilde{\Hh}\colon \tau \mapsto (\det\Hh(\Lc(\tau)))^{-\frac 12}\Hh(\Lc(\tau))$, $\det\widetilde{\Hh} = 1$ on $\R_+$. 
The function $\tau \mapsto \Theta(\Lc(\tau),z)$ then solves Cauchy problem \eqref{cs} for $\widetilde{\Hh}$. The previous reasoning shows  that for the corresponding coefficient $\widetilde\alpha_\tau$ we have $\widetilde{\alpha}_\tau = 1$, $\tau \in \R_+$. But $\widetilde{\alpha}_\tau = \alpha_{\Lc(\tau)}$, and we see that $\alpha_\tau = 1$ for all diagonal Hamiltonians $\Hh \in \szcs$ such that $\det\Hh > 0$ almost everywhere on $\R_+$. 
Then, the general case follows via an approximation argument by considering Hamiltonians of the form $\Hh_\eps = \eps\idm + \Hh$, and letting $\eps \to 0$. \qed  

\medskip

\subsection{Scattering and wave operators} The  following theorem answers the question: does the asymptotics $Y_{X,\pm}$ of a state $X$ under the evolution $U_t$ determine the state $X$ itself?  It also strengthens Theorem~\ref{t1}.
\begin{Thm}\label{t3}
Let $\Hh \in \szcs$, let $\Di_{\Hh}$ be the corresponding self-adjoint operator \eqref{eq42} on $H$, and let $\mu = w\,dx + \mus$ be the main spectral measure for $\Hh$. Then, the strong wave operators
\begin{equation}\label{eq35}
W_{\pm} = \lim_{t\to\pm\infty}U_{t}^{-1}\widetilde U^{0}_{\gamma, t}
\end{equation}
exist and are complete, i.e., they are correctly defined (the limits are understood in the strong operator topology) and unitary as operators from $L^2(\Hh_0)$ onto the absolutely continuous subspace $H_{\bf ac}$ of $\Di_{\Hh}$. Moreover, if $X\in H$, $Y_{X,\pm}$ are defined by \eqref{sa9} and $P_{\bf ac}$ denotes the orthogonal projection in $H$ onto the absolutely continuous subspace of $\Di_{\Hh}$, then we have $W_{\pm}Y_{X,\pm} = P_{\bf ac}X$. Hence, $Y_{X,\pm}$ determine $P_{\bf ac}X$ uniquely and we have $Y_{X_1,\pm} = Y_{X_2,\pm}$ in Theorem \ref{t1} if and only if $P_{\bf ac}X_1 = P_{\bf ac}X_2$. 
The scattering operator 
$$
S = W_+^{-1}W_{-}, \qquad S\colon Y_{X,-} \mapsto Y_{X,+},
$$ 
is a unitary operator on $L^2(\Hh_0)$, and its spectral representation takes the form
\begin{equation}\label{eq5}
\W_{\Hh_0}S\W_{\Hh_0}^{-1} f_0 = \frac{\ov D_{\mu}}{D_{\mu}}f_0, \qquad f_0 \in L^2(\R),
\end{equation}
where $D_{\mu}$ is the Szeg\H{o} function of $\mu$. {In particular, the operator $S$ does not depend on the choice of the phase function $\gamma$ in Theorem \ref{t1}.}
\end{Thm}

\medskip

\beginpf Let us first prove that the limit  in \eqref{eq35} exists as $t \to+\infty$.  The argument for $t\to -\infty$ is similar.  In fact, we claim that for an arbitrary $X$, we have
\begin{equation}\label{sa13}
\lim_{t\to+\infty}U_{t}^{-1}\widetilde U^{0}_{\gamma, t} Y_{X, +}=P_{\bf ac} X.
\end{equation}
Indeed, denoting $\Delta_{b,t} = [\Lc(|t|-b), \Lc(|t|+b)]$ for some positive numbers $b$, we have
\begin{align*}
U_{t}^{-1}\widetilde U^{0}_{\gamma,t}Y_{X, +} 
&= 
U_{t}^{-1}\chi_{\Delta_{b,t}}\widetilde U^{0}_{\gamma,t} Y_{X, +} + U_{t}^{-1}\chi_{\R_+\setminus\Delta_{b,t}}\widetilde U^{0}_{\gamma,t} Y_{X, +}\\
&= 
U_{t}^{-1} \chi_{\Delta_{b,t}} U_t X + U_{t}^{-1}(\chi_{\Delta_{b,t}}\widetilde U^{0}_{\gamma,t}Y_{X, +} - \chi_{\Delta_{b,t}}U_t X) + U_{t}^{-1}\chi_{\R_+\setminus\Delta_{b,t}}\widetilde U^{0}_{\gamma,t} Y_{X, +}.
\end{align*}
 To get our claim, it is enough to prove that 
\begin{itemize}
\item[$(a)$] $
\lim_{b\to+\infty}\limsup_{t \to +\infty}\|\chi_{\R_+\setminus\Delta_{b,t}}\widetilde U^{0}_{\gamma,t} Y_{X, +}\|_{L^2(\Hh)}=0,$
\item[$(b)$]$\lim_{t \to +\infty}\|\chi_{\Delta_{b,t}}\widetilde U^{0}_{\gamma,t} Y_{X, +} - \chi_{\Delta_{b,t}}U_t X\|_{L^2(\Hh)} = 0$ for every $b \ge 0$,
\item[$(c)$] $\lim_{b\to+\infty}\limsup_{t \to +\infty}\|U_{t}^{-1} \chi_{\Delta_{b,t}} U_t X - P_{\bf ac}X\|_{L^2(\Hh)} = 0$.
\end{itemize}
Clearly, $(b)$ is just a restatement of Theorem \ref{t1}. To check $(a)$, observe that
\begin{align*}
\|\chi_{\R\setminus\Delta_{b,t}}\widetilde U^{0}_{\gamma,t}Y_{X, +}\|_{L^2(\Hh)}
&=\int_{\R_+\setminus\Delta_{b,t}}\sqrt{\det\Hh(\tau)}\|U^{0}_{t}Y_{X, +}(\Tc(\tau))\|_{\C^2}^{2}\,d\tau\\
&=\int_{\R_+\setminus [|t|-b,|t|+b]}\|U^{0}_{t}Y_{X, +}(\tau)\|_{\C^2}^{2}\,d\tau. 
\end{align*}
Now, \eqref{sa11} yields $(a)$.  It remains to prove $(c)$. First, notice that
\[
\|U_{t}^{-1} \chi_{\Delta_{b,t}} U_t X - P_{\bf ac}X\|^2_{L^2(\Hh)}=\|\chi_{\Delta_{b,t}} (U_t X - U_tP_{\bf ac}X)\|_{L^2(\Hh)}^2+\| \chi_{\R_+\backslash \Delta_{b,t}}  U_tP_{\bf ac}X\|_{L^2(\Hh)}^2.
\]
By Theorem \ref{t1}, we have $Y_{X,+}=Y_{P_{\bf ac}X,+}$ and, therefore,
$
\lim_{t\to +\infty}\|\chi_{\Delta_{b,t}} (U_t X - U_tP_{\bf ac}X)\|_{L^2(\Hh)}=0.
$
Finally, 
\begin{align*}
\| \chi_{\R_+\backslash \Delta_{b,t}}  U_tP_{\bf ac}X\|^2_{L^2(\Hh)}
&=\| U_tP_{\bf ac}X\|^2_{L^2(\Hh)}-\| \chi_{ \Delta_{b,t}}  U_tP_{\bf ac}X\|^2_{L^2(\Hh)}\\
&=\|P_{\bf ac}X\|^2_{L^2(\Hh)}-\| \chi_{ \Delta_{b,t}}  U_tP_{\bf ac}X\|^2_{L^2(\Hh)},
\end{align*}
and, by Theorem \ref{t1},
\begin{eqnarray*}
\limsup_{t\to +\infty}(\|P_{\bf ac}X\|^2_{L^2(\Hh)}-\| \chi_{ \Delta_{b,t}}  U_tP_{\bf ac}X\|^2_{L^2(\Hh)})=\|P_{\bf ac}X\|^2_{L^2(\Hh)}-\liminf_{t\to +\infty}\| \chi_{ \Delta_{b,t}}  \widetilde U^{0}_{\gamma,t}Y_{P_{\bf ac}X, +}   \|^2_{L^2(\Hh)}\\
=\|P_{\bf ac}X\|^2_{L^2(\Hh)}-\liminf_{t\to +\infty}(\|   \widetilde U^{0}_{\gamma,t}Y_{P_{\bf ac}X, +}   \|^2_{L^2(\Hh)}-\| \chi_{ \R_+\backslash\Delta_{b,t}}  \widetilde U^{0}_{\gamma,t}Y_{P_{\bf ac}X, +}   \|^2_{L^2(\Hh)}).
\end{eqnarray*}
By \eqref{sa9}, we get
\begin{equation}\label{rab1}
\|  \widetilde U^{0}_{\gamma,t} Y_{P_{\bf ac}X, +}   \|^2_{L^2(\Hh)}=\|P_{\bf ac}X\|^2_{L^2(\Hh)}
\end{equation}
 and  \[
\lim_{b\to+\infty}\limsup_{t \to +\infty}\|\chi_{\R_+\setminus\Delta_{b,t}}\widetilde U^{0}_{\gamma,t}Y_{P_{\bf ac}X, +}\|_{L^2(\Hh)}=0\] follows from  $(a)$. Hence, we get $(c)$. 

Then, \eqref{sa9} implies, in particular, that the map $\mathcal{Y}_+\colon X\mapsto Y_{X,+}$ is the unitary map from $H_{\bf ac}$ onto $L^2(\Hh_0)$. So, $\lim_{t\to+\infty}U_{t}^{-1}\widetilde U^{0}_{\gamma,t} Y$ exists for every $Y\in L^2(\Hh_0)$ and $\lim_{t\to+\infty}U_{t}^{-1}\widetilde U^{0}_{\gamma,t} Y=\mathcal{Y}_+^{-1}Y$.

Summarizing, we have proved that the strong wave operator $W_+$ in \eqref{eq35} exists and $W_+ Y_{X, +} = P_{\bf ac} X$ for every $X$,  where $Y_{X, +}$ is defined as in Theorem \ref{t1}. Analogously, one can check the existence of the wave operator $W_{-}$ and prove the formula $W_- Y_{X, -} = P_{\bf ac} X$. All other assertions of the theorem are simple consequences of these two facts. \qed

\medskip
The following corollary implies, in particular, that if Szeg\H{o} measure $\mu$ is purely a.c., then every $X\in H$ propagates and the global $L^2(\Hh)$ asymptotics holds for $U_tX$. The reader can compare it to Theorem \ref{t1} which establishes the asymptotics over the finite interval.

\begin{Cor}\label{c2}
Let $\Hh$ be a Hamiltonian of class $\szcs$. Then, there exists a phase function $\gamma$ such that the following assertion holds. For every $X \in H_{\bf ac}$, there are uniquely defined  $Y_{X,\pm} \in L^2(\Hh_0)$ such that
\begin{equation}\label{eq000}
\lim_{t \to \pm\infty}\|U_t X - \widetilde U^{0}_{\gamma,t} Y_{X,\pm}\|_{L^2(\Hh)} =0, 
\end{equation}
Moreover, if $\Hh$ is diagonal, then one can take $\gamma = 1$ on $\R_+$.

\end{Cor} 
\beginpf Indeed, if $X \in H_{\bf ac}$, then 
$X = P_{\bf ac}X$, and \eqref{eq000} is equivalent to 
$$
\lim_{t \to \pm\infty}\|P_{\bf ac}X - U_t^{-1} \widetilde U^{0}_{\gamma,t} Y_{X,\pm}\|_{L^2(\Hh)} =0,
$$
that holds by Theorem \ref{t3}. \qed\medskip

\subsection{Dynamical classification of spectral types} Our analysis allows to detect the spectral types of $\Di_\Hh$ by observing the long-time dynamics of $U_t$.

\medskip
Suppose $X\in H$ is given. Denote the orthogonal projections to absolutely continuous, singular continuous, and pure point subspaces of $\Di_{\Hh}$ by $P_{\bf ac}, P_{\bf sc},$ and $P_{\bf pp}$, respectively. 
Our next result gives the dynamical characterization of whether $X$ has nontrivial projections to any of these subspaces.

\begin{Thm}\label{t4}
Let $\Hh \in \szcs$. Then, for every $X \in H$ we have
\begin{align}
\lim_{b \to +\infty}\lim_{\tb \to +\infty}\frac{1}{\tb}\int_{0}^{\tb}\|U_{t}X\|_{L^2(\Hh, [0, \Lc(b)])}^{2}\,dt 
&= \|P_{\bf pp}X\|_{L^2(\Hh)}^{2}, \label{eq83bis}\\
\lim_{b \to +\infty}\lim_{\tb \to +\infty}\frac{1}{\tb}\int_{0}^{\tb}\|U_t X\|_{L^2(\Hh, [\Lc(b), \Lc(t-b)])}^{2}\,dt &= \|P_{\bf sc} X\|_{L^2(\Hh)}^{2},\label{eq89bis}\\
\lim_{b \to +\infty}\lim_{t \to +\infty}\|U_t X\|_{L^2(\Hh, [\Lc(t-b), \Lc(t+b)])}^{2} &= \|P_{\bf ac} X\|_{L^2(\Hh)}^{2}, \label{eq82bis}\\
\lim_{b \to +\infty}\lim_{t \to +\infty}\|U_t X\|_{L^2(\Hh, [\Lc(t+b), +\infty))} & = 0. \label{eq91bis}
\end{align} 
The analogous statements hold when $t\to +\infty$ is replaced by $t\to -\infty$.
\end{Thm}
\beginpf
We start with proving \eqref{eq91bis}. Given $X$ and $\eps>0$, we can find $X_\eps\in H_c$ such that $\|X-X_\eps\|_{L^2(\Hh)}\le \eps$. From Theorem \ref{p0bb}, we get
\[
\lim_{b \to +\infty}\limsup_{t \to +\infty}\|U_t X_\eps\|_{L^2(\Hh, [\Lc(b+t), +\infty))}= 0.
\]
Therefore, 
\begin{align*}
\limsup_{b \to +\infty}\limsup_{t \to +\infty}\|U_t X\|_{L^2(\Hh, [\Lc(b+t), +\infty))}
\le &\limsup_{b \to +\infty}\limsup_{t \to +\infty}\|U_t (X-X_\eps)\|_{L^2(\Hh, [\Lc(b+t), +\infty))}\\
&+\limsup_{b \to +\infty}\limsup_{t \to +\infty}\|U_t X_\eps\|_{L^2(\Hh, [\Lc(b+t), +\infty))}\le \eps,
\end{align*}
because $\|U_t(X-X_\eps)\|_{L^2(\Hh)}=\|X-X_\eps\|_{L^2(\Hh)}\le \eps$. Since $\eps$ is arbitrary, we get \eqref{eq91bis}.

\smallskip

Formula \eqref{eq82bis} is immediate from Theorem \ref{t1}.

\smallskip

To prove \eqref{eq83bis}, we apply Lemma \ref{lemo1} from Section \ref{app5}. Take $\Lambda>0$ and $b\in \R_+ \setminus \bigcup_{I \in \mathfrak{I}(\Hh)} I.$ Let $P_{[-\Lambda,\Lambda]}$ be the orthogonal projection associated with the spectral decomposition of~$\Di_{\Hh}$. We claim that the operator
$
\chi_{[0,b]}P_{[-\Lambda,\Lambda]}
$ is compact in $H$. Indeed, this follows from the formula
\[
\chi_{[0,b]}(\tau)(P_{[-\Lambda,\Lambda]}X)(\tau)=\frac{\chi_{[0,b]}(\tau)}{\sqrt\pi}\int_{[-\Lambda,\Lambda]}\Theta(\tau,x)(\W_{\Hh}X)(x)\,d\mu(x)
\]
and the fact that the set $\{\int_{[-\Lambda,\Lambda]}\Theta(\tau,x)(\W_{\Hh}X)\,d\mu\colon \|X\|_{H}\le 1\}$ is precompact in $C[0,b]$ by Arzela-Ascoli theorem. Hence, by Lemma \ref{lemo1}  applied to Hilbert space $H$, operator $\Di_{\Hh}$, and  $A=\chi_{[0,b]}$, one has
\[
\lim_{\tb\to\infty}\frac 1\tb\int_0^\tb \|\chi_{[0,b]}U_tP_{[-\Lambda,\Lambda]}X\|_{L^2(\Hh)}^2dt=\sum_{j}\|\chi_{[0,b]}P_{\{E_j\}}P_{[-\Lambda,\Lambda]}X\|^2_{L^2(\Hh)},
\]
where $P_{\{E_j\}}$ is orthogonal projection corresponding to eigenvalue $E_j$ of $\Di_{\Hh}$ and the sum is done over all eigenvalues.
Taking $b$ to infinity (see Corollary 2 in \cite{Robinson}), we have
\begin{equation}\label{des6}
\lim_{b\to \infty}\lim_{\tb\to\infty}\frac 1\tb \int_0^\tb \|\chi_{[0,b]}U_t P_{[-\Lambda,\Lambda]}X\|_{L^2(\Hh)}^2dt=\|P_{\bf pp}P_{[-\Lambda,\Lambda]}X\|^2_{L^2(\Hh)},
\end{equation}
for every $X$. Now, taking $\Lambda\to \infty$, we get \eqref{eq83bis}.\smallskip

We are left with showing \eqref{eq89bis}. Fix any $X$ and $b>0$. Then, 
\begin{align*}
\|U_tX\|^2_{L^2(\Hh)}
=&\|U_tX\|^2_{L^2(\Hh,[0,\Lc(b)])}+\|U_tX\|^2_{L^2(\Hh,[\Lc(b),\Lc(t-b)])}\\
&+\|U_tX\|^2_{L^2(\Hh,[\Lc(t-b),\Lc(t+b)])}+\|U_tX\|^2_{L^2(\Hh,[\Lc(t+b),\infty))}.
\end{align*}
We also have 
\[
\|U_tX\|^2_{L^2(\Hh)}=\|X\|^2_{L^2(\Hh)}=\|P_{\bf ac}X\|^2_{L^2(\Hh)}+\|P_{\bf sc}X\|^2_{L^2(\Hh)}+\|P_{\bf pp}X\|^2_{L^2(\Hh)}.
\]
Subtracting one identity from the other and taking the Cesaro mean, we get
\[
\lim_{b\to\infty}\lim_{\tb\to\infty}\left(\frac{1}{\tb}\int_0^\tb \|U_tX\|^2_{L^2(\Hh, [\Lc(b),\Lc(t-b)])}dt-\|P_{\bf sc}X\|^2_{L^2(\Hh)}\right)=0,
\]
as follows from already established \eqref{eq83bis}, \eqref{eq82bis}, and \eqref{eq91bis}. The arguments for $t\to -\infty$ are identical.
\qed

\medskip
We will need the following technical lemma later in the text.
\begin{Lem}\label{lod}Suppose $\Hh$ is a proper Hamiltonian,  $r\in \R_+ \setminus \bigcup_{I \in \mathfrak{I}(\Hh)} I$, $X\in H_c$ is  supported on $[0,r]$, and $(X, \left(\begin{smallmatrix} 1\\ 0\end{smallmatrix}\right))_{L^2(\Hh)}\neq 0$. If $Y$ is defined by
\begin{eqnarray}
Y(\tau)=0,\quad \tau>r,   \nonumber \\
Y(\tau)=J\int_\tau^r \Hh(s)X(s)ds, \quad \tau<r,   \label{des1}
\end{eqnarray}
then $X$ and $P_HY$ satisfy $
|(\W_{\Hh}X)(x)|^2+|(\W_{\Hh}P_{H}Y)(x)|^2>0$ for all $x\in \R$. Such $X$ and $Y$ exist.
\end{Lem}
\beginpf First, observe that 
\[
  \langle  Y(0), \left(\begin{smallmatrix} 0\\ 1\end{smallmatrix}\right)\rangle_{\C^2}=  \langle  J\int_0^r \Hh(s)X(s)ds, \left(\begin{smallmatrix} 0\\ 1\end{smallmatrix}\right)\rangle_{\C^2}=(X, \left(\begin{smallmatrix} 1\\ 0\end{smallmatrix}\right))_{L^2(\Hh)}\neq 0.
\]
Second, notice that $\Theta(\tau,0)=\left(\begin{smallmatrix} 1\\ 0\end{smallmatrix}\right)$ and so $(\W_{\Hh}X)(0)\neq 0$ given assumptions of the lemma.
Suppose $x$ is such that $x\neq 0$ and $(\W_{\Hh}X)(x)=0$. Observe that, by Lemma \ref{l5}, $\widetilde \W_{\Hh}Y=\W_{\Hh}P_HY$. Then, 
\begin{align*}
\sqrt{\pi} \cdot \widetilde\W_{\Hh}Y
&=\int_0^r \langle \Hh Y,\Theta(\tau,x)\rangle_{\C^2} d\tau=-x^{-1}\int_0^r  \langle J Y,\Theta'(\tau,x)\rangle_{\C^2} d\tau\\
&=
x^{-1}\left(  -\langle  Y(0), \left(\begin{smallmatrix} 0\\ 1\end{smallmatrix}\right)\rangle_{\C^2}+\int_0^r  \langle J Y',\Theta(\tau,x)\rangle_{\C^2} d\tau \right)\\
&=x^{-1}\left(   -\langle  Y(0), \left(\begin{smallmatrix} 0\\ 1\end{smallmatrix}\right)\rangle_{\C^2}+\int_0^r  \langle \Hh(\tau) X(\tau),\Theta(\tau,x)\rangle_{\C^2} d\tau \right)\\
&=x^{-1}\left(-\langle  Y(0), \left(\begin{smallmatrix} 0\\ 1\end{smallmatrix}\right)\rangle_{\C^2}+\sqrt{\pi} (\W_{\Hh}X)(x) \right)
\neq 0,
\end{align*}
where we used \eqref{des1} and our other assumptions. Finally, since $\Hh$ is proper, we can always find $X$ that satisfies all conditions and define $Y$ accordingly. For example, 
\[
X=\left(\begin{smallmatrix} 1\\ 0\end{smallmatrix}\right)\cdot \chi_{[0,r]}, \quad Y=\left(\begin{smallmatrix} -\int_\tau^r h(s)ds\\ \int_\tau^r h_1(s)ds\end{smallmatrix}\right)\cdot \chi_{[0,r]}
\]
is one possible choice.
\qed
\medskip

Theorem \ref{t4} gives a dynamical description of spectral types for each element $X$ but it does not tell how to detect the presence of pure point, singular continuous, and absolutely continuous spectral types for $\Di_{\Hh}$ itself. We will address it in the next theorem. Recall that $\dlt$ is defined as $\dlt = \bigl[\Lc(|t| -b), \Lc(|t| + b)\bigr]$ for $b>0$.
\begin{Thm}\label{c3}
Let $\Hh$ be a Hamiltonian of class $\szcs$. Then, the following holds true.
\begin{itemize}
\item[$(A)$] If the singular spectrum of $\Hh$ is empty, then
\begin{equation}\label{eq0001}
\lim_{b\to +\infty}\lim_{t \to +\infty}\|U_t X\|_{L^2(\Hh,\dlt)} =\|X\|_{L^2(\Hh)}
\end{equation}
for every $X\in H$.
\item[$(B)$] If there is some $X\in H_c$ for which
\begin{equation}\label{lrt1}
\lim_{b \to +\infty}\lim_{\tb \to +\infty}\frac{1}{\tb}\int_{0}^{\tb}\|U_t X\|_{L^2(\Hh, [\Lc(b), \Lc(t-b)])}^{2}\,dt=0,\end{equation}
 then $\Di_{\Hh}$ has no singular continuous spectrum.
\item[$(C)$]Let vectors $X$ and $Y$ be defined as in Lemma~\ref{lod}.  If both of the equalities
\begin{align}
\lim_{b\to +\infty}\liminf_{t \to +\infty}\|U_t X\|_{L^2(\Hh,\dlt)} &=\|X\|_{L^2(\Hh)}, \label{eq0001e} \\ 
\lim_{b\to +\infty}\liminf_{t \to +\infty}\|U_t P_{H}Y\|_{L^2(\Hh,\dlt)} &=\|Y\|_{L^2(\Hh)}, \label{eq0002e}
\end{align}
hold, then the singular spectrum of $\Di_{\Hh}$ is empty.
\item[$(D)$] Let vectors $X$ and $Y$ be defined as in Lemma~\ref{lod}. If  both of the equalities
\begin{align}
\lim_{b \to +\infty}\lim_{\tb \to +\infty}\frac{1}{\tb}\int_{0}^{\tb}\|U_{t}X\|_{L^2(\Hh, [0, \Lc(b)])}^{2}\,dt&=0, \label{des5}\\
\lim_{b \to +\infty}\lim_{\tb \to +\infty}\frac{1}{\tb}\int_{0}^{\tb}\|U_{t}Y\|_{L^2(\Hh, [0, \Lc(b)])}^{2}\,dt&=0, \label{des5bis}
\end{align}
hold, then $\Di_{\Hh}$ has no bound states. 
\end{itemize}
We get the same conclusions if  the limits $t\to +\infty$ are replaced by $t\to -\infty$.
\end{Thm}

\beginpf Suppose the singular spectrum is empty, then $X=P_{\bf ac} X$ and our claim follows from \eqref{eq91bis}.

\smallskip

Then, suppose $X \in H_c$ is such that \eqref{lrt1} holds. Recalling Theorem \ref{t1}, consider $f=\W_{\Hh} X$. 
Represent the measure $\mu=\mu_{\bf ac}+\mu_{\bf s}$ as a sum of absolutely continuous and singular components and further write $\mu_{\bf s}=\mu_{\bf sc}+\mu_{\bf pp}$ as a sum of singular continuous and pure point parts. Then, 
 \eqref{eq89bis} gives 
\[
\int_\R |f|^2\,d\mu_{\bf sc}=0.
\]
On the other hand, $f$ is an entire function that can have only countably many zeroes in $\C$. Therefore, $|f|>0$ a.e. with respect to $\mu_{\bf sc}$ and so $\mu_{\bf sc}=0$.

\smallskip

To show $(C)$, we only need to prove that \eqref{eq0001e}, \eqref{eq0002e} imply that the spectrum of $\Hh$ is purely absolutely continuous.  If $f:=\W_{\Hh}X$ and $g:=\W_{\Hh}P_H Y$, then \eqref{eq0001e} and \eqref{eq0002e} give
\[
\int_{\R} (|f|^2+|g|^2)\,d\mu_{\bf s}=0.
\]
That, however, contradicts Lemma \ref{lod} unless $\mu_{\bf s}(\R)=0$.

\smallskip

Finally, to get $(D)$, we notice that \eqref{des5}, \eqref{des5bis} and \eqref{eq83bis} give $P_{\bf pp}X=P_{\bf pp}Y=0$ which can be rewritten as
\[
\int_{\R} (|f|^2+|g|^2)\,d\mu_{\bf pp}=0,
\]
where, again, $f=\W_{\Hh}X$ and $g=\W_{\Hh}P_H Y$. Since $f$ and $g$ are entire functions that have no common zeroes by Lemma \ref{lod}, we get $\mu_{\bf pp}=0$.

\smallskip

The arguments for $t\to -\infty$ are identical. \qed

\section{Krein strings}\label{ks}
The theory of Krein strings goes back to works by M.~Krein \cite{Krein} and Feller \cite{Feller}. In this section, we recall some basic definitions and facts, explain the connection between Krein strings and diagonal canonical systems, and use it to translate some results obtained in the previous section to the new setting.

\smallskip

\subsection{Krein strings}

Let $0 < L \le \infty$ and $M$ be a non-decreasing right-continuous function on $(-\infty,L)$, satisfying $M(\xi)=0, \xi<0$. The Lebesgue-Stieltjes measure $\mf$ on $[0,L)$ is defined by $\mf[0,\xi]=M(\xi)$. We write its decomposition into the absolutely continuous and singular parts as $\mf=\mf_{\bf ac}+\mf_{\bf s}=\rho(\xi)\,d\xi+\mf_{\bf s}$. Recall that in our notation  $M(L-)=\lim_{\xi\uparrow L}M(\xi)$ and we call the $[M,L]$ pair {\bf proper} if $M$ and $L$ satisfy the following conditions
\begin{align} 
&L+M(L-)=\infty,\label{sa17}\\
&0<M(\xi)<M(L-), \quad \forall\xi\in (0,L).\label{sa15}
\end{align}
These two conditions are very natural from the point of view of spectral theory \cite{KK68}. They guarantee that the spectral measure $\sigma$ of the string operator is unique in the class of spectral measures with non-negative support. Additionally, they make sure that the map $[M,L]\mapsto \sigma$ is injective.
In this paper, we will work with proper $[M,L]$ pairs only.  
Let us consider functions $\phi$, $\psi$ defined by the integral equations  
\[	
\phi(\xi,z)=1-z\int_{[0,\xi]}(\xi-s)\phi(s,z)\,d\mf (s),
\]
\[
\psi(\xi,z)=\xi-z\int_{[0,\xi]}(\xi-s)\psi(s,z)\,d\mf (s), 
\]
where $\xi\in [0,L)$, $z \in \C$.
It is customary to extend $\phi$ and $\psi$ to $(-\infty,0)$ by $\phi(\xi,z)=1$ and $\psi(\xi,z)=\xi$ where $\xi<0$.
These functions are uniquely determined by the string $[M,L]$ and they define the  Titchmarsh-Weyl function $q$ of $[M,L]$ by
\begin{equation}\label{ssnm}
q(z)=\lim_{\xi\to L} \frac{\psi(\xi,z)}{\phi(\xi,z)}, \qquad z\in \mathbb{C}\backslash [0,\infty),
\end{equation}
see formula (2.21) in \cite{KWW} or Theorem 10.1 in \cite{KK68}. That function $q$ has the unique integral representation 
\begin{equation}\label{nnsm}
q(z)=\int_{\R_+} \frac{d\sigma(\lambda)}{\lambda-z}, \qquad z\in \mathbb{C}\backslash [0,\infty),
\end{equation}
where  $\sigma$, the {\bf main (or orthogonal) spectral measure} of the string $[M,L]$, is a nonnegative Borel measure on $\R_+$ satisfying condition
\[
\int_{\R_+}\frac{d\sigma(\lambda)}{1+\lambda}<\infty.
\]
We emphasize (see \cite{KK68}) that a proper string  is in the limit-point case if and only if 
\begin{equation}\label{secondmoment}
\int_{[0,L)}\xi^2d\mf =+\infty.
\end{equation}
However, when the integral in \eqref{secondmoment} is finite and we are in limit-circle case, the main spectral measure  with {\it non-negative} support is unique and is given by \eqref{ssnm}. Later in the text, we will focus on strings $[M,L]$ in Szeg\H{o} class. For this type of strings, the condition \eqref{secondmoment} is always satisfied.

\medskip

\noindent Similarly to \eqref{wt}, one can define the generalized Fourier transform associated with the string $[M,L]$:

\begin{equation}\label{wts}
\Us\colon v \mapsto \int_{0}^{L}v(\xi)\phi(\xi,z)\,d\mf(\xi), \qquad z \in \C,  
\end{equation}
starting with functions  $v\in L_c^2(\mf)$ that have compact support in $[0,L)$.  It is known (see Section 10 in \cite{KK68}) that $\Us$ can be extended to the unitary operator from $L^2(\mf)$ onto $L^2(\sigma)$.
The inverse map is given by (see formula (2.25) in \cite{KWW})
\[
v=\Us^{-1}(\Us v)=\int_0^\infty \phi(\xi,\lambda)(\Us v)(\lambda)\,d\sigma(\lambda), \qquad \xi \in [0, L),
\]
where the last integral can be first densely defined on $L_c^2(\sigma)$ and then extended to all of $L^2(\sigma)$. Let us define the Krein string operator $\cal{S}_M$ by
\[
\cal{S}_M=\Us^{-1}M_\lambda \Us,
\]
where $\dom \cal{S}_M:=\{v\in L^2(\mf): M_\lambda \Us v\in L^2(\sigma)$\} and we recall that $M_\lambda f$ is a function in $\lambda$ which is equal to $\lambda f(\lambda)$. Clearly, $\cal{S}_M$ is a self-adjoint operator in $L^2(\mf)$.

\bigskip

\subsection{Connection between Krein strings and canonical systems with diagonal Hamiltonians}

Suppose  $[M,L]$ is a proper string. 
Consider the increasing function  $N\colon \xi \mapsto \xi+M(\xi)$ on $[0, L)$ and let $\nf$ denote the corresponding measure, $\nf[0,\xi]=N(\xi)$ for $\xi \in [0,L)$. Condition  \eqref{sa17} is equivalent to $N(L-)=+\infty$.   Define  the function $N^{(-1)}$ as generalized inverse of $N$, see Section \ref{notation}. Using the fact that $N$ is strictly increasing, one can show that $N^{(-1)}$ is continuous on $\R_+$, and we have $N^{(-1)}(N(\xi)) = \xi$ for every $\xi \in [0, L)$. Recall that  $\rho$ is the density of the absolutely continuous part of $\mf$, so that $\mf=\rho \,d\xi+\mf_s$.  Define two functions on $\R_+$:
\begin{equation}
\htwo(\tau)= 
\begin{cases}
1, & \mbox{if } N^{(-1)}(\tau)\in \Omega_{\bf s}(\mf),\\
\frac{\rho(N^{(-1)}(\tau))}{1+\rho(N^{(-1)}(\tau))}, & \mbox{otherwise}.
\end{cases}
\label{sd1}
\end{equation}
and
\begin{equation}
\hone(\tau)= 
\begin{cases}
0, & \mbox{if } N^{(-1)}(\tau)\in \Omega_{\bf s}(\mf),\\
\frac{1}{1+\rho (N^{(-1)}(\tau))}, & \mbox{otherwise},
\end{cases}
\label{sd2}
\end{equation}
Given $\htwo,\hone$, define $\Hh_{\ast}=\diag(\htwo,\hone)$.
If $\Hh_{\ast}$ is proper, we let
\begin{equation}\label{lk3}
\xi(\tau)=\int_0^\tau \hone(s)ds, \quad L=\int_0^\infty \hone(s)ds, \quad 
M(\xi)=\int_0^{\sup\{x\colon \xi(x)=\xi\}} \htwo(s)ds,\quad \xi<L.
\end{equation}
We now collect some facts related to the well-known connection between Krein strings and diagonal canonical systems. The first of them can be found in \cite{GohKr} (see Section~8 in Chapter~6), p. 239 in \cite{DymMcKean}, or \cite{KWW}. 
\begin{Lem}\label{l20}
Formulas \eqref{sd1}, \eqref{sd2}, and \eqref{lk3} establish the bijection $[M ,L] \mapsto \diag(\htwo,\hone)$ between proper $[M, L]$ pairs  and  proper Hamiltonians $\Hh_{\ast}= \diag(\htwo,\hone)$ with unit trace. 
\end{Lem}
We want to make one comment here. The references  \cite{DymMcKean}, \cite{GohKr}, an \cite{KWW}  explain that connection for the general  strings and diagonal Hamiltonians. However, one can see that the proper strings correspond to proper Hamiltonians. Indeed, the assumption that the left end of the string is heavy is equivalent to the condition that the Hamiltonian $\Hh_{\ast} = \diag(\htwo,\hone)$ is not equal to $\excl$ on $[0,\varepsilon]$ with some  $\varepsilon>0$. Moreover,
 making the assumption that
$ L+M(L-)=+\infty$ and the right end is heavy is equivalent to saying that $\Hh_\ast$ is not equal to either $\excl$ or $\left(\begin{smallmatrix} 1&0\\0&0\end{smallmatrix}\right)$ on $(\tau_0,\infty)$ for some $\tau_0>0$. 

\medskip

For the proof of the following result, check  Theorem~4.2 in \cite{KWW}.
\begin{Lem}\label{sal16}
Let $[M,L]$ and $\Hh_{*}$ be the string and its corresponding Hamiltonian obtained via the bijection in Lemma \ref{l20}.
 Then, for the corresponding Titchmarsh-Weyl functions $q$, $m_*$, we have 
\begin{equation}\label{sdk1}
zq(z^2)=m_\ast(z), \quad z\in \mathbb{C}^+.
\end{equation} 
Consequently, the spectral measures $\sigma$, $\mu_\ast$ of $[M, L]$, $\Hh_*$ satisfy
\begin{equation}\label{sa20}
\mu_\ast([E_1,E_2])=\frac{\pi}{2} \sigma([E_1^2,E_2^2]), \quad \mu_\ast(\{0\})=\pi\sigma(\{0\})
\end{equation}
for all $0<E_1<E_2$.
\end{Lem}

Relation \eqref{sa20} shows that the operators $\Di^2_{\Hh_{\ast}}$ and $\cal{S}_{M}$ are unitarily equivalent. The unitary equivalence is given via the explicit operator $\Upsilon$ in the lemma below. 

\medskip

\begin{Lem}Let $[M,L]$ and $\Hh_{*}$ be the string and its corresponding Hamiltonian obtained via the bijection in Lemma \ref{l20}. Then, the map $\Upsilon\colon v\in L^2(\mf)\mapsto X=(v\circ N^{(-1)}, 0)^t$ is a unitary map onto the subspace $\{X=(X_1,X_2)^t\in H: X_2=0\}$ of the space $L^2(\Hh_{\ast})$. Moreover,  
\begin{equation}\label{lk6}
\Upsilon^{-1}\Di^2_{\Hh_{\ast}}\Upsilon=\cal{S}_M.
\end{equation}
\end{Lem}
Let us give a sketch of the proof of this well-known fact. The map $\Upsilon$ is correctly defined and unitary due to \eqref{sd1}, \eqref{sd2} and the change of variables in the Lebesgue-Stieltjes integral:
\[
\int_0^L |v(\xi)|^2d\mf=\int_0^\infty |v(N^{(-1)}(\tau))|^2\htwo(\tau)d\tau.
\]
To prove \eqref{lk6}, it is convenient to work on the spectral side of both $\cal{S}_M$ and $\Di_{\Hh_{\ast}}$. 
We will check that $\Di^2_{\Hh_{\ast}}\Upsilon\cal{U}_M^{-1}g=\Upsilon\cal{S}_M\cal{U}_M^{-1}g$ for every $g \in L^2(\sigma)$ such that $\lambda g\in L^2(\sigma)$, where $\sigma$ is the main spectral measure of $[M,L]$. The monodromy matrix of $\Hh_*$ has the form
\begin{equation}\label{loy1}
\left(\begin{smallmatrix}\Theta_{\ast}^+(\tau,z)&\Phi_{\ast}^+(\tau,z)\\\Theta_{\ast}^-(\tau,z)&\Phi_{\ast}^-(\tau,z)
\end{smallmatrix}\right)=\left(\begin{smallmatrix} \phi(\xi,z^2) &z\psi(\xi,z^2) \\  z^{-1}\phi'(\xi-,z^2) &\psi'(\xi-,z^2)
\end{smallmatrix}\right), \qquad \xi=N^{(-1)}(\tau),
\end{equation}
for details, see, e.g., Lemma 4.1 in \cite{KWW}.  Then, we obtain
\begin{align*}
(\Upsilon (\cal{U}_M^{-1}g))(\tau)
&=\left(\int_{\R_+}g(\lambda)\phi(N^{(-1)}(\tau),\lambda)\,d\sigma, 0\right)^t\\
&=\left(\frac 1\pi\int_{\R}g(x^2)\Theta^+_{\ast}(\tau,x)\,d\mu_\ast, 0\right)^t=\W^{-1}_{\Hh_{\ast}}(\pi^{-\frac 12}g(x^2)),
\end{align*}
using the fact that $g(x^2)$ is even, measure $\mu_\ast$ in \eqref{sa20} is even, and $\Theta_{\ast}^-(\tau,x)$ is odd in $x$. 
We see that 
\begin{equation}\label{lkj}
(\W_{\Hh_{\ast}}\Upsilon\,\, \cal{U}_M^{-1}g)(x)=\frac{g(x^2)}{\sqrt{\pi}}.
\end{equation}
Notice that, thanks to \eqref{sa20}, the map $g(\lambda)\mapsto g(\lambda^2)/\sqrt{\pi}$ is a unitary map of $L^2(\sigma)$ onto to the set of even functions in $L^2(\mu_\ast)$. Hence, $\Upsilon \cal{U}_M^{-1}g$ belongs to the domain of $\Di^2_{\Hh_{\ast}}$ if $g, \lambda g \in L^2(\sigma)$. Moreover, for $v=\cal{U}_M^{-1}g$ we have
\begin{equation}\label{des20}
\W_{\Hh_{\ast}}\Di^2_{\Hh_{\ast}}\Upsilon v=(\W_{\Hh_{\ast}}\Di^2_{\Hh_{\ast}}\W^{-1}_{\Hh_{\ast}})\W_{\Hh_{\ast}}\Upsilon v\stackrel{\eqref{lkj}}{=}x^2\cdot \frac{g(x^2)}{\sqrt\pi}.
\end{equation}
On the other hand, \eqref{lkj} also yields
\begin{equation}\label{lkj_bis}
\W_{\Hh_{\ast}}\Upsilon\,\, \cal{S}_Mv=x^2\cdot \frac{g(x^2)}{\sqrt\pi}.
\end{equation}
The comparison of \eqref{des20} and \eqref{lkj_bis} now gives $\Di^2_{\Hh_{\ast}}\Upsilon\cal{U}_M^{-1}g=\Upsilon\cal{S}_M\cal{U}_M^{-1}g$, as desired. 

\bigskip

\subsection{Wave equation for Krein strings}\label{weks}
The vibration of the proper string with parameters $[M,L]$  is governed by the following formal hyperbolic Cauchy problem:
\begin{equation}\label{e1}
u_{tt}+\cal{S}_Mu=0, \quad u(\xi,0)=u_0(\xi), \quad u_t(\xi,0)=0, \quad u'(0-,t)=0,
\end{equation}
where $u_0$ is the initial displacement of the string, its initial velocity is equal to zero, and the Neumann boundary condition $u'(0-,t)=0$ indicates that its left end is ``loose''. In this paper, we will only study solutions to \eqref{e1} given by the formula
\begin{equation}\label{e44}
u(\xi,t)=\cos (t\sqrt{\cal{S}_M})u_0,  \quad u_0\in L^2(\mf),
\end{equation}
where $\cos (t\sqrt{\cal{S}_M})$ is defined via spectral theorem:
\[
\cos (t\sqrt{\cal{S}_M})u_0=\cal{U}_M^{-1}\Bigl(\cos (t\sqrt \lambda)(\cal{U}_M u_0)(\lambda)\Bigr).
\]
Function $u$ is understood as element in $C(\R,L^2(\mf))\cap L^\infty(\R,L^2(\mf))$ in the standard mixed-norm notation. Let us notice that assumption $u_0\in {\rm dom }\,\cal{S}_M$ implies that $u$ is {\bf strong solution} which means that it is twice strongly continuously differentiable function of $t$ in the topology of Hilbert space $L^2(\mf)$ and that it satisfies equation for every $t\ge 0$ and initial conditions for $t=0$ (see \cite{Berezanskii}, p. 225). The uniqueness of such strong solution follows immediately from the self-adjointness of $\cal{S}_M$ (see \cite{Berezanskii}, Theorem 6.2 on p.~229). Assumption $u_0\in {\rm dom }\,\cal{S}_M$ for real-valued $u_0$ guarantees that the  energy is finite since
\begin{align*}
E(t)=\frac 12\int_{[0,L)} (u_t^2+(\cal{S}_Mu)u)\,d\mf=\frac 12\int_0^\infty \lambda(\cos^2(t\sqrt \lambda)+\sin^2(t\sqrt \lambda)) \cdot(\cal{U}_M u_0)^2d\sigma=\\
=\frac 12\|\sqrt{\cal{S}_M}u_0\|_{L^2(\mf)}^2<\infty.
\end{align*}
Thus, one can argue that real initial data $u_0\in {\rm dom}\, \cal{S}_M$ give rise to  solutions that make physical sense. Since $u$ is a linear operator of $u_0$, we can assume that $u_0$ is real when studying the dynamics of $u$.

\medskip

Using spectral theorem, formula \eqref{e44} can be rewritten as follows. Let $u$ be defined by \eqref{e44} and $U_t$ be an evolution for the canonical system with the Hamiltonian  $\Hh_{*} = \diag(\htwo,\hone)$ in which $\htwo$ and $\hone$ are obtained by formulas \eqref{sd1} and \eqref{sd2}. Then,
\begin{equation}\label{nw1}
\left(\begin{smallmatrix}u(\xi,t)\\0\end{smallmatrix}\right)=\frac 12\Bigl( (U_tX)(\tau)+ (U_{-t}X)(\tau) \Bigr),
\end{equation}
where $X(\tau)=\left(\begin{smallmatrix}u_0(N^{(-1)}(\tau))\\0\end{smallmatrix}\right), \,\xi=N^{(-1)}(\tau), \, \tau\ge 0$.

\medskip

\noindent{\bf Proof of Theorem \ref{tswf}.}
A change of variables in Lebesgue-Stieltjes integral gives
\begin{equation}\label{sa30}
\Ts(\xi)=\int_0^\tau \sqrt{\htwo(l)\hone(l)}dl=\Tcd(\tau), \quad \xi=N^{(-1)}(\tau),
\end{equation}
if the string $[M,L]$ and the Hamiltonian $\Hh_{*}$ are related as in Lemma \ref{l20}. The proof is now immediate from \eqref{nw1}, \eqref{sa30}, and Theorem~\ref{p0bb}. 
\qed
\begin{Rema}If the number $\Lc(t+a)$ in \eqref{add1} is an endpoint of indivisible interval, it must be its left endpoint by definition. Hence, Theorem \ref{p0bb} and the formula for  operator $\Upsilon$ show that $\|u(\cdot,t)\|_{L^2(\mf,\{\sff_t\})}=0$, i.e., the wavefront as a point never carries a positive $L^2(\mf)$--norm of solution.
\end{Rema}

The Theorem \ref{tswf} can be applied to many models. The vibration of the classical infinite Stieltjes string with beads of equal masses connected by massless wire exhibits infinite speed of propagation (see, e.g., \cite{TeschlBook}, p.25) and that example corresponds to $\rho=0$ and $\mf_{\bf s}=\sum_{j=0}^\infty \delta(\xi-j)$ where $\delta(\xi)$ denotes the unit point-mass at zero. The formula \eqref{ffwf} for the front of the wave confirms our intuition that the wave propagates instantaneously through the intervals on which $\rho=0$. In fact, it shows that the presence of nontrivial $\mf_{\bf s}$ on such an interval plays no role in that phenomenon. For example, if $\rho=1$ for $\xi\notin [\alpha,\beta], 0<\alpha<\beta<\infty$ and $\rho=0$ on $[\alpha,\beta]$, then the formula \eqref{ffwf} yields
\[
\sff_t=\left\{
\begin{array}{cc}
\sff_0+t, & t\le t_{\bf cr},\\
\beta+t, & t>t_{\bf cr}
\end{array}
\right., \quad t_{\bf cr}=\alpha-\sff_0
\]
as along as $\sff_0<\alpha$. Observe that we have $\sff_{t_{\bf cr}}=\alpha$ at the critical time $t_{\bf cr}$ since $\Ls$ is left-continuous.


\medskip

Recall that the measure
$\sigma = \upsilon\, m + \sigma_{\bf s}$ on $\R_+$ with the density $\upsilon$ and the singular part $\sigma_{\bf s}$  belongs to the Szeg\H{o} class $\szp$ if $(x+1)^{-1} \in L^1(\sigma)$ and 
$$
\int_{\R_+}\frac{\log \upsilon(x)}{\sqrt x(x+1)}\,dx > -\infty.
$$
A simple change of variables shows that $\sigma\in \szp$ if and only if $\mu_\ast\in \sz$, where $\mu_\ast$ is taken from Lemma~\ref{sal16}. In the Introduction, the class of strings for which the spectral measure is Szeg\H{o} was called $\szks$ and it was characterized in Theorem \ref{ch1sz}. The following result gives its dynamical description and has Theorem \ref{ts2bis} as a corollary.

\begin{Thm}\label{ts2}
Suppose $[M,L]$ is a proper string. If there is  $u_0 \in L_c^2(\mf)$ such that  
\begin{equation}\label{eqs14_cot}
\limsup_{t \to +\infty} \|u\|_{L^2(\mf, \dlto)} > 0, \qquad
\dlto := \bigl[\Ls(t+ a - \ell), \Ls(t + a)\bigr],
\end{equation}
then $[M,L]\in \szks$. Here,  $a=\Ts(\sff_0)$. Conversely, if $[M,L]\in \szks$, then 
\begin{equation}\label{eqs147}
\liminf_{t \to +\infty} \|u\|_{L^2(\mf, \dlto)} > 0
\end{equation}
for every $u_0\in L_c^2(\mf)$ not equal to zero identically and for every $\ell>0$. Here, again, $a=\Ts(\sff_0)$.
\end{Thm}

\beginpf That follows from Corollary \ref{t2bis} and \eqref{nw1}.\qed\medskip

\begin{Rema} Combined with the Theorem \ref{aps2} below, we conclude that \eqref{eqs147} can be strengthened to 
\[
\lim_{t \to +\infty} \|u\|_{L^2(\mf, \dlto)}> 0
\]
for every $\ell>0$.\end{Rema}

\smallskip

\noindent The solution to \eqref{e1} for the homogeneous string $[M,L]=[d\xi,\infty]$ is given explicitly via d'Alembert's formula: 
\begin{equation}\label{e2}
V_t^0u_0=\frac{u_0(\xi + t) + u_0(\xi - t)}{2},
\end{equation}
where $u_0$ is extended to $\R$ as even function. Note that $V_t^0u_0=\frac{u_0(\xi-t)}{2}+o(1)$, $t \to \infty$, with $o(1)$ in $L^2(\R_+)$. Hence, the evolution is equivalent to translation when $t\to +\infty$. For general strings, we need to introduce the modified dynamics. Given $y\in L^2(\R)$, we let
\begin{equation}\label{eq12s}
\widetilde V_{t}^0 y\colon \xi \mapsto \frac{1}{2}\chi_{\Omega_{\bf ac}(\mf)}(\xi)\cdot \rho^{-\frac 14}(\xi)\cdot y(\Ts(\xi)-t), \qquad \xi \in {[0, L)} .
\end{equation}
A change of variables gives $\|\widetilde V_{t}^0 y\|_{L^2(\mf)}\le \|y\|_{L^2(\R)}$. 
 If $\sigma\in \szp$ and $\sigma = \upsilon\,dx + \sigma_{\bf s}$, then its Szeg\H{o} function is defined by  
\[
D_{\sigma }(z):=\frac{{D}_{\mu_\ast}(\sqrt z)}{\sqrt\pi \sqrt[4] z}, \qquad z \in \C \setminus\R_{+},
\]
where the measure $\mu_\ast$, given by \eqref{sa20}, is the spectral measure of Hamiltonian $\Hh_{\ast}$. Notice that $|D_{\sigma }(\lambda)|^2=\upsilon(\lambda)$ for a.e. $\lambda>0$ in the sense of non-tangential boundary values. 
In the case when $[M,L]\in \szks$, we can obtain the asymptotics of $u$ near its wavefront. The following result implies Theorem \ref{aps2bis} from the Introduction. 

\begin{Thm}\label{aps2}Suppose $[M,L]\in \szks$.  Then, there is a map $u_0 \mapsto G_{u_0}$ from $L^2(\mf)$ to $L^2(\R)$, such that for every fixed positive $a$ we have
\begin{eqnarray}\label{e110}
\lim_{t\to +\infty}\|u-\widetilde V_{t}^0 G_{u_0}\|_{L^2(\mf,[\Ls(t-a), \Ls(t + a)])}= 0.
\end{eqnarray}
The function $G_{u_0}$ satisfies $\|G_{u_0}\|_{L^2(\R)}^2=2\|u_{0, \bf ac}\|^2_{L^2(\mf)}$, where $u_{0, \bf ac}$ is the orthogonal projection of $u_0$ to the absolutely continuous subspace $H_{\bf ac}(\cal{S}_M)$ of $\cal{S}_M$. Moreover, if $u_0\in H_{\bf ac}(\cal{S}_M)$, then
\begin{eqnarray}\label{e111}
\lim_{t\to +\infty}\|u-\widetilde V_{t}^0 G_{u_0}\|_{L^2(\mf)}= 0.
\end{eqnarray}
\end{Thm}
\beginpf We will prove that \eqref{e110} and \eqref{e111} hold with the following choice of the function $G_{u_0}$:
\begin{equation}\label{vopr1}
G_{u_0}(\eta)=\frac{1}{\sqrt\pi}\int_{\R_+}\frac{\Re(D_{\sigma }(\alpha) e^{i\eta\sqrt{\alpha}}) g_{\bf ac}(\alpha)}{ \sqrt[4]{\alpha}}\, d\alpha,\quad  g:=\cal{U}_Mu_0, \quad g_{\bf ac}:=g\cdot \chi_{\Omega_{\bf ac}(\sigma)}\,.
\end{equation}
Fix $a> 0$ and choose $b > a$. 
Define $\Hh_{\ast}$ by formulas \eqref{sd1}, \eqref{sd2} and note that  $\Hh_{\ast}\in \szcs$. Set $X=\Upsilon u_0$ and let, as above, $\xi = N^{(-1)}(\tau)$ for $\tau \ge 0$. Formula \eqref{nw1} and Theorem \ref{t1} give 
\begin{equation}\label{eq3}
\begin{pmatrix}u(\xi,t)\\0\end{pmatrix} = \frac 12\Bigl( (U_tX)(\tau)+ (U_{-t}X)(\tau) \Bigr) = \frac{1}{2}\Hh_{\ast,\mathfrak{n}}^{-\frac 12}(\tau)(U_t^0Y_{X,+}+U_{-t}^0Y_{X,-})(\Tcd(\tau)) + o(1), 
\end{equation}
as $t \to +\infty$, with $o(1)$ in $L^2(\Hh_{\ast},\dlt)$, where 
$$
\Hh_{\ast,\mathfrak{n}}^{-\frac 12}(\tau) = \begin{cases}\begin{pmatrix}
\left(\frac{\hone(\tau)}{\htwo(\tau)}\right)^{\frac 14} & 0 \\ 0 & \left(\frac{\htwo(\tau)}{\hone(\tau)}\right)^{\frac 14}
\end{pmatrix}, &\tau\colon\; \hone(\tau)\htwo(\tau) > 0,\\ 
0, & \tau\colon\; \hone(\tau)\htwo(\tau) = 0,
\end{cases}
$$
is the matrix from \eqref{refe7}, and $f=\W_{\Hh_{\ast}}X$,  $f_{\bf ac}=f\cdot \chi_{\Omega_{\bf ac}(\mu_\ast)}$, $\W_{\Hh_{0}}Y_{X, +} =  f_{\bf ac}\ov{D_{\mu_\ast}}$, $\W_{\Hh_{0}}Y_{X, -} =  f_{\bf ac} D_{\mu_\ast}$.  Applying \eqref{sa11}, we obtain
$$
\frac{1}{2}(U_t^0Y_{X,+}+U_{-t}^0Y_{X,-})(\tau)=
\frac{1}{2\sqrt 2} \left( \widehat h_+(\tau-t)   \left( \begin{matrix}1 \\ -i \end{matrix}\right) +\widecheck h_-(\tau-t)   \left( \begin{matrix}1 \\ {i} \end{matrix}\right)\right)+o(1),
$$
as $t \to +\infty$, with $o(1)$ in $L^2(\Hh_0)$ and $h_+ = f_{\bf ac}\overline{D_{\mu_\ast}}$, $h_- = f_{\bf ac}D_{\mu_\ast}$. In other words, we have
\begin{align*}
\frac{1}{2}(U_t^0Y_{X,+}+U_{-t}^0Y_{X,-})(\tau)=
\left(\begin{matrix}A_f(\tau-t)\\0\end{matrix}\right)+o(1), \qquad t \to +\infty,
\end{align*}
where $o(1)$ is in $L^2(\Hh_0)$ and 
\[
A_f(\eta)=\frac{1}{2\sqrt\pi}\int_{\R}   f_{\bf ac}(x)     \Re (D_{\mu_\ast}(x)e^{ix\eta}) \,dx, \qquad \eta \in \R.
\]
The formula  \eqref{eq3} implies 
\begin{equation}\label{eq04}
\begin{pmatrix}
u(\xi, t)\\
0
\end{pmatrix} 
= \Hh_{\ast,\mathfrak{n}}^{-\frac 12}(\tau)\left(\begin{matrix}A_f(\Tcd(\tau)-t)\\0\end{matrix}\right) + \Hh_{\ast,\mathfrak{n}}^{-\frac 12}(\tau) R_{1,t}(\Tcd(\tau)) + R_{2,t}(\tau),
\end{equation}
where $\|R_{1,t}\|_{L^2(\Hh_0)} + \|R_{2,t}\|_{L^2(\Hh, \dlt)} \to 0$ as $t \to +\infty$. Note that considering $\xi\in \Omega_{\bf ac}(\mf)$ for which $\rho(\xi)>0$ is the same as considering those $\tau \ge 0$ for which $\det \Hh_{\ast}(\tau)=\htwo(\tau)\hone(\tau)>0$. Moreover, for such $\xi$ and $\tau$ we have 
$$
\rho^{-\frac 14}(\xi)=\left(\frac{\htwo(\tau)}{\hone(\tau)}\right)^{\frac 14}, \qquad \Ts(\xi) = \Tcd(\tau).
$$
So, one can rewrite relation \eqref{eq04} in the form
$$
\begin{pmatrix}
u(\xi, t) - \chi_{\Omega_{\bf ac}(\mf)}(\xi) \rho^{-\frac 14}(\xi) A_f(\Ts(\xi) - t)\\
0
\end{pmatrix} 
= \Hh_{\ast,\mathfrak{n}}^{-\frac 12}(\tau) R_{1,t}(\Tcd(\tau)) + R_{2,t}(\tau),
$$
or in the form
$$
\Upsilon \bigl(u(\cdot, t) - \widetilde V_{t}^0 y\bigr)(\tau) = \Hh_{\ast,\mathfrak{n}}^{-\frac 12}(\tau) R_{1,t}(\Tcd(\tau)) + R_{2,t}(\tau),
$$
with $y = 2A_f$. Since $b > a$, the mapping $\Upsilon$ sends $L^2(\mf, [\Ls(t-a), \Ls(t + a)])$ into a subset of $L^2(\Hh, \dlt)$. Noting that 
$$
\|\Hh_{\ast,\mathfrak{n}}^{-\frac 12}\cdot (R_{1,t}\circ \Tcd )+ R_{2,t}\|_{L^2(\Hh, \dlt)} \le \|R_{1,t}\|_{L^2(\Hh_0)} + \|R_{2,t}\|_{L^2(\Hh, \dlt)} \to 0, \qquad t \to +\infty,
$$
we see that \eqref{e110} holds with $G_{u_0} = y = 2A_f$. If we put $g=\cal{U}_Mu_0$, then 
$g$ and $f$ are related to each other by $f(x)=g(x^2)/\sqrt\pi$ according to \eqref{lkj}. 
By construction, the function $f$ is even, and 
$D_{\mu_\ast}(-x) = \overline{D_{\mu_\ast}(x)}$ almost everywhere on $\R$ since $\mu_\ast$ is even. 
After changing variables $x^2=\alpha$, we obtain 
\begin{align*}
G_{u_0}(\eta)
&=\frac{1}{\sqrt\pi}\int_{\R}   f_{\bf ac}(x)    \Re (D_{\mu_\ast}(x)e^{ix\eta})\,dx
=\frac{1}{\sqrt{\pi}}\int_{\R_+}\frac{ g_{\bf ac}(\alpha)\Re (D_{\sigma }(\alpha)e^{i\sqrt{\alpha}\eta})}{\sqrt[4]{\alpha}}\,d\alpha,
\end{align*}
as in \eqref{vopr1}. We also have
$\|G_{u_0}\|_{L^2(\R)}^2= 4 \|A_f\|_{L^2(\R)}^2 = 2 \|  f_{\bf ac}\|_{L^2(\mu_\ast)}^{2} = 2 \|P_{\bf ac} X\|_{L^2(\Hh)}^{2} = 2 \|u_{0, \bf ac}\|_{L^2(\mf)}^{2} = 2\|g_{\bf ac}\|^2_{L^2(\sigma)}$, where in the second identity we have used the formula 
$$
A_f(\eta)=\frac{1}{2\sqrt\pi}\int_{\R} f_{\bf ac}(x)\Re (D_{\mu_\ast}(x)e^{ix\eta})\,dx = \frac{1}{2\sqrt\pi}\int_{\R} f_{\bf ac}(x)D_{\mu_\ast}(x)e^{ix\eta}\,dx,
$$ 
in which the integrals are understood in the $L^2(\R)$-sense. Finally, Corollary~\ref{c2} implies \eqref{e111}.\qed

\medskip

\subsection{Examples} \label{ssexamples} In this subsection, we explain how the general results can be applied to two examples considered in the Introduction.\medskip

\noindent {\bf Example \ref{ex2bis}: strings for which $\rho=1$.}
 Consider the case when $L=\infty$ and $\rho=1$. For the associated measure $\mf$, we get
\begin{equation}\label{e0}
d\mf=d\xi+d\mf_{\bf s},
\end{equation}
where $\mf_{\bf s}$ is any singular measure.  If $\mf_{\bf s}=0$, then the solution $u$ is given by \eqref{e2}.  The models described by our choice of $M$ are numerous, e.g., think about the beads with masses $\{m_j\}$ placed at points $\{\xi_j\}, \xi_0<\xi_1<\ldots$ connected by the string with a uniform density equal to one.\medskip

One might want to know how the presence of ``impurities'' encoded by $\mf_{\bf s}$ changes the character of wave propagation. The general results from the previous subsection can be reformulated as follows. From Theorem \ref{tswf}, we immediately get

\begin{Prop}If $u_0 \in L_c^2(\mf)$ is nonzero, then $\sff_t=\sff_0+t$.\label{ex1}
\end{Prop}
\smallskip

Clearly, the front propagates with the same linear speed regardless of the nature of $\mf_{\bf s}$. For $M$ that satisfy $\rho=1$, it was established (see \cite{BD2017}), that
\begin{equation}\label{s_a2}
\mf_{\bf s}(\R_+)<\infty \Longleftrightarrow \sigma\in \szp.
\end{equation}
In the next two statements, we describe how the dynamics of $u$ depends on $\mf_{\bf s}$.

\begin{Prop}\label{ap1}If $\|u_0\|_{L^2(\mf)}>0, \sff_0<\infty$, and
$\mf_{\bf s}(\R_+)=\infty$, then 
\[
\lim_{t\to\infty}\|u\|_{L^2(\mf,[\sff_t-a,\sff_t])}=0
\]
for every fixed $a>0$. Conversely, if there is $u_0$ that satisfies  $\|u_0\|_{L^2(\mf)}>0, \sff_0<\infty$, and
\[
\limsup_{t\to\infty}   \|u\|_{L^2(\mf,[\sff_t-a,\sff_t])}>0
\]
for some fixed $a>0$, then $\mf_{\bf s}(\R_+)<\infty$.\label{ex2}
\end{Prop}
\beginpf 
That follows from the Theorem \ref{ts2} and \eqref{s_a2}. \qed

\medskip

The next result shows that
the condition $\mf_{\bf s}(\R_+)<\infty$ guarantees that part of the wave propagates like a traveling wave in \eqref{e2}.  In that theorem, $P_{\bf ac}$ denotes the orthogonal projection to $H_{\bf ac}(\cal{S}_M)$.

\begin{Prop}\label{ap2}
In the case when $\mf_{\bf s}(\R_+) < \infty$,  we have 
\begin{align}
\lim_{t\to+\infty}\|u(\cdot,t)\|_{L^2(\mf_{\bf s}, [t-a,t+a])} = 0,\label{first1}\\
\lim_{t\to+\infty}\|u(\cdot,t)-G_{u_0}(\cdot-t)\|_{L^2[t-a,t+a]}= 0,\nonumber
\end{align}
for some  $G_{u_0} \in L^2(\R)$ and all $ a>0$.
Moreover, $\|G_{u_0}\|_{L^2(\R)}>0$ if and only if $P_{\bf ac}u_0\neq 0$. If $u_0$ is not identically equal to zero and has compact support, then $P_{\bf ac}u_0\neq 0$.
\end{Prop}
\beginpf These results are contained in Theorem \ref{ts2}, Theorem \ref{aps2}, and \eqref{s_a2}.\qed\smallskip

The statements made in Example \ref{ex2bis} in Introduction now follow.\bigskip

\noindent {\bf Example \ref{ex19}:  strings made of two types of material.} For another example, let us consider a string 
 $\mf = \rho\,d\xi$ on $\R_+$ with no singular part whose density $\rho$ takes two positive values: $a$ and $b$. Specifically,
\begin{equation}\label{eq68bis}
\rho(\tau) 
= \begin{cases}
a, & \tau \in E,\\
b, & \tau \in \R_+ \setminus E,
\end{cases}
\end{equation}
for some Lebesgue-measurable set $E\subseteq \R_+$. We interpret such strings as those made from two types of material. Despite the relative simplicity, the model when $\rho$ takes only two positive values can have a nontrivial spectrum, e.g., a spectrum with gap structure if $\rho$ is periodic (see also \cite{DLS2020} for analysis of related problems on the graphs).

\smallskip

In Example \ref{ex19} of Introduction, we claimed

\begin{Prop}Suppose $a\neq b$. We have $\sigma\in \szp$ if and only if either $|\{\xi\colon \rho(\xi)=a\}|<\infty$ or $|\{\xi\colon \rho(\xi)=b\}|<\infty$.\label{svk77}
\end{Prop}

\beginpf We will apply Theorem \ref{ch1sz} with properly chosen $\{\eta_n\}$. Define $\eta_n$ by 
\[
\eta_n=\int_0^n\sqrt{\rho(\xi)}\,d\xi.
\]
Thus, $\xi_n=n, n=0,1,\ldots$. Since $\rho$ takes values $a$ and $b$, we always have condition 
\[
0<C_1<\eta_{n+1}-\eta_n<C_2, \quad n=0,1,2,\ldots
\]
satisfied. For each $n \ge 0$, we have 
\begin{equation}\label{eq68}
\rho(\tau) 
= \begin{cases}
a, & \tau \in E_n,\\
b, & \tau \in F_n,
\end{cases}
\end{equation}
where $|E_n| = \delta_n, |F_n| = 1 - \delta_n$, and $E_n\subseteq [n,n+1), F_n\subseteq [n,n+1)$. Then,
\begin{align*}
2\int_{n}^{n+2}\rho(\tau)\,d\tau
&= 2a(\delta_n + \delta_{n+1}) + 2b(2-\delta_n - \delta_{n+1})\\
&= 4b + 2(a-b)(\delta_n + \delta_{n+1})
\end{align*}
and 
\begin{align*}
\left(\int_{n}^{n+2}\sqrt{\rho(\tau)}\,d\tau\right)^2 
&=
(2\sqrt{b} + 
(\sqrt{a} - \sqrt{b})(\delta_n + \delta_{n+1}))^2\\
&=4b + 4\sqrt{b}(\sqrt{a} - \sqrt{b})(\delta_n + \delta_{n+1}) + 
(\sqrt{a} - \sqrt{b})^2(\delta_n + \delta_{n+1})^2.
\end{align*}
If we denote 
$$
A_n = 2\int_{n}^{n+2}\rho(\tau)\,d\tau - \left(\int_{n}^{n+2}\sqrt{\rho(\tau)}\,d\tau\right)^2, \qquad n \ge 0,
$$
then  the straightforward calculation shows
\begin{align*}
A_n&=2(a-b)(\delta_n + \delta_{n+1}) - (4\sqrt{b}(\sqrt{a} - \sqrt{b})(\delta_n + \delta_{n+1}) + 
(\sqrt{a} - \sqrt{b})^2(\delta_n + \delta_{n+1})^2)\\
&=(\sqrt{a} - \sqrt{b})^2(2-\delta_n - \delta_{n+1})(\delta_n + \delta_{n+1}).
\end{align*}
Then, the string satisfies conditions in the left-hand side of \eqref{sdk7} if and only if either $a = b$ or
\begin{equation}\label{rev2_sa1}
\sum_{n=0}^\infty (2-\delta_n - \delta_{n+1})(\delta_n + \delta_{n+1}) < \infty.
\end{equation}
 Next,  if either 
\begin{equation}\label{nnk1}
|E|=\sum_{n=0}^\infty|E_n|=\sum_{n=0}^\infty \delta_n<\infty
\end{equation}
 or 
\begin{equation}\label{nnk2}
|E^c|=\sum_{n=0}^\infty|F_n|=\sum_{n=0}^\infty (1-\delta_n)<\infty,
\end{equation}
then \eqref{rev2_sa1} converges. Conversely, 
the convergence of the series \eqref{rev2_sa1} implies that \begin{equation}\lim_{n\to \infty}(2-\delta_n - \delta_{n+1})(\delta_n + \delta_{n+1})=0. \label{rev2_sa2}
\end{equation}
 Since $\delta_n\in [0,1]$ for each $n$, one can not have $\delta_l+\delta_{l+1}<\eps$ and $2-(\delta_{k}+\delta_{k+1})<\eps$ simultaneously if $|k-l|=1$ and $\eps<\frac 12$. Hence, \eqref{rev2_sa2} gives that
 either $\lim_{n\to\infty}\delta_n=0$ or $\lim_{n\to\infty}(1-\delta_n)=0$. In the former case, \eqref{rev2_sa1} is equivalent to \eqref{nnk1} and in the latter case, it is equivalent to \eqref{nnk2}.\qed

\begin{Rema}Taking $b=0$ in \eqref{eq68}, one gets a string made of one type of material which can be distributed with ``gaps'' over $\R_+$. The application of Theorem \ref{ch1sz} with $\eta_n=n, n=0,1,\ldots$ yields the similar result. Namely,  $\sigma\in \szp$ if and only if  $|\{\xi\colon \rho(\xi)=0\}|<\infty$. 
\end{Rema}

Indeed, in that case condition \eqref{sdk7} reads (recall that $\xi_n=\Ls(n)$ and $F=\{\xi\colon \rho(\xi)=0\}$)
\begin{equation}\label{amend1}
\sum_{n=0}^\infty \Big(2\sqrt a(\xi_{n+2}-\xi_n)-4\Big)<\infty.
\end{equation}
Since
\[
n=\sqrt a\int_0^{\xi_n}\chi_Ed\xi=\sqrt a (\xi_n-|[0,\xi_n]\cap F|),
\]
one has
$
\sqrt a|[0,\xi_n]\cap F|=\sqrt a\xi_n-n
$
so the sum in \eqref{amend1} converges if and only if 
$
\sum_{n=0}^\infty |[\xi_n,\xi_{n+2}]\cap F|<\infty.
$
Since $\lim_{n\to\infty}\xi_n=+\infty$, the last condition is equivalent to $|F|<\infty$. Notice that if $|F|=\infty$ in the last example, adding singular measure $\mf_{\bf s}$ can not place $\sigma$ in $\szp$. Indeed, inserting $\mf_{\bf s}$ does not change the grid $\{\xi_n\}$ but it increases $M(\xi_{n+2}) - M(\xi_n)$ in \eqref{sdk7} making the total sum diverge.

\medskip

\section{Dirac operators}\label{dir}

We start this section by recalling the definition of the one-dimensional Dirac operator. Then, we make the connection to canonical systems and explain how the results from the second section can be applied to prove the  theorems stated in Section \ref{s1point2}.

\subsection{Dirac operators}\label{s4}
\noindent Recall that the one-dimensional Dirac operator $\Dd_Q$ on $\R_+$ is defined by 
\begin{equation}\label{do}
\Dd_Q \colon Z \mapsto JZ' + Q Z, \qquad J = \jm, \qquad Q = \left(\begin{smallmatrix}q_1 & q_2\\ q_2  & -q_1\end{smallmatrix}\right). 
\end{equation}
Here the functions $q_1, q_2$ are real and belong to $L^1_{\loc}(\R_+)$. 
The ``free'' Dirac operator with potential $Q = 0$ will be denoted by $\Dd_0$.  The domain of $\Dd_Q$ is given by 
$$
\dom\Dd_{Q} = \left\{Z \in L^2(\R_+,\C^2)\colon\; 
\begin{aligned}
&Z \mbox{ is locally absolutely continuous on  }\R_+,\\ 
&JZ' + QZ  \in L^2(\R_+,\C^2),\\
&\langle Z(0), \zo\rangle_{\C^2} = 0.
\end{aligned}
\right\}.
$$
With this domain, the operator $\Dd_Q$ is a densely defined self-adjoint operator on $L^2(\R_+,\C^2)$, see Section~8.6 in \cite{LSb} or \cite{Malamud} for recent developments. Let $\Psi$ denote the generalized eigenvector of $\Dd_Q$: 
\begin{equation}\label{cp}
J\Psi'(\tau, z) + Q \Psi(\tau, z) = z \Psi(\tau,z), \quad \Psi(0, z) = \oz, \qquad \tau \ge 0, \quad z \in \C,
\end{equation} 
where the derivative is taken with respect to $\tau$. Then, there exists a unique Borel measure $\mu_D$ on $\R_+$ such that the generalized Fourier transform
\begin{equation}\label{ft}
\F_Q\colon Z \mapsto \frac{1}{\sqrt{\pi}}\int_{\R_+}\langle Z(\tau), \Psi(\tau, \bar z)\rangle_{\C^2}\,d\tau, \qquad z \in \C,
\end{equation} 
densely defined on $L_c^2(\R_+,\C^2)$, can be extended to a unitary operator from $L^2(\R_+,\C^2)$ to $L^2(\mu_D)$. That measure is called the main spectral measure of $\Dd_{Q}$.  

One can see that $\Dd_0$ and $\Di_{\Hh_0}$ are the same operators acting on the same Hilbert space $L^2(\R_+,\C^2)$ and giving rise to identical generalized eigenvectors, generalized Fourier transforms, and main spectral measures (cf. \eqref{cp}, \eqref{ft} and \eqref{cs}, \eqref{wt}).

\medskip

\subsection{The reduction of Dirac operator to a canonical system and the Szeg\H{o} condition}
The following result is well-known, see, e.g., Section 2.4 in \cite{B2018}. 
\begin{Lem}\label{ueq}
Let $Q \in L^1_{\loc}(\R_+)$ be as in \eqref{do}, and let the matrix-valued function $N_0$ be the solution of the Cauchy problem
\begin{equation}\label{s_a3}
JN_0'(\tau) + Q(\tau)N_0(\tau) = 0, \quad \tau \ge 0, \quad N_0(0) = \idm. 
\end{equation}
Denote by $\Di_{\Hh}$ the self-adjoint operator on $H$ corresponding to the canonical system generated by the Hamiltonian $\Hh = N_0^*N_0$. Then, the main spectral measures of the operators $\Dd_{Q}\colon L^2(\R_+,\C^2) \to L^2(\R_+,\C^2)$ and $\Di_{\Hh}\colon H \to H$  coincide. In particular, the operators $\Dd_{Q}$ and $\Di_{\Hh}$ are unitary equivalent and the unitary equivalence is given by the operator  $V\colon X\mapsto N_0 X$, which is a unitary map from $H$ to $L^2(\R_+,\C^2)$. Moreover, we have $\W_{\Hh} = \F_{Q}V$. 
\end{Lem}
An important property of the locally absolutely continuous Hamiltonian $\Hh = N_0^*N_0$ in the previous lemma is that it has unit determinant everywhere on $\R_+$. Indeed, the Wronskian in problem \eqref{s_a3} is constant so $\det N_0(\tau) = \det N_0(0) = 1$. Hence, 
\begin{equation}\label{eq59}
\det\Hh(\tau) = 1, \quad \Tc(\tau) = \tau, \qquad \Lc(\eta)=\eta \qquad \tau,\eta \in \R_+,
\end{equation}
for the corresponding functions $\Tc$ and $\Lc$ introduced in \eqref{sa21}.  
The identity $\det N_0(\tau)=1$ has other important implications. First, in the polar decomposition  $N_0=O|N_0|$ 
the matrix $|N_0|$ satisfies
\[
|N_0|=\sqrt{N_0^*N_0}=\Hh^{\frac 12}, \quad \det |N_0|=\det \Hh^{\frac 12}=(\det \Hh)^{\frac 12}=1
\]
and the rotation matrix $O$ is defined uniquely and is locally absolutely continuous. 
Second, the space $H$ coincides with $L^2(\Hh)$ defined in \eqref{eq31}. Moreover, if function $X$ is compactly supported, then $VX$ is also compactly supported and their supports coincide.

\medskip

The spectral measures $\mu_D$ of the Dirac operators define a proper subset of all Poisson-finite measures on the real line that generate the canonical systems as discussed in the second section. Some of them belong to the Szeg\H{o} class. Thanks to the characterization~\eqref{s_a5}, we have the following proposition (see Corollary~1.4 in \cite{BD2019}):
\begin{Prop}\label{p41}
The condition $N_0^*N_0\in \szcs$ is necessary and sufficient for the spectral measure $\mu_D$ of the Dirac operator to satisfy $\mu_D\in \sz$.
\end{Prop}
Checking that $N^*_0N_0\in \szcs$ is not always easy. However, in many cases, the application of our proposition is straightforward. 
\begin{Prop}\label{prop3}
Suppose that the potential $Q$ with entries in $L^1_{\loc}(\R_+)$ has the form 
\begin{equation}\label{eq60}
 Q = \begin{pmatrix}q & 0\\ 0 & -q\end{pmatrix} \quad \mbox{or}\quad Q = \begin{pmatrix}0 & q\\ q & 0\end{pmatrix}.
\end{equation}
Then, for the spectral measure $\mu_D$ of the corresponding Dirac operator $\Dd_Q$ we have
\begin{equation}\label{eq61}
\mu_D\in \sz\quad  \Longleftrightarrow \quad N_0^*N_0\in \szcs\quad \Longleftrightarrow\quad  \sum_{n \ge 0}\left(\int_{n}^{n+2}h(\tau)\,d\tau\int_{n}^{n+2}\frac{d\tau}{h(\tau)} - 4\right) < \infty,
\end{equation}
where $h(\tau) = e^{2\int_{0}^{\tau}q(s)\,ds}$, $\tau \ge n$.
\end{Prop}
\beginpf The first equivalence has already been discussed and we need to show the second one.  For the potentials $Q$ of the form \eqref{eq60}, define $g(\tau)=\int_{0}^{\tau}q(s)\,ds$. Then, solving the problem \eqref{s_a3} to find $N_0$ is easy.  That gives 
$$
N_0(\tau) = \begin{pmatrix}\cosh g(\tau) & \sinh g(\tau)\\\sinh g(\tau) & \cosh g(\tau)\end{pmatrix} \quad \mbox{or} \quad
N_0(\tau) = \begin{pmatrix}e^{-g(\tau)} & 0\\0 & e^{g(\tau)}\end{pmatrix},
$$
respectively. Therefore, for $\Hh=N_0^*N_0$, we have
$$
\Hh(\tau) = \begin{pmatrix}\cosh 2g(\tau) & \sinh 2g(\tau)\\\sinh 2g(\tau) & \cosh 2g(\tau)\end{pmatrix} \quad 
\quad \mbox{or} \quad
\Hh(\tau) = \begin{pmatrix}e^{-2g(\tau)} & 0\\0 & e^{2g(\tau)}\end{pmatrix}.
$$
In both cases $\det \Hh = 1$ on $\R_+$, and the straightforward calculation yields
$$
\det\int_{n}^{n+2}\Hh(\tau)\,d\tau = \int_{n}^{n+2}e^{2g(\tau)}\,d\tau\int_{n}^{n+2}e^{-2g(\tau)}\,d\tau = \int_{n}^{n+2}h(\tau)\,d\tau\int_{n}^{n+2}\frac{1}{h(\tau)}\,d\tau.
$$
So, $\Hh \in \szcs$ if and only if \eqref{eq61} holds. The result follows. \qed

\medskip

\begin{Rema}
Propositions \ref{p41} and \ref{prop3} imply Theorem \ref{t3-i}.
\end{Rema}

\medskip

\subsection{The evolution for Dirac equation and M{\o}ller wave (modified wave) operators}
The self-adjoint operator $\Dd_Q$ defines a unitary evolution $e^{it\Dd_Q}$. Lemma \ref{ueq} above gives the connection between $e^{it\Dd_Q}$ and evolution $e^{it\Di_{\Hh}}$ for canonical systems. 
\begin{Lem}If $\Hh = N_0^*N_0$, then 
$e^{it\Dd_Q}=Ve^{it\Di_{\Hh}}V^{-1}$
for all $t\in \R$.\label{reduc1}
\end{Lem}

For every $Z$ in $L^2(\R_+,\C^2) = L^2(\Hh_0)$ we again define the {\bf front} as
$$
\sff[Z]=\inf\{\tau \ge 0\colon\; Z(s) = 0 \; \mbox{for almost every }s > \tau\}.
$$
Notice that all elements of $N_0$ are real-valued so $Z$ has real components if and only if $X=V^{-1}Z$ has real components.
\begin{Prop}\label{p01}
For every real $Z \in L_c^2(\R_+,\C^2)$ and every $t\in \R$, we have 
$$\sff[e^{it\Dd_Q}Z] = |t|+\sff[Z].$$
\end{Prop}  
\beginpf That follows from Theorem \ref{p0bb} and formula \eqref{eq59}.\qed

\begin{Rema} For an arbitrary $Z \in L_c^2(\R_+,\C^2)$, we can write each of its components as a sum of real and imaginary parts. Then, Proposition \ref{p01} gives
$\sff[e^{it\Dd_Q}Z] \le |t|+\sff[Z]$ for all $t \in \R$. 
\end{Rema}

\medskip

We will also need the following proposition. \begin{Prop}\label{p011}
Let $\mu$ be the spectral measure of the Dirac operator $\Dd_Q$. Suppose $\mu\notin \sz$. Then, for every $Z\in L^2(\R_+,\C^2)$, we have
\[
\lim_{t\to\pm\infty} \|e^{it\Dd_Q}Z\|_{L^2(\C^2, [|t|-b,|t|+b])}=0
\]
for every $b>0$.
\end{Prop}  
\beginpf That follows from Corollary \ref{t2bis1} and formula \eqref{eq59}.\qed

\medskip

\noindent{\bf Proof of Theorem \ref{t4-i}.} Take a nonzero function $Y \in L^2(\R_+,\C^2)$. Set
$$
\Zplus=W_+(\Dd_{Q}, \Dd_0, \gamma)Y=\lim_{t \to +\infty} e^{-it\Dd_Q}M_{\gamma} e^{it\Dd_0}Y,
$$ 
and notice that $\|\Zplus\|_{L^2(\R_+,\C^2)}=\|Y\|_{L^2(\R_+,\C^2)}>0$.
That yields
\begin{equation}\label{s_a7}
\lim_{t \to +\infty} \| M_{\gamma}e^{it\Dd_0}Y-e^{it\Dd_Q}\Zplus\|_{L^2(\R_+,\C^2)}=0.
\end{equation}
Formula \eqref{sa12} for $Y = \begin{pmatrix}Y_1 \\ Y_2\end{pmatrix}$ can be recast as
\begin{equation}\label{eq32bisbis}
e^{it\Dd_{0}}\begin{pmatrix}Y_1 \\ Y_2\end{pmatrix} 
= 
\frac{1}{2}
\begin{pmatrix}
Y_1(\tau-t) + Y_1(\tau+t)\\
-i(Y_1(\tau-t) - Y_1(\tau+t))
\end{pmatrix}
+
\frac{1}{2}
\begin{pmatrix}
i(Y_2(\tau-t) - Y_2(\tau+t))\\
Y_2(\tau-t) + Y_2(\tau+t)
\end{pmatrix}, \qquad \tau \in \R_+,
\end{equation}
where $Y_1 \in L^2(\R_+)$  is extended to the whole real line $\R$ as an even function and $Y_2$ is extended as odd function.
That gives
$$\liminf_{t\to +\infty}\|e^{it\Dd_0}Y\|_{L^2(\C^2,[t-b,t+b])}>0$$ for large enough $b$. From \eqref{s_a7}, one gets $\liminf_{t\to +\infty}\|e^{it\Dd_Q}\Zplus\|_{L^2(\C^2,[t-b,t+b])}>0$. Now we have $\mu\in \sz$ by Proposition \ref{p011}. The case when $t\to -\infty$ can be handled similarly. \qed

\medskip

\begin{Rema}
Notice that we have used the existence of $\lim_{t \to +\infty} e^{-it\Dd_Q}M_{\gamma} e^{it\Dd_0}Y$ for just one nonzero element $Y \in L^2(\R_+,\C^2)$ in the proof of Theorem \ref{t4-i}.
\end{Rema}

\medskip

\noindent {\bf Proof of Theorem \ref{t5-i}.} Assume that the main spectral measure $\mu_D$ of $\Dd_{Q}$ is in the Szeg\H{o} class. By Lemma \ref{ueq}, $\mu_D$ coincides with the spectral measure of the Hamiltonian $\Hh = N_0^* N_0$ generated by the solution of equation $J N_0' + QN_0 = 0$, $N_0(0) = \idm$. Taking into account \eqref{eq59}, Theorem \ref{t3} tells us that for some function $\gamma_0\colon \R_+ \to \T$ and for every $Y\in L^2(\Hh_0)=L^2(\R_+,\C^2)$,  the limits
$$
\lim_{t\to +\infty}e^{-it\Di_{\Hh}}M_{\gamma_0}\Hh^{-\frac 12}e^{it\Di_{\Hh_0}}Y, \quad \lim_{t\to -\infty}e^{-it\Di_{\Hh}}M_{\overline{\gamma}_0}\Hh^{-\frac 12}e^{it\Di_{\Hh_0}}Y
$$
exist in the norm of $L^2(\Hh)$.  Since  $\Di_{\Hh_0}= \Dd_0$, that implies existence of the limit 
$$
\lim_{t\to +\infty}Ve^{-it\Di_{\Hh}}M_{\gamma_0}\Hh^{-\frac 12}e^{it\Di_{\Hh_0}}Y
=
\lim_{t\to +\infty}e^{-it\Dd_{Q}}VM_{\gamma_0}\Hh^{-\frac 12}e^{it\Dd_{0}}Y
$$
in $L^2(\R_+,\C^2)$. 
Note that $(N_0(\tau)\Hh^{-\frac 12}(\tau))^* (N_0(\tau)\Hh^{-\frac 12}(\tau))$ is the identity matrix for each $\tau \in \R_+$. Since \[
\det (N_0(\tau)\Hh^{-\frac 12}(\tau))=1, \qquad \tau \in \R_+,\] it follows that the operator $Y \mapsto V\Hh^{-\frac 12}Y$ on $L^2(\R_+,\C^2)$ coincides with the multiplication operator by a $2 \times 2$ matrix-valued function of the form 
$$
\Sigma_{\phi(\tau)} = \begin{pmatrix}\cos\phi(\tau) & \sin\phi(\tau) \\ -\sin\phi(\tau) & \cos\phi(\tau)\end{pmatrix}, \qquad \phi(\tau) \in [0, 2\pi).
$$  
We have
$$
\Sigma_{\phi}\begin{pmatrix}1 \\ -i\end{pmatrix}  
= \begin{pmatrix}\cos\phi - i\sin\phi \\ -\sin\phi - i\cos\phi\end{pmatrix} 
= 
e^{-i\phi}\begin{pmatrix}1 \\ -i\end{pmatrix}. 
$$
Formula \eqref{eq32bisbis} shows that when $t \to +\infty$, for every $Y\in L^2(\R_+,\C^2)$ we have
$$
e^{it\Dd_{0}}\begin{pmatrix}Y_1 \\ Y_2\end{pmatrix} 
= 
\frac{Y_1(\tau-t)}{2}
\begin{pmatrix}
1\\
-i
\end{pmatrix}
+
\frac{iY_2(\tau-t)}{2}
\begin{pmatrix}
1\\
-i
\end{pmatrix} + o(1),
$$
where $o(1)$ is with respect to $L^2(\R_+,\C^2)$--norm. Therefore, 
$$e^{-it\Dd_{Q}}VM_{\gamma_0}\Hh^{-\frac 12}e^{it\Dd_{0}}Y = e^{-it\Dd_{Q}}M_{\gamma_0e^{-i\phi}}e^{it\Dd_{0}}Y + o(1), \qquad t \to +\infty$$ and  the limit 
$$
W_+(\Dd_{Q}, \Dd_0, \gamma)Y=\lim_{t\to+\infty}e^{-it\Dd_{Q}}M_{\gamma}e^{it\Dd_{0}}Y
$$
exists in $L^2(\R_+,\C^2)$ for all $Y \in L^2(\R_+,\C^2)$ if we take $\gamma =\gamma_0 e^{-i\phi}$. The existence of the wave operator $W_{+}(\Dd_Q, \Dd_0, \gamma)$ follows. Arguing similarly, one can prove the existence of $W_{-}(\Dd_Q, \Dd_0, \gamma)$ with the modification to the dynamics given by  $M_{\overline{\gamma}}$. Moreover, the proof shows that  
$$
\Ran W_{\pm}(\Dd_Q, \Dd_0,\gamma) = V(\Ran W(\Di_{\Hh}, \Di_{\Hh_0},\gamma)) = \F_{Q}^{-1}\W_{\Hh}(H_{\bf ac}(\Dd_{\Hh})) =  \F_{Q}^{-1}(L^2(\mu_{\bf ac})) = H_{\bf ac}(\Dd_{Q}),
$$ 
where $\mu_{\bf ac} = w\,dx$ is the absolutely continuous part of the measure $\mu$ and $\mu=\mu_D$. In other words, the wave operators $ W_{\pm}(\Dd_Q, \Dd_0,\gamma)$ are complete. It is also clear from the proof that our construction gives $\gamma = 1$ in the case where $Q$ is anti-diagonal ($q_1 = 0$). \qed

\medskip

\noindent{\bf Proof of Theorem \ref{t05-i}.} Given Lemma \ref{reduc1}, Theorem \ref{t05-i} is a direct consequence of Theorem \ref{t4} and Theorem \ref{c3}. \qed

\medskip

\noindent{\bf Proof of Theorem \ref{t7-i}.} Given Lemma \ref{reduc1},  Theorem \ref{t7-i} follows from  Theorem \ref{t4}. \qed

\subsection{Wiegner-von Neumann potentials}\label{wfp}
In this subsection, we prove Theorem \ref{t6-i}. Let us recall its statement for convenience:

\medskip

\noindent {\it 
Suppose that a potential $Q \in L^1_{\loc}(\R_+)$ has the form \eqref{eq60}, with $q = \frac{\sin \tau^\alpha}{\tau^\beta}$ 
on $[\tau_0, +\infty)$ for some $\tau_0 > 0$ and $\alpha, \beta \in \R$. Then, $Q \in \szd$ if and only $(\alpha, \beta) \in A$, where 
$$
A= \{\alpha \le 0,\; \beta-\alpha > \frac 12\} \cup  \{\alpha \in (0,1),\; \beta > \frac 12\} \cup \{\alpha \ge 1, \;\alpha+\beta > \frac 32\}
$$ 
is the open set depicted on Figure~\ref{fig:plot-i}.
}

\medskip

We will need a variant of Korey's estimate from \cite{Korey}. Recall that we use notation $\langle f \rangle_{I} = \frac{1}{|I|}\int_{I}f(x)\,dx$.
\begin{Lem}
Suppose $I=[a,b]$ and  measurable function $f\colon I \to \R$ satisfies
\begin{equation}\label{des11}
\langle e^{f}\rangle_{I} \cdot \langle e^{-f}\rangle_I = 1 + \eps, \qquad \eps \in [0, 1].
\end{equation}
Then, 
\begin{equation}\label{eq63}
\langle \left|f - \langle f\rangle_I\right|\rangle_I\le c\sqrt{\eps}, 
\end{equation}
for a universal constant $c$.
\end{Lem}
\beginpf In \cite{Korey}, formula (3.7), Korey shows that 
$$
\langle e^{f}\rangle_{I} \cdot  \exp\left(-\langle f\rangle_I\right) = 1 + \eps', \qquad \eps' \in [0, 1],
$$
implies 
$
\langle |f-m_I(f)|\rangle\le C\sqrt{\eps'},
$
where $m_I(f)$ denotes a median of $f$ over $I$. By Jensen's inequality, 
$
\exp\left(-\langle f\rangle_I\right)\le  \langle e^{-f}\rangle_I.
$
Therefore, \eqref{des11} implies
\[
\langle |f-m_I(f)|\rangle\le C\sqrt\eps.
\]
Now, for every $c' \in \R$, we have
\begin{align*}
\langle \left|f - \langle f\rangle_I\right|\rangle_I\le \langle \left|f - c'\right|\rangle_I+|c'-\langle f\rangle_I|
\le 2\langle \left|f - c'\right|\rangle_I.
\end{align*}  
Taking $c'=m_I(f)$ finishes our proof.\qed

\medskip

For integer $n \ge 0$, we let $I_n = [n, n+2]$. Given real-valued $q\in L^1_{\rm loc}(\R_+)$, denote  $g_n(\tau) = 2\int_{n}^{\tau}q(s)\,ds$.
\begin{Prop}\label{prop4}
If   $\lim_{n \to +\infty} \|g_n\|_{L^\infty(I_n)} = 0$, then  
\begin{equation}\label{des16}
 \sum_{n \ge 0}\left(\int_{I_n}e^{g_n}\,d\tau\int_{I_n}e^{-g_n}\,d\tau - 4\right) < \infty\quad  \Longleftrightarrow\quad 
\sum_{n \ge 0}\int_{I_n} |g_n-\langle g_n\rangle_{I_n}|^2\,d\tau < +\infty.
\end{equation}
\end{Prop}

\beginpf Set $\widetilde g_n = g_n-\langle g_n\rangle_{I_n}$ and notice that
\[\langle e^{g_n}\rangle_{I_n} \cdot \langle e^{-g_n}\rangle_{I_n}
=\langle e^{\widetilde g_n}\rangle_{I_n} \cdot \langle e^{-\widetilde g_n}\rangle_{I_n}.
\]
 Then, since $\lim_{n\to\infty}\|\widetilde g_n\|_{L^\infty(I_n)}=0$, we use Taylor expansion to get
\begin{align*}
\int_{I_n} e^{\pm\widetilde g_n}\,d\tau 
&= 2 \pm \int_{I_n}\widetilde g_n\,d\tau + \frac{1}{2}\int_{I_n}\widetilde g^2_n\,d\tau + O\left(\int_{I_n}|\widetilde g_n|^3\,d\tau\right),
\end{align*}
as $n \to +\infty$. It follows that
$$
\int_{I_n}e^{\widetilde g_n}\,d\tau\int_{I_n}e^{-\widetilde g_n}\,d\tau =
4 + 2 \int_{I_n}\widetilde g^2_n(\tau)\,d\tau + o\left(\int_{I_n}\widetilde g^2_n(\tau)\,d\tau\right),
$$
which proves the required claim. \qed\medskip

\medskip

\noindent{\bf Proof of Theorem \ref{t6-i}.}  Recall Proposition \ref{prop3}, Theorems \ref{t4-i} and \ref{t5-i}. To prove our result, we only need to establish the range of  parameters $\alpha$ and $
\beta$ for which the condition
\begin{equation}\label{help}
\sum_{n \ge 0}\left(\int_{I_n}e^{g_n}\,d\tau\int_{I_n}e^{-g_n}\,d\tau - 4\right) < \infty
\end{equation}
is satisfied, where $g_n$ is defined right before the Proposition \ref{prop4}.  For an integer $n \ge \max(1, \tau_0)$ and $\gamma = (\alpha+\beta-1)/\alpha$, we have
\begin{equation}\label{eq100}
\frac{\alpha g_n(x)}{2} = \alpha\int_{n}^{x}\frac{\sin \tau^\alpha}{\tau^\beta}\,d\tau 
= \int_{n^\alpha}^{x^\alpha}\frac{\sin y}{y^{\gamma}}\,dy =\left. -\frac{\cos y}{y^{\gamma}}\right|_{n^\alpha}^{x^\alpha} -\gamma \left. \frac{\sin y}{y^{\gamma +1}}\right|_{n^\alpha}^{x^\alpha} -
\gamma(\gamma+1)
\int_{n^\alpha}^{x^\alpha}\frac{\sin y}{y^{\gamma+2}} \,dy.
\end{equation}
Let us consider several cases.

\smallskip

\noindent {\bf Case \boldmath$\alpha \ge 1$, \boldmath$\alpha + \beta > 1$.} In this case $\gamma > 0$ and we are in the setting of Proposition \ref{prop4}. Note that 
$$
\sup_{x\in I_n}\left|
-\gamma \left. \frac{\sin y}{y^{\gamma +1}}\right|_{n^\alpha}^{x^\alpha} -\gamma(\gamma+1) \int_{n^\alpha}^{x^\alpha}\frac{\sin y}{y^{\gamma+2}} \,dy \right|= O\left(\frac{1}{n^{\alpha(\gamma+1)}}\right)\in \ell^2(\mathbb{N}).
$$
Thus, we only need to control the sum of dispersions 
\begin{equation}\label{des18}
\sum_{n\ge 1}\Bigl\langle  \left|\frac{\cos x^\alpha}{x^{\alpha\gamma}}-   \Bigl\langle\frac{\cos x^\alpha}{x^{\alpha\gamma}}\Bigr\rangle_{I_n} \right|^2\Bigr\rangle_{I_n}=\sum_{n\ge 1}\left(\Bigl\langle  \left|\frac{\cos x^\alpha}{x^{\alpha\gamma}}\right|^2\Bigr\rangle_{I_n}-  \Bigl\langle  \frac{\cos x^\alpha}{x^{\alpha\gamma}}\Bigr\rangle_{I_n}^2  \right).
\end{equation}
Set $\eta = \alpha\gamma + \alpha - 1$, and note that $\eta \ge \gamma > 0$. Similarly to \eqref{eq100}, we have
\[
\int_{I_n}\frac{\cos (x^\alpha)}{x^{\alpha\gamma}}\,dx=\frac{1}{\alpha}\left(\frac{\sin(n+2)^\alpha}{(n+2)^{\eta}}-\frac{\sin n^\alpha}{n^{\eta}}\right)+O(n^{-\eta-\alpha})
\]
and $n^{-\eta-\alpha}\in \ell^2(\mathbb{N})$. So, the question reduces to the convergence of the series
\begin{equation}\label{eq62}
\sum_{n \ge 1}\left(2\alpha^2\int_{I_n}\frac{\cos^2 x^\alpha}{x^{2\alpha\gamma}}\,dx -  \left(\frac{\sin(n+2)^\alpha}{(n+2)^{\eta}}-\frac{\sin n^\alpha}{n^{\eta}}\right)^2  \right).
\end{equation}
Notice first that  $\alpha\ge 1, \alpha+\beta>\frac 32$ implies $2\alpha\gamma>1$ and $2\eta>1$. We get convergence in that situation and \eqref{help} holds.

\smallskip

We claim that for $\alpha>1, 1<\alpha+\beta\le \frac 32$, the series diverges and \eqref{help} does not hold.
Indeed, in that case, one has
\[
\int_{I_n}\frac{\cos^2 x^\alpha}{x^{2\alpha\gamma}}\,dx\sim n^{-2\alpha\gamma}, \quad n^{-2\eta}=o(n^{-2\alpha\gamma}), \quad n\to\infty,
\]
and we get the claim because $2\alpha\gamma\le 1$. \smallskip

It is only left to consider $\alpha=1$ and $0<\beta\le \frac 12$. In that case, $\gamma=\beta$, and
\[
\frac{\cos x}{x^\gamma}-\Bigl\langle \frac{\cos x}{x^\gamma}\Bigr\rangle_{I_n}= n^{-\beta}(\cos x-\langle \cos x\rangle_{I_n})+O(n^{-\beta-1}).
\]
Since $\langle (\cos x-\langle \cos x\rangle_{I_n})^2\rangle_{I_n}\sim 1$ and $\beta\le \frac 12$, the series \eqref{des18} diverges.\smallskip

To summarize, if $\alpha \ge 1$ and $\alpha + \beta > 1$, then \eqref{help} holds if and only if $\alpha + \beta > \frac 32$.

\smallskip

\noindent {\bf Case \boldmath$\alpha \ge 1$,  \boldmath$\alpha +\beta \le 1$.} We will show that \eqref{help} {fails} for this range of parameters. Suppose, on the contrary, that \eqref{help} holds. Then, estimate \eqref{eq63} gives 
\begin{equation}\label{eq64}
\sum_{n \ge 1}\left(\int_{I_n}\left|g_n-\langle g_n\rangle_{I_n}\right|\,dx\right)^2 < \infty. 
\end{equation}
Recall that $\gamma = (\alpha+\beta-1)/\alpha\le 0$. Integration by parts gives
$$
\frac{\alpha g_n(x)}{2} = \alpha\int_{n}^{x}\frac{\sin \tau^\alpha}{\tau^\beta}\,d\tau
= n^{-\alpha\gamma}\cos n^{\alpha} - x^{-\alpha\gamma}\cos x^{\alpha} + O(n^{-\alpha(\gamma + 1)}).
$$
and
$$
\frac{\alpha}{2}\langle g_n\rangle_{I_n}
= n^{-\alpha\gamma}\cos n^{\alpha}
+ O(n^{-\alpha(\gamma + 1) + 1}).
$$
Hence, uniformly in $x \in I_n$, we have
$$
\frac{\alpha}{2}\left|g_n(x) - \langle g_n\rangle_{I_n} \right| = 
x^{-\alpha\gamma}|\cos x^{\alpha}| + O(n^{-\alpha(\gamma + 1) + 1}).
$$
Since 
$$
\int_{I_n} x^{-\alpha\gamma}|\cos x^{\alpha}|\,dx + O(n^{-\alpha(\gamma + 1) + 1}) \ge n^{-\gamma\alpha}\left(\int_{I_n} |\cos x^{\alpha}|\,dx + O(n^{-\alpha + 1})\right) 
$$
and 
\begin{equation}
\inf_{n}\int_{I_n} |\cos x^{\alpha}|\,dx \ge C>0
\end{equation}
for every $\alpha \ge 1$, the terms in series \eqref{eq64} can be estimated from below by $c n^{-2\gamma\alpha} \ge c > 0$ for large $n$. So, \eqref{help} does not hold if $\alpha \ge 1$ and $\alpha +\beta \le 1$. 

\smallskip

\noindent {\bf Case \boldmath$\alpha \in (0,1)$, \boldmath$\beta \le \frac 12$.} We are going to show that \eqref{help} fails for this range of parameters. Assume, as before, that \eqref{help} holds, so that the sum in \eqref{eq64} is finite. Applying the mean-value theorem, we see that
$$
\sum_{n \ge 1}\left(\int_{I_n}\left|g_n(x) - g_n(x_n)\right|\,dx\right)^2 < \infty
$$
for some points $x_n \in I_n$. For each $n \ge 0$, one can choose an interval $\Delta_n \subset I_n$ of length $\frac 14$ such that $\dist(x_n, \Delta_n) \ge \frac 14$. Then, again by the mean-value theorem (this time -- on the interval $\Delta_n$), there exist points $\widetilde x_n \in \Delta_n$ such that
$$
\sum_{n \ge 1}\left|g_n(\widetilde x_n) - g_n(x_n)\right|^2 \lesssim \sum_{n \ge 1}\left(\int_{I_n}\left|g_n(x) - g_n(x_n)\right|\,dx\right)^2 < \infty.
$$
We see that
\begin{equation}\label{s_a20}
\sum_{n \ge 1} \left|\int_{x_n}^{\widetilde x_n}\frac{\sin \tau^\alpha}{\tau^\beta}\,d\tau\right|^2 < \infty. 
\end{equation}
Since $\sin \tau^{\alpha} - \sin n^\alpha = O(n^{\alpha - 1})$ for $\tau \in I_n$, we have 
$$
\left|\int_{x_n}^{\widetilde x_n}\frac{\sin \tau^\alpha}{\tau^\beta}\,d\tau\right| \ge |\sin n^\alpha| \cdot J_n - O( n^{\alpha-1}J_n), \qquad J_n = \left|\int_{x_n}^{\widetilde x_n}\frac{d\tau}{\tau^\beta}\,\right|.
$$
Given that $\beta\le \frac 12$, we have
$
J_n\sim n^{-\beta}
$
because $x_n$ and $\widetilde x_n$ are at least $\frac 14$-distance apart.
The sequence $\{\frac{1}{2\pi}n^{\alpha}\}$ is uniformly distributed mod $1$ (see, e.g.,
Section I.3 in \cite{KN74}). Therefore, there is $N_0$ so that  for every $N>N_0$, 
we will have a bound
\[
|\sin n^\alpha|>0.01
\]
for at least $\frac N2$ integer numbers $n\in [N,2N]$. Therefore, since $\alpha\in (0,1)$,
\[
\sum_{n=N}^{2N} \Bigl||\sin n^\alpha| \cdot J_n - O( n^{\alpha-1}J_n)\Bigr|^2\gtrsim N^{1-2\beta}.
\]
Since $\beta\le \frac 12$, $\lim_{N\to\infty}N^{1-2\beta}$ is either equal to $1$ or is infinite and we have a contradiction with \eqref{s_a20}. Therefore, \eqref{help} fails.

\smallskip

\noindent {\bf Case \boldmath$\alpha \in (0,1)$, \boldmath$\beta > \frac 12$.} We have
$$
\sum_{n \ge 1}\int_{I_n}g_n^2\,dx 
\le
4\sum_{n \ge 1}\int_{I_n}\left(\int_{n}^{x}\frac{d\tau}{\tau^\beta}\,\right)^2\,dx \le
2\sum_{n \ge 1}\left(\int_{I_n}\frac{d\tau}{\tau^\beta}\,\right)^2\,dx \le 
8\sum_{n \ge 1} \frac{1}{n^{2\beta}} < 
\infty.
$$
Hence, \eqref{help} is true for this range of parameters by Proposition \ref{prop4}.

\smallskip

\noindent {\bf Case \boldmath$\alpha \le 0$, \boldmath$\beta-\alpha > \frac 12$.} For this range of parameters, we have
$$
\sum_{n \ge 1}\int_{I_n}g_n^2\,dx  \lesssim   \sum_{n \ge 1}\frac{1}{n^{2(\beta-\alpha)}}
$$
due to the fact that $\sin y$ is comparable to $y$ when $y \in [0,1]$. Using Proposition \ref{prop4}, we conclude that \eqref{help} holds if  $\beta-\alpha > \frac 12$.

\smallskip

\noindent {\bf Case \boldmath$\alpha \le 0$, \boldmath$\beta-\alpha \le \frac 12$.} We claim that \eqref{help} fails in that situation. Arguing as in the case when $\alpha \in (0,1)$ and $\beta \le \frac 12$, we assume that \eqref{help} is true and obtain
\begin{equation}\label{eq66}
\sum_{n \ge 1}\left|g_n(\widetilde x_n) - g_n(x_n)\right|^2 < \infty, 
\end{equation}
for some points $x_n, \widetilde x_n \in I_n$ such that $|x - \widetilde x_n| > 1/4$. Then, 
$$
\sum_{n \ge 0}n^{2(\alpha-\beta)} \lesssim \sum_{n \ge 0}\left(\int_{x_n}^{\widetilde x_n}\tau^{\alpha-\beta}\,d\tau\right)^2 \le   \sum_{n\ge 0}\left|\int_{x_n}^{\widetilde x_n}\frac{\sin \tau^\alpha}{\tau^\beta}\,d\tau\right|^2
< \infty 
$$
which leads to a contradiction. \qed

\medskip

We end this section with one more example of an oscillating potential $Q \in \szd$. As we mentioned in Section \ref{s1point2} of the Introduction, a simplest ``physical'' interpretation of the Dirac operator $\Dd_{Q}$ is that 
it describes the one-dimensional  particle of unit mass moving in a  field defined by a potential $Q$ of the form
\begin{equation}\label{eq102}
Q = \begin{pmatrix}q_1 & q_2\\ q_2  & -q_1\end{pmatrix}, 
\qquad q_1=\cos\left(2\int_0^\tau q\,ds\right),
\qquad q_2=-\sin\left(2\int_0^\tau q\,ds\right),
\end{equation}
where $q\colon \R_+ \to \R$. It is natural to ask if there exist potentials $Q \in \szd$ that can be written that way.  Since $q_1^2 + q_2^2 = 1$ on $\R_+$, we cannot have $Q \in L^p(\R_+)$ for any $p<\infty$. In particular, the standard summability test (see the discussion after Theorem \ref{t3-i})
$$
Q \in \bigcup_{p \in [1,2]}L^p(\R_+) \Longrightarrow Q \in \szd
$$ 
does not give any example of such a potential $Q$. However, Proposition \ref{prop4} can be used to construct such examples. Indeed, let us consider
$$
q_1(\tau) = \cos\left(\frac{\pi}{2}(-1)^{[e^\tau]}\right) = 0, \qquad 
q_2(\tau) = -\sin\left(\frac{\pi}{2}(-1)^{[e^\tau]}\right), \qquad \tau \in \R_+, 
$$
where $[e^{\tau}]$ stands for the integer part of $e^{\tau}$. Due to high oscillation of $q_2$, Proposition \ref{prop4} applies and  the corresponding potential  $Q = \left(\begin{smallmatrix}q_1 & q_2\\ q_2  & -q_1\end{smallmatrix}\right) = \left(\begin{smallmatrix}0 & q_2\\ q_2  & 0 \end{smallmatrix}\right)$ belongs to $\szd$. 
Then, one can construct its small perturbation of the form \eqref{eq102}  which belongs to $\szd$ (the Theorem~\ref{t3-i} can be used for that approximation argument). 
\section{Appendix}
\addtocontents{toc}{\protect\setcounter{tocdepth}{0}}
 
In this Appendix, we collect a few auxiliary results and  prove some statements made in the main text.
 
\subsection{Proof of Proposition \ref{str1}}\label{app1} The modification of the  proofs in \cite{BD2019} yields the statement. Alternatively, one can argue as follows.  First, we claim that \eqref{s_a5} implies  \eqref{fix1}  for $\alpha_n=\lambda n$ with any $\lambda>0$. 
Indeed, if $m_\lambda$ and $\mu_\lambda$ denote Titchmarsh-Weyl function and spectral measure, respectively, of canonical system with Hamiltonian $\Hh(\lambda \tau)$, then $m_\lambda(z)=m_1(\lambda^{-1} z)$ as  follows from (1-5) in \cite{BD2019}. 
Now, it is enough to observe that  $\mu_\lambda\in \sz \Longleftrightarrow \mu_1\in \sz$.\smallskip

Second, we claim that, given intervals $I^-\subseteq I, |I|=1$, and $\eps\in (0,1]$, the following implication holds
\begin{equation}\label{fix2}
\det \int_I \Hh(\tau) d\tau- 1=\eps \Rightarrow\det \int_{I^-} \Hh(\tau) d\tau- |I^-|^2\lesssim \eps, 
\end{equation}
for every non-negative Hamiltonian $\Hh$ that satisfies $\det \Hh=1$ a.e. on $I$. Indeed, denote $A=\int_{I^-} \Hh d\tau$ and  $B=\int_{I^+} \Hh d\tau$ where $I=I^-\cup I^+$. Then, we get (see, e.g., (A-1) in \cite{BD2019}):
\[
\det A\ge |I^-|^2, \quad \det B\ge |I^+|^2
\]
and
\[
\det (A+B)= 1+\eps.
\]
Minkowski inequality for determinants yields
\[
\det(A+B)\ge (\sqrt{\det A}+\sqrt{\det B})^2.
\]
Denoting $\sqrt{\det A}=x$ and $\sqrt{\det B}=y$, we get
\[
x+y\le (1+\eps)^{\frac 12}, \quad|I^-|+|I^+|=1,\quad |I^-|\le x, \quad |I^+|\le y.
\]
That implies (draw the corresponding domains on the plane), that $|I^-|\le x\le |I^-|+C\eps$. Taking the square of the last bound yields the estimate on the right-hand side of \eqref{fix2}.\smallskip

Now, if $\det \Hh=1$ a.e., the statement in  \eqref{fix1} holds by combining these two claims. Indeed, suppose $\mu\in \szp$ and we are given sequence $\{\alpha_n\}$. Then, there is $\lambda$ such that every interval $[\alpha_n, \alpha_{n+2}]$ is inside one of the intervals $[\lambda l, \lambda (l+2)]$ or $[\lambda (l-1),\lambda (l+1)]$ for some $l$. Since the sum in \eqref{fix1} converges for $\{\alpha_l\}=\{\lambda l\}$, we can apply \eqref{fix2} (with dilated and translated interval $I$) to get condition in the left-hand side of \eqref{fix1} satisfied for $\{\alpha_n\}$. Conversely, if the sum in \eqref{fix1} converges for some $\{\alpha_n\}$, then there is suitable $\lambda$ such that each interval $[\lambda n, \lambda (n+2)]$ is covered by either $[\alpha_l, \alpha_{l+2}]$ or $[\alpha_{l-1}, \alpha_{l+1}]$ for some $l$. Thus, applying \eqref{fix2} again, we get that the sum in \eqref{fix1} converges with $\{\alpha_n\}=\{\lambda n\}$ and $\mu\in \szp$. The case of general $\Hh$ follows by making the change of variables in $\tau$ and using an approximation argument.\qed\medskip

\subsection{Free evolution for canonical systems and Dirac operators.}\label{app2} Recall that $\Dd_0=\Di_{\Hh_0}$. We now show that the free evolution for these operators is, in fact, equivalent to the shift on the real line and that relation is algebraic.  To this end, we work in terms of Dirac operator and perform two elementary unitary transformations
\[
\widetilde \Dd_0=-Z^{-1}\Dd_0 Z= \left(\begin{smallmatrix} -i\partial_\tau & 0\\ 0 & i\partial_\tau\end{smallmatrix}\right), \quad Z=\tfrac{1}{\sqrt 2} \left(\begin{smallmatrix} i & -1\\ 1 & -i\end{smallmatrix}\right).
\]
Here, $\partial_\tau$ stands for the differentiation operator. Operator $\widetilde \Dd_0$, taken with suitable boundary condition at $0$: $f_1(0)=if_2(0)$, is self-adjoint on
the same Hilbert space $\left(\begin{smallmatrix} f_1\\f_2\end{smallmatrix}\right)\in L^2(\R_+,\C^2)$. If one further maps 
\[
\left(\begin{smallmatrix} f_1\\ f_2\end{smallmatrix}\right)\mapsto g(x)=\begin{cases} f_1(x), &x>0,\\ if_2(-x), &x<0,\end{cases}
\]
then, $\Dd_0$ becomes unitary equivalent to $-i\partial_x$ on $L^2(\mathbb{R})$ with $e^{t\partial_x}g=g(x+t)$, which is the standard shift operator.
\medskip

\subsection{A formula for exponential type}\label{app3}  
\begin{Lem}If entire function $f$ has bounded type both in $\C_+$ and  $\C_-$, then its exponential type can be computed by the formula
\begin{equation}\label{eq21bis}
\type f = \limsup_{y \to +\infty}\frac{\log\max(|f(iy)|, |f(-iy)|)}{y}.
\end{equation}\label{appe1}
\end{Lem}
\beginpf Let us apply Theorem 2 in Lecture 16 of \cite{Levin}. It says that for every entire function $f$ of bounded type in $\C_+$ and $\C_-$ we have
\begin{align*}
\log |f(z)| &= \sigma_+ y  + o(|z|),  \qquad y \ge 0,\\
\log |f(z)| &= \sigma_- y  + o(|z|),  \qquad y \le 0,
\end{align*}
outside of a set of disks $\{z \in \C\colon |z - a_j| < r_j\}$ of finite view (the latter means that $\sum_j \frac{r_j}{|a_j|} < \infty$). Here $\sigma_{\pm}  \in \R$ and $y = \Im z$. Take $\eps > 0$ and denote $\sigma = \max(\sigma_{+}, -\sigma_-)$. By the maximum principle for subharmonic functions, we have
$$
\log |f(z)| \le (\sigma + \eps)|y| + o(|z|)
$$
everywhere in $\C$ as $z \to \infty$. Therefore, we have $\type f \le \sigma$. On the other hand, the set of disks of finite view cannot fill any half-axis, hence 
$$
\sigma_+ = \limsup_{y \to +\infty}\frac{\log|f(iy)|}{y}, \qquad -\sigma_- = \limsup_{y \to -\infty}\frac{\log|f(iy)|}{|y|},
$$ 
which proves the statement. \qed

\medskip

\subsection{Rotation matrices}\label{app4}  The following result in the linear algebra has been used in the main text.
\begin{Lem}\label{l3}
For every real $2\times 2$ matrix $A$ with non-negative determinant,   
there is a rotation matrix $\Sigma_\phi$ of the form
$$
\Sigma_\phi = \begin{pmatrix}\cos\phi & \sin\phi \\ -\sin\phi & \cos\phi\end{pmatrix}, \qquad \phi \in [0, 2\pi),
$$
such that $\Sigma_\phi A \ge 0$.
\end{Lem}
\beginpf 
This is immediate from the proof of the polar decomposition given in  \cite{gant}, p.\ 276.
\qed \bigskip

\subsection{Robinson's theorem}\label{app5} In the main text, we used the following variation of a result by Robinson \cite{Robinson}, which is based on ideas dating back to Ruelle's work \cite{Ruelle}.
\begin{Lem} Suppose $H$ is a Hilbert space, $\Di$ is a densely defined self-adjoint operator, and $P_{\Delta}$ denotes the orthogonal projector for $\Di$ relative to a set $\Delta\subseteq \R$. If $A$ is a bounded  operator on $H$ and $AP_{[-\Lambda,\Lambda]}$ is compact for some $\Lambda\ge 0$, then
\[
\lim_{\tb\to +\infty}\frac 1\tb\int_0^\tb \|Ae^{it\Di}P_{[-\Lambda,\Lambda]}\psi\|^2dt=\sum_{j}\|AP_{E_j}P_{[-\Lambda,\Lambda]}\psi\|^2
\]
for every $\psi\in H$, where $P_{ E_j}$ denotes the orthogonal projection on the eigenspace that corresponds to eigenvalue $E_j$ and the sum is over all eigenvalues $\{E_j\}$ of $\Di$.\label{lemo1}
\end{Lem}
\beginpf The proof is an application of Theorem 2 in \cite{Robinson} to operator $\Di P_{[-\Lambda,\Lambda]}$ with the perturbation taken as $AP_{[-\Lambda,\Lambda]}$ where  both $\Di P_{[-\Lambda,\Lambda]}$ and $AP_{[-\Lambda,\Lambda]}$  are considered as operators acting on $H$. \qed

\bibliographystyle{plain} 
\bibliography{bibfile}
\end{document}